\documentclass[11pt,letterpaper]{amsart}
\usepackage[utf8]{inputenc}
\usepackage{array,longtable,booktabs}
\usepackage{amssymb,graphicx,mathrsfs,mathtools}
\usepackage{tikz-cd}
\usepackage[abbrev,nobysame,alphabetic]{amsrefs}
\usepackage[hidelinks]{hyperref}

\theoremstyle{plain}
\newtheorem{theorem}{Theorem}[section]
\newtheorem{proposition}[theorem]{Proposition}
\newtheorem{lemma}[theorem]{Lemma}
\newtheorem{corollary}[theorem]{Corollary}
\theoremstyle{definition}
\newtheorem{definition}[theorem]{Definition}
\newtheorem{notation}[theorem]{Notation}
\newtheorem{remark}[theorem]{Remark}
\theoremstyle{remark}
\newtheorem{example}[theorem]{Example}
\newtheorem{claim}[theorem]{Claim}
\numberwithin{equation}{section}

\DeclareMathOperator{\Hom}{Hom}
\DeclareMathOperator{\RHom}{RHom}
\DeclareMathOperator{\Ext}{Ext}
\DeclareMathOperator{\Spec}{Spec}
\DeclareMathOperator{\Pic}{Pic}
\DeclareMathOperator{\codim}{codim}
\DeclareMathOperator{\rk}{rk}
\DeclareMathOperator{\Sym}{Sym}
\DeclareMathOperator{\Bl}{Bl}
\DeclareMathOperator{\Ker}{Ker}
\DeclareMathOperator{\cD}{\mathcal{D}}
\DeclareMathOperator{\cM}{\mathcal{M}}
\DeclareMathOperator{\cE}{\mathcal{E}}
\DeclareMathOperator{\cH}{\mathcal{H}}
\DeclareMathOperator{\cF}{\mathcal{F}}
\DeclareMathOperator{\cO}{\mathcal{O}}
\DeclareMathOperator{\cP}{\mathcal{P}}
\DeclareMathOperator{\bP}{{\mathbb{P}}}
\DeclareMathOperator{\bZ}{{\mathbb{Z}}}
\DeclareMathOperator{\bLambda}{{\boldsymbol\Lambda}}
\DeclareMathOperator{\Supp}{Supp}
\DeclareMathOperator{\depth}{depth}
\newcommand{\scrL}{\mathscr L}
\newcommand{\StacksTag}[1]{\href{https://stacks.math.columbia.edu/tag/#1}{Stacks Project, Tag~#1}}
\newcommand{\ft}{{\lfloor t\rfloor}}
\newcommand{\fr}{{\lfloor r\rfloor}}
\def\arrow{\mathop{\longrightarrow}\limits}

\title{Braid and Phantom}
\author{Jenia Tevelev}

\begin{document}

\begin{abstract}
Let $N$ be the moduli space of stable rank-$2$ vector bundles on a smooth
projective curve of genus $g\ge2$ with fixed determinant of odd degree. In
\cite{TT}, we constructed a semiorthogonal decomposition of $D^b(N)$ into
bounded derived categories of symmetric powers of the curve and a possible
phantom component. In this paper, we use weaving patterns to eliminate the
phantom and complete the proof of the semiorthogonal decomposition conjectured by
Narasimhan and, independently, by Belmans--Galkin--Mukhopadhyay.
 \end{abstract}

\maketitle

\section{Introduction}

Let $C$ be a smooth projective curve of genus $g\ge2$, and let $N$ be the
moduli space of stable rank-$2$ vector bundles on $C$ with a fixed
determinant line bundle of odd degree.  Let $\theta$ be the ample generator
of $\Pic(N)\simeq\mathbb Z$, and let $\cE$ be the Poincar\'e bundle on
$C\times N$, normalized by
$\det(\cE|_{\{p\}\times N})\simeq\theta$ for every $p\in C$.
The tensor bundle $\cE^{\boxtimes k}$ on $\Sym^kC\times N$ has
rank $2^k$ \cite[Definition~2.3]{TT}; the corresponding
Fourier--Mukai functor
\[
 \cP_{\cE^{\boxtimes k}}:D^b(\Sym^kC)\longrightarrow D^b(N)
\]
is fully faithful for $0\le k\le g-1$ \cite[Theorem~1.1]{TT}.

For smooth projective varieties $X$ and $Y$ and a fully faithful
Fourier--Mukai functor
$
 \mathcal P_{\mathcal A}:D^b(X)\longrightarrow D^b(Y)$ with kernel
$\mathcal A\in D^b(X\times Y)$,
we write $\langle\mathcal A\rangle$ for its  essential image and use similar notation for semi-orthogonal decompositions.
When $X$ is a point, the notation records an exceptional object or collection;
for example, $D^b(\mathbb P^n)=\langle\mathcal O,\ldots,\mathcal O(n)\rangle$.

\begin{theorem}\label{MainTheorem}
$D^b(N)$ has a semiorthogonal decomposition into blocks equivalent to
$D^b(\operatorname{Sym}^kC)$: two blocks for $0\le k\le g-2$ and one
block for $k=g-1$. These blocks are
embedded by the following Fourier--Mukai functors,
$$
\Bigl\langle
\bigl\langle 
{\theta}^{1-g+k}\otimes\cE^{\boxtimes g-2-2k}
\bigr\rangle_{0\le k\le \lfloor{g-2\over2}\rfloor},\quad
\bigl\langle 
{\theta}^{2-g+k}\otimes\cE^{\boxtimes g-3-2k}
\bigr\rangle_{0\le k\le  \lfloor{g-3\over2}\rfloor},\quad{}$$
$${}\quad\bigl\langle 
{\theta}^{2-g+k}\otimes\cE^{\boxtimes g-2-2k}
\bigr\rangle_{0\le k\le \lfloor{g-2\over2}\rfloor},
\quad \bigl\langle 
{\theta}^{2-g+k}\otimes\cE^{\boxtimes g-1-2k}
\bigr\rangle_{0\le k\le  \lfloor{g-1\over2}\rfloor}
\Bigr\rangle,
$$ 
arranged into four mega-blocks.
Within the mega-blocks, the 
blocks are arranged in decreasing order of~$k$.
\end{theorem}

We refer to \cite{TT} for the history of the problem and the proof of semi-orthogonality. 
For the alternative decomposition
used in the paper \cite{TT}, see Theorem~\ref{ComparingSODs}.
Tensoring bundles by a line bundle identifies the fixed-odd-determinant
moduli spaces, so we  choose the determinant $\Lambda$ to have
degree $2g-1$. As in \cite{TT}, let $M$ be the Fano moduli space of pairs
$(F,s)$, where $F$ is a stable bundle with determinant $\Lambda$ and
$0\ne s\in H^0(C,F)$.

We~play the two-ray game. The variety $M$ admits the forgetful birational morphism
$\zeta:M\to N$, $(F,s)\mapsto F$, and a birational map
$\psi:M\dashrightarrow\mathbb P^{3g-3}$, which can be described as follows.  On the dense open locus of $M$ on which
the section $s$ is nowhere vanishing, the section defines an exact sequence
\[
0\longrightarrow\mathcal O_C\xrightarrow{s}F\longrightarrow\Lambda
\longrightarrow0.
\]
The map $\psi$ sends $(F,s)$ to the corresponding extension class in
$\mathbb P\Ext^1(\Lambda,\mathcal O_C)$.
Furthermore, 
the map $\psi$ factors as a sequence of standard flips followed by a divisorial contraction \cite{thaddeus}:
\begin{equation}\label{THDiagram}
M=M_{g-1}\dashrightarrow M_{g-2}\dashrightarrow\ldots\dashrightarrow M_1\to M_0=\bP^{3g-3}.
\end{equation}
At each step $M_k\dashrightarrow M_{k-1}$, one projective bundle over $\Sym^kC$ is flipped into another of smaller rank.
By \cite[Proposition~3.18]{TT}, this wall-crossing construction gives a
semiorthogonal decomposition of $D^b(M)$ with $3g-2-3k$ blocks equivalent
to $D^b(\Sym^kC)$ for each $0\le k\le g-1$. 
 
The blocks with $k=0$ are the line-bundle blocks inherited from the
Beilinson collection on $\mathbb P^{3g-3}$.  For $k>0$, the blocks are
introduced by the successive wall crossings and their Fourier--Mukai kernels
have torsion supports.
We will mutate this decomposition into one in which the desired blocks
belong to $L\zeta^*D^b(N)$ and all other blocks belong to
\[
 \Ker(R\zeta_*)=\bigl(L\zeta^*D^b(N)\bigr)^\perp.
\]
Here the equality follows from adjunction $L\zeta^*\dashv R\zeta_*$, and
$L\zeta^*$ is fully faithful because $R\zeta_*\mathcal O_M\simeq\mathcal O_N$.
The projection onto $L\zeta^*D^b(N)$ then kills the remaining blocks and
shows that the desired blocks generate $D^b(N)$.

\begin{figure}[htbp]
\centering
\includegraphics[width=\textwidth,height=0.94\textheight,keepaspectratio]{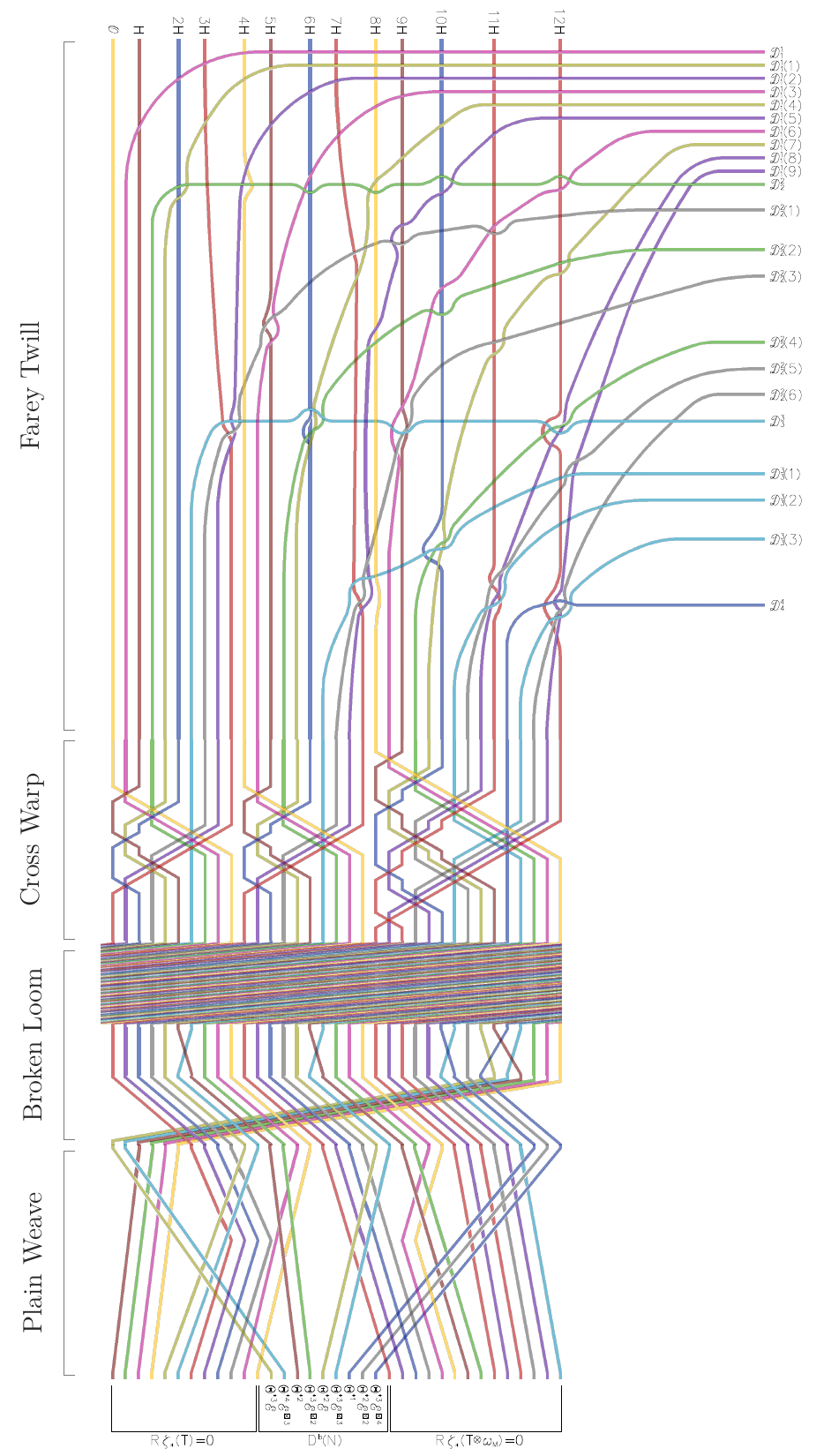}
\caption{All weaving patterns in genus $5$.}\label{genus5total}
\end{figure}

To construct the mutation, we combine vanishing theorems from \cite{TT} with an analysis of ``weaving patterns'' related to various geometric constructions. The result is illustrated in Figure~\ref{genus5total} (for genus 5). We use the term ``weaving'' as opposed to ``braiding'' to emphasize the importance of invisible wefts, which break  the mutation into a sequence of standard steps. 

\subsection*{Notation and conventions}

All schemes are over $\mathbb C$.

For an admissible subcategory
$\mathcal A\subset D^b(Y)$, we use the notation
\[
 \mathcal A^\perp=\{T:\RHom(\mathcal A,T)=0\},
 \qquad
 {}^\perp\mathcal A=\{T:\RHom(T,\mathcal A)=0\}.
\]
Thus left mutation $\mathbf L_{\mathcal A}$  through $\mathcal A$ takes values in
$\mathcal A^\perp$, and right mutation $\mathbf R_{\mathcal A}$ takes values in
${}^\perp\mathcal A$. Sequences of mutations are arranged into weaving diagrams.
Every weaving diagram is read from top to bottom.  A strand represents an admissible
subcategory; labels at the top and bottom record its Fourier--Mukai kernel
before and after the indicated mutations.  A transverse crossing means that
the two corresponding subcategories are completely orthogonal, so their
order may be interchanged without changing either embedding.  At a bend, the
strand passing over the bend keeps the same embedding, whereas the strand
passing under it is replaced by the displayed left or right mutation.  
The figures record the order of operations;
the semiorthogonality and mutation triangles are proved in the text.

Pullbacks, tensor products, and pushforwards between derived categories
are derived.  We omit $L$ or $R$ when no confusion can arise.  

\begin{longtable}{@{}p{0.22\textwidth}p{0.73\textwidth}@{}}
\toprule
\textbf{Symbol} & \textbf{Meaning and reference}\\
\midrule
\endfirsthead
\toprule
\textbf{Symbol} & \textbf{Meaning and reference}\\
\midrule
\endhead
$C$, $g$, $\Lambda$ & The smooth projective curve, its genus, and the
fixed determinant line bundle of degree $2g-1$ (Section~1).\\[2pt]
$N$, $\theta$, $\mathcal E$ & The moduli space of stable bundles, 
ample generator of $\Pic N$, and the normalized Poincar\'e bundle
(Section~1).\\[2pt]
$M_i$, $M$ & The Thaddeus spaces of stable pairs for determinant $\Lambda$ and
$M=M_{g-1}$ (Notation~\ref{notation:stable-pairs}).\\[2pt]
$M_i(\Lambda')$, $M_i(d)$ & The corresponding Thaddeus spaces for other
determinants $\Lambda'$; $d=\deg\Lambda'$ (Notation~\ref{notation:other-determinants}).\\[2pt]
$\psi$, $\zeta$ & The birational map $M\dashrightarrow M_0$ and the
forgetful morphism $M\to N$ (Section~1).\\[2pt]
$Z$ & The exceptional divisor
of  $\zeta$ (Notation~\ref{notation:universal-pair}).\\[2pt]
$(\mathcal F,\Sigma)$, $\boldsymbol\Lambda$ & The universal stable pair on $C\times M_i$
and the line bundle on $M_i$ with
$\det\mathcal F=\Lambda\boxtimes\boldsymbol\Lambda$
(Notation~\ref{notation:universal-pair}).  
\\[2pt]
$H$, $E$, $\mathcal O(m,n)$ & The divisor classes and line bundles of
Notation~\ref{notation:stable-pairs}.\\[2pt]
$D_i^k$, $\mathcal D_i^k$ & The  variety $\{(D,F,s)\,:\,s|_D=0\}$ and its structure
sheaf (Notation~\ref{notation:stable-pairs}).  For $i=g-1$ the subscript
is omitted.\\[2pt]
$\mathcal P_{\mathcal K}$, $\langle\mathcal K\rangle$ & The
Fourier--Mukai functor with kernel $\mathcal K$ and, when it is fully
faithful, its essential image (Section~1).\\[2pt]
$\tau_k$, $\operatorname{Alt}_k$ & The  quotient $\tau_k:C^k\times M_i\to\Sym^kC\times M_i$ and alternating
direct image  $(\tau_{k*}(-\otimes\operatorname{sgn}_k))^{S_k}$ (Notation~\ref{notation:ordered-quotient}).\\[2pt]
$\mathcal F^{\boxtimes k}$,
$\overline{\mathcal F}^{\boxtimes k}$ & The ordinary and alternating
tensor bundles on $\Sym^kC\times M_i$ (Notation~\ref{notation:ordered-quotient}).  A hat denotes
the corresponding object on the ordered product $C^k\times M_i$.\\[2pt]

$\Delta_k$, $\mathscr L_k$ & The divisor of nonreduced divisors on
$\Sym^kC$ and its sign line bundle
$\mathscr L_k=\mathcal O_{\Sym^kC}(-\Delta_k/2)$
(Lemma~\ref{amazinggrace}).\\[2pt]
$\mathcal F^\bullet$, $\mathcal F^{\bullet\boxtimes k}$ & The $2$-term
evaluation complex and its alternating tensor product on $\Sym^kC\times M_i$
(Definitions~\ref{def:universal-two-term} and~\ref{,aENFVkejhfv}).\\[2pt]
$\widehat D_i^k$, $\widehat{\mathcal D}_i^k$ & The pullback of the
incidence variety and its structure sheaf to $C^k\times M_i$
(Definition~\ref{def:ordered-incidence}).\\[2pt]
$\mathcal G_{[0,i)}$, $\iota_i$ & The window and its induced embedding
in~\eqref{windowsembedding}.  
\\
\bottomrule
\end{longtable}
The notation used only in the proof of
Theorem~\ref{sdfgnsfghsftjt}, including $Y_{r,l}$,
$Z_{r,l}$, $\Gamma_{r,l}$, and the local rings in
Proposition~\ref{prop:split-comparison}, is introduced in Section~\ref{asrgarsgargaerhaehae}.

\subsection*{Acknowledgments}
This project would have been impossible without our collaboration
with Sebasti\'an Torres \cite{TT}.
I am grateful to An\'ibal Aravena and Alex Perry
for helpful discussions, and to Elias Sink and the anonymous referees
for pointing out several mistakes in earlier versions of the paper and
suggesting how to correct them.
Most of the paper was written during a visit to Imperial College
in January 2023, and I thank Paolo Cascini, Alessio Corti, and
Yank{\i} Lekili for providing an inspiring environment.
This research was partially supported by NSF grants DMS-2101726
and DMS-2401387.
The graphics were created by \url{http://www.plainformstudio.com}.

\tableofcontents

\section{Farey Twill}\label{TwillSection}

\begin{figure}[htbp]
\includegraphics[width=\textwidth]{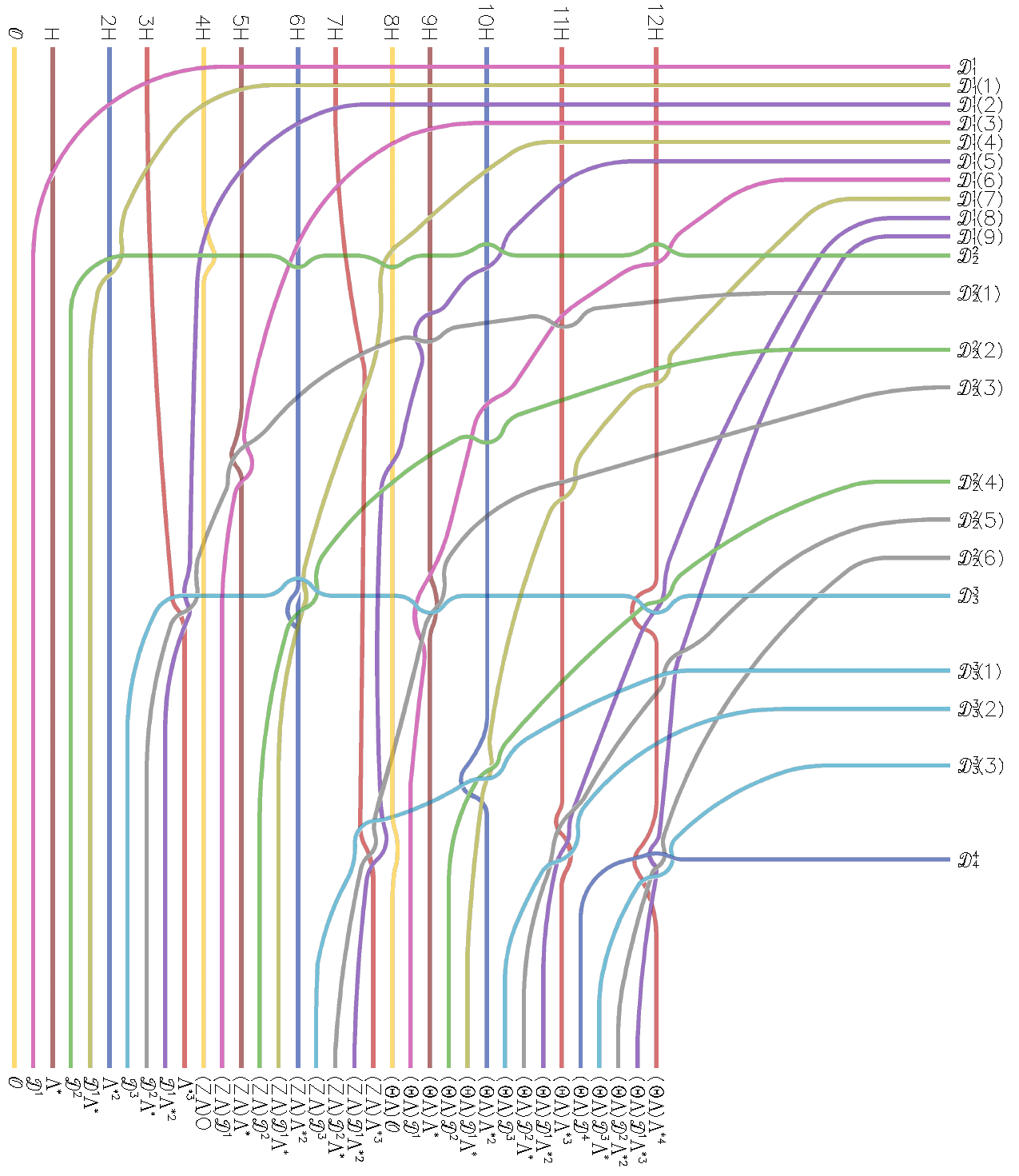}
\caption{Farey Twill in genus $5$}\label{genus5twill}
\end{figure}

The Farey Twill, illustrated in Figure~\ref{genus5twill}, reorganizes the
semiorthogonal decomposition of $D^b(M)$ constructed in
\cite[Proposition~3.18]{TT}.  Its blocks are equivalent to
$D^b(\Sym^kC)$ for $0\le k\le g-1$.  The $k=0$ blocks are given by line bundles,
whereas the Fourier-Mukai kernels with $k\ge1$ are supported on the incidence loci
created by the successive wall crossings.

The goal of this section is to prove Theorem~\ref{TwillTheorem} below,
except for  several lemmas,
which will be proved in Section~\ref{CrossWarpSection}.
We first introduce some notation, mostly from \cite{TT} and \cite{thaddeus},  that will be used throughout the paper.

\begin{notation}\label{notation:universal-pair}
Let $(\cF,\Sigma)$ be the universal stable pair on $C\times M$. Then 
$\cF=\cE(-Z)$,
$\det\cF=\Lambda\boxtimes\bLambda$,
and $\theta=\bLambda(2Z)$, 
where $Z$ is the exceptional divisor of the  contraction $\zeta:\,M\to N$.
 Let $\cD^k$ be the structure sheaf of 
 a smooth subscheme  
 $D^k=\{(D,F,s):\,s|_D=0\} \subset \Sym^kC\times M$
 of codimension $2k$,
 where we view $D\in \Sym^kC$ as a closed subscheme of $C$. 
\end{notation}

\begin{theorem}\label{TwillTheorem}
The Fourier--Mukai functor $\cP_{\cD^k}:\,D^b(\Sym^kC)\to D^b(M)$ is fully faithful for $k\le g-1$ and 
$D^b(M)$ has a   semiorthogonal decomposition
into admissible subcategories arranged into three mega-blocks, as follows:
$$
\Bigl\langle
\bigl\langle {\bLambda^*}^j\otimes\cD^k\bigr\rangle_{ j+k\le g-2\atop j,k\ge0},
\ \bigl\langle (Z\bLambda){\bLambda^*}^{j}\otimes\cD^k\bigr\rangle_{ j+k\le g-2\atop j,k\ge0},
\ \bigl\langle (\theta\bLambda){\bLambda^*}^{j}\otimes\cD^k\bigr\rangle_{ j+k\le g-1\atop j,k\ge0}
\Bigr\rangle.
$$ 
Within each of the three mega-blocks, the blocks are ordered first by increasing $j+k$ and, for fixed $j+k$, by increasing $j$.
\end{theorem}

In genus $5$, the blocks of Theorem~\ref{TwillTheorem} appear at the bottom of Figure~\ref{genus5twill}.

 \begin{notation}\label{notation:stable-pairs}
Let $E$ be the exceptional divisor of the birational morphism $M_1\to M_0=\bP^{3g-3}$
and let $H$ be the pullback of the hyperplane divisor.
The varieties $M_1,\ldots,M_{g-1}=M$ appearing in the sequence of flips \eqref{THDiagram} are isomorphic in codimension~$1$.
This allows us to use the same notation for their line bundles, including an important line bundle
\cite{thaddeus}
\[
\cO(m,n):=\cO_{M_i}\bigl((m+n)H-nE\bigr).
\]
We use the moduli interpretation of $M_i$ throughout.  If $(F,s)\in M_i$
and $D=Z(s)$ is the zero divisor of $s$, the saturation of
$\mathcal O_C\xrightarrow{s}F$ is the line subbundle $\mathcal O_C(D)$.
The stability inequality gives $\deg D\le i$;
see~\cite[Definition~1.1 and Section~3]{thaddeus}.
For $0\le k\le i$, put
\[
 D_i^k=\{(D,F,s)\in\operatorname{Sym}^kC\times M_i:s|_D=0\},
 \qquad \mathcal D_i^k=\mathcal O_{D_i^k}.
\]
The notation $s|_D=0$ includes vanishing to the prescribed multiplicities
when $D$ is nonreduced.  For $i=g-1$ this agrees with $D^k$ and
$\mathcal D^k$ above.
\end{notation}

The following lemma, as well as several others, will be proved in Section~\ref{CrossWarpSection}.
\begin{lemma}\label{asgasrhare}
The Fourier--Mukai functor $\cP_{\cD_i^k}:\,D^b(\Sym^kC)\to D^b(M_i)$ is fully faithful for $0\le k\le i\le g-1$.
\end{lemma}

\begin{definition}
We introduce admissible (by Lemma~\ref{asgasrhare}) subcategories 
$$\langle\cD^{k,s}_t\rangle\subset D^b(M_\ft)$$ for 
integer parameters
$0\le k\le g-1$, $0\le s\le 3g-3-3k$, and a real parameter $k<t<g$. Namely, $\langle\cD^{k,s}_t\rangle$ is  
the image of $D^b(\Sym^k C)$
under
the Fourier--Mukai functor with the  kernel $\cD^{k,s}_t=\cD^k_\ft\otimes L^{k,s}_t$, where
\begin{equation}\label{wonderingblocks}
L^{k,s}_t=\begin{cases}
\cO\left(
\left\lfloor{s\over t-k}\right\rfloor, s+\left\lfloor{s\over t-k}\right\rfloor(k-1)\right)&\hbox{\rm if}\quad \ft>k,\cr
\cO(s,sk) &\hbox{\rm if}\quad \ft=k.\cr
\end{cases}
\end{equation}
Here $\lfloor x\rfloor$ denotes the greatest integer not exceeding $x$.
\end{definition}

\begin{figure}[htbp]
\includegraphics[width=\textwidth]{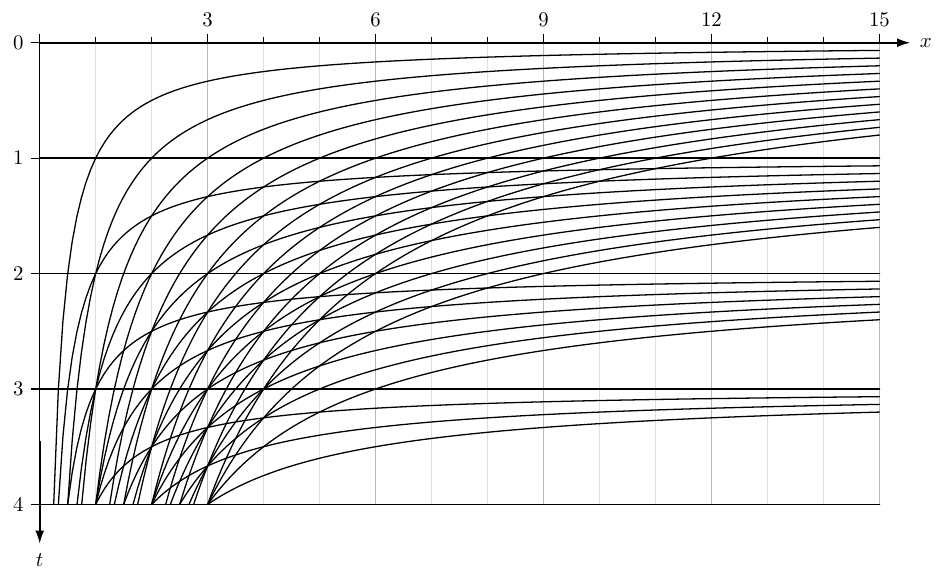}
\caption{Genuine paths of blocks 
in the Farey Twill (genus $5$)}\label{trajectories}
\end{figure}

We interpret the variable $t$ as time. In the Farey Twill, named after the Farey fractions, 
the admissible subcategory $\langle\cD^{k,s}_t\rangle$ 
``moves'' in the $(x,t)$-plane along the trajectory 
\begin{equation}\label{skhjgb,sHEG}
x_{k,s}(t)={s\over t-k},
\end{equation}
with the $t$-axis pointing down and the $x$-axis  to the right.

For the construction of Theorem~\ref{TwillTheorem}, we use these
trajectories only up to level $t=g-1$, followed by the terminal
chain mutation at that level. At each stage, the blocks
$\langle\cD^{k,s}_t\rangle\subset D^b(M_{\ft})$ are ordered by the value
of $x_{k,s}(t)$.
The mutations at crossings in this range will be justified below.
The paths of the blocks with $s=0$ are modified 
to allow them to participate in the mutations occurring at integer values of $t$.
At an integer boundary $t=k$, we use only the block with $s=0$ and set
$
\cD_k^{k,0}:=\cD_k^k
$.
Its~modified trajectory contains the horizontal segment at level $t=k$ and
then the short vertical segment used to order the blocks
$\cD_k^{k,0}$.  
The paths \eqref{skhjgb,sHEG} (and modified paths of blocks with $s=0$) are plotted in Figure~\ref{trajectories} in genus~$5$.


\begin{example} We analyze the Farey Twill for small values of~$t$, see Figure~\ref{twillstarts}. 
In~this and other illustrations, we apply an $(x,t)$-plane transformation to better visualize 
intersections of trajectories. For example, paths of the line bundles $H,2H,3H,\ldots$
should come from infinity  (see the top of Figure~\ref{trajectories}),
but in Figures~\ref{genus5total}, \ref{genus5twill}, and \ref{twillstarts} we draw them as parallel vertical lines.

\begin{figure}[htbp]
\includegraphics[width=\textwidth]{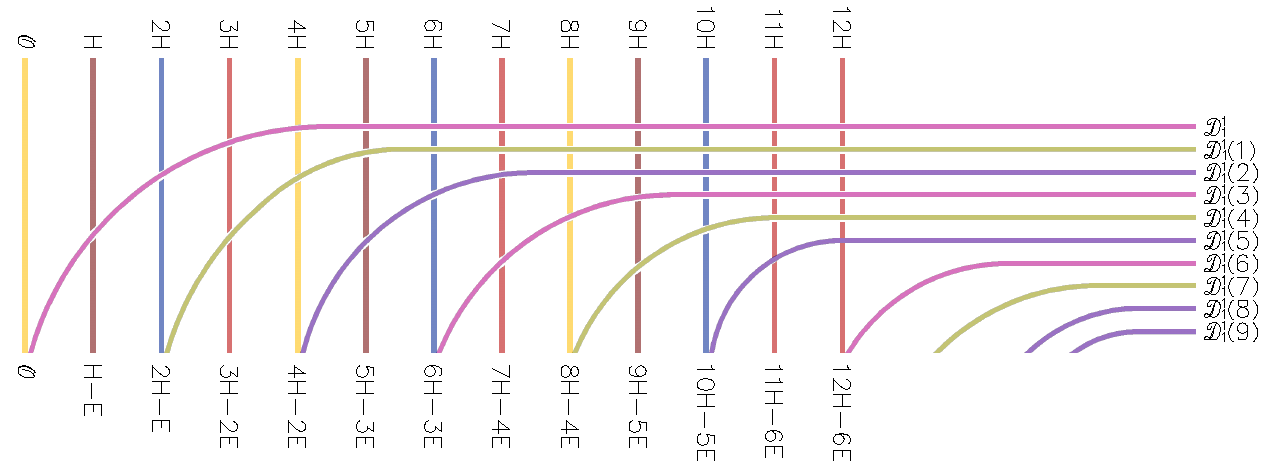}
\caption{Starting the Farey Twill in genus $5$}\label{twillstarts}
\end{figure}

If $\boxed{0<t<1}$,
 $\cD^{0,s}_t=\cO(sH)$ and
we have the Beilinson decomposition:
$$D^b(M_0)=D^b(\bP^{3g-3})=\langle\cO,H,\ldots,(3g-3)H\rangle.$$

At $\boxed{t=1}$, the same blocks $\cD^{0,s}_1=\cO(s,0)=\cO(sH)$ are regarded as pullbacks to $M_1$.
This gives an admissible subcategory 
$$\pi_1^*D^b(\bP^{3g-3})=\langle\cO,H,\ldots,(3g-3)H\rangle\subset D^b(M_1),$$
where $ \pi_1:\,M_1=\operatorname{Bl}_C(M_0)\longrightarrow M_0$ is the blow-up.

For the remainder of this example, assume $g\ge3$.
Let $\boxed{1<t<2}$. By \cite{thaddeus}, $M_1\simeq\Bl_C\bP^{3g-3}$
(where the curve $C$ is embedded in $\bP^{3g-3}$ by the linear system $|K_C+\Lambda|$).
Let 
\[
D_1^1=\bigl\{(p,(F,s))\in C\times M_1\mid s(p)=0\bigr\}.
\]
The second projection identifies $D_1^1$ with the exceptional divisor
$E\subset M_1$, and the first projection exhibits it as a projective bundle
over $C$.

The line bundles $\cO(s,s)$ restrict to $\cO(s)$ on the fibers $\bP^{3g-5}$ of the projective bundle.
So the sheaves $\cD^{1,s}_t=\cD^1_1(s,s)$ are the Fourier--Mukai kernels that appear in the 
Orlov decomposition \cite{orlov}
\begin{equation}\label{aasfbasfbasn}
D^b(M_1)=D^b(\Bl_C\bP^{3g-3})=\langle\cO,H,\ldots,(3g-3)H,\cD^{1,0}_t,\ldots,\cD^{1,3g-6}_t\rangle.
\end{equation}
The paths \eqref{skhjgb,sHEG} of the new blocks $\cD^{1,s}_t$ come from infinity, corresponding to the fact that these blocks appear on the right in \eqref{aasfbasfbasn}.

The sheaves $\cD^{1,s}_t$ do not change for $1<t<2$, but the line bundles in \eqref{aasfbasfbasn} undergo mutations as $t$ increases from $1$ to $2$.
Note that
$\cD^{0,s}_{1+\epsilon}=\cO(s-1,1)=\cO(sH-E)$ for $1\le s\le3g-3$ and $0<\epsilon\ll1$, so 
these  blocks have already changed compared to $t=1$. The corresponding mutations are encoded in the 
convention that paths of blocks with $s=0$ have to be modified to give the  asymptotes of the hyperbolas~\eqref{skhjgb,sHEG}.
Namely, we add the horizontal line $t=1$ to the trajectory of the block $\cD^{1,0}_t$ and
get a sequence of mutations on level $t=1$ induced by the standard short exact sequences (for $s>0$), 
\begin{equation}\label{wrgarhaerh}
0\to\cO(sH-E)\to\cO(sH)\to\cO_E(sH)\to0.
\end{equation}
This gives a semiorthogonal decomposition of $D^b(M_1)$ at level $t=1+\epsilon$:
$$\langle\cO,\cD^{1,0}_t,H-E,2H-E,\ldots,(3g-3)H-E,\cD^{1,1}_t,\ldots,\cD^{1,3g-6}_t\rangle.
$$
Note that we keep the block $\cD^{1,0}_t$ to the right of $\cO$. Formally speaking, we modify the vertical part of its path as well (make it $x=\epsilon$
for $0<\epsilon\ll1$).

The next trajectory to cross other trajectories in the Farey Twill is that of $\cD^{1,1}_t$, followed by those of $\cD^{1,2}_t$, and so on.
Mutations are given by the line bundle twists of \eqref{wrgarhaerh}, eventually
giving a semi-orthogonal decomposition of $D^b(M_1)$  at level $t=2-\epsilon$:
\begin{equation}\label{aEGSGS}
\langle\cO,\cD^{1,0}_t,H-E,2H-E,\cD^{1,1}_t,3H-2E,4H-2E,\cD^{1,2}_t,\ldots,\cD^{1,3g-6}_t\rangle.
\end{equation}
Many blocks $\cD^{1,s}_t$ remain at the end of the decomposition.

The advantage of~\eqref{aEGSGS} over~\eqref{aasfbasfbasn}
is its compatibility with the window embedding
$
 \iota:D^b(M_1)\simeq\mathcal G_{[0,2)}
 \subset D^b(\mathcal M)\longrightarrow D^b(M_2)
$
of~\cite[Proposition~3.18]{TT}. Here $\mathcal M$ is the common
quotient stack for the wall crossing.

The window consists of complexes
whose derived restriction to the center of the unstable stratum has
cohomology of weights $0,1$; see~\cite[Theorem~3.21]{TT} and
\eqref{windowsembedding} below.
The line bundles $\cO,H-E,2H-E,\ldots$ appearing
in \eqref{aEGSGS}
have weights $0,1,0,1,0,\ldots$ (we refer to \cite[Section~3]{TT} for the calculation of weights of all standard vector bundles), and so the windows
embedding $\iota$ maps them to the same line bundles on $M_2$. 
Furthermore, $\iota$ takes the decomposition  \eqref{aEGSGS}
into the following admissible subcategory:

$\boxed{t=2}$,
$\langle\cO,\cD^{1,0}_2,H-E,2H-E,\cD^{1,1}_2,3H-2E,4H-2E,\cD^{1,2}_2,\ldots,\cD^{1,3g-6}_2\rangle$.
At this level, complete this admissible subcategory of $D^b(M_2)$
to a semiorthogonal decomposition by adding the blocks
$\cD^{2,0}_{2+\epsilon},\ldots,\cD^{2,3g-9}_{2+\epsilon}$ at the end.
If $g\ge4$, the Farey Twill may then be continued through
$\boxed{2<t<3}$ using the mutations encoded by intersections of
trajectories. 
\end{example}

To realize this program in general, we consider the ``windows'' embedding 
from \cite[Proposition 3.18]{TT},
\begin{equation}\label{windowsembedding}
\iota_i:\,D^b(M_{i-1})\cong\mathcal G_{[0,i)}\subset D^b(\cM)\mathop{\to}\limits^{r} D^b(M_i),
\end{equation}
where $\cM$ is an appropriate quotient stack that contains  $M_{i-1}$ and $M_i$ as open substacks, $r$
is the restriction, and  $\mathcal G_{[0,i)}\subset D^b(\cM)$ is a full subcategory  of complexes in the equivariant derived category
whose derived restriction to the center of the unstable stratum has
cohomology of weights $0,\ldots,i-1$ for the wall
$M_{i-1}\dashrightarrow M_i$; see \cite[Theorem~3.21]{TT}.

\begin{lemma}\label{dfbadnsetn}
For $0\le k<i$, $\cP_{\cD^k_{i}}=\iota_i\circ \cP_{\cD^k_{i-1}}$.
Furthermore, objects in the subcategory $\langle\cD^k_{i-1}\rangle\subset D^b(M_{i-1})\cong\mathcal G_{[0,i)}$ have weights in the range $\{0,\ldots,k\}$.
\end{lemma}

Assume Lemma~\ref{dfbadnsetn}.  For $k<i$ and sufficiently small
$\epsilon>0$, we claim that there is a natural isomorphism
\[
 \mathcal P_{\mathcal D_i^{k,s}}
 \simeq \iota_i\circ\mathcal P_{\mathcal D_{i-\epsilon}^{k,s}}.
\]
Indeed, $\cD^{k,s}_{i-\epsilon}=\cD^k_{i-1}\otimes L^{k,s}_i$ and $\cD^{k,s}_{i}=\cD^k_{i}\otimes L^{k,s}_i$, where
$$L^{k,s}_i
=\begin{cases}
\cO\left(
\left\lfloor{s\over i-k}\right\rfloor, s+\left\lfloor{s\over i-k}\right\rfloor(k-1)\right)&\hbox{\rm if}\quad i-k>1,\cr
\cO(s,s(i-1)) &\hbox{\rm if}\quad i-k=1.\cr
\end{cases}
$$
The wall weight of
$L_i^{k,s}$ is equal to 
\[
 w_{k,s}:=
 \begin{cases}
 s\bmod(i-k),&i-k>1,\\
 0,&i-k=1,
 \end{cases}
\]
where $s\bmod(i-k)$ is the least nonnegative residue, so
$
 0\le w_{k,s}\le i-k-1$.
By Lemma~\ref{dfbadnsetn}, the weights of
$\mathcal D_{i-1}^k$ lie in $[0,k]$.  Therefore
\[
 \operatorname{weights}
 \bigl(\mathcal D_{i-1}^k\otimes L_i^{k,s}\bigr)
 \subset[w_{k,s},w_{k,s}+k]
 \subset[0,i-1].
\]
Restricting the same lifted kernel
to the two chambers shows that
$
 \mathcal P_{\mathcal D_i^{k,s}}
 \simeq
 \iota_i\circ\mathcal P_{\mathcal D_{i-\epsilon}^{k,s}}
$.
This analysis shows that the Farey Twill is compatible with the windows embeddings
$\iota:\,D^b(M_{i-1})\hookrightarrow D^b(M_i)$.

Next, we analyze what happens to the blocks $\langle\cD^{k,s}_{t}\rangle$ and  $\langle\cD^{k',s'}_{t}\rangle$ 
whose trajectories meet  at level $t$. 
This happens when
${s\over t-k}={s'\over t-k'}$ with one exception: as explained above, we modify
the paths of blocks $\langle\cD^{k,0}_{t}\rangle$
to be the horizontal line $t=k$ followed by the vertical line $x=k\epsilon$, $t>k$ for $0<\epsilon\ll1$.

\begin{lemma}\label{sgsrhsRHJ}
For two distinct blocks defined by~\eqref{wonderingblocks}, if
\[
 t\le g-1,\qquad
 \frac{s}{t-k}=\frac{s'}{t-k'}\notin\mathbb Z,
\]
then $\langle\cD^{k,s}_{t}\rangle$ and
$\langle\cD^{k',s'}_{t}\rangle$ are mutually orthogonal. Hence their
interchange does not change either embedded subcategory.
The same conclusion holds at $t=i$ for
$\langle\cD_i^{i,0}\rangle=\langle\cD_i^i\rangle$ and
$\langle\cD_i^{k,s}\rangle$, where $0\le k<i\le g-1$,
$0\le s\le3g-3-3k$, and $s/(i-k)\notin\mathbb Z$.
\end{lemma}

Next we analyze an integral crossing
$x_{k,s}(t)=s/(t-k)=n\in\bZ_{>0}$ with $t\le g-1$.
All intermediate pairs $(h,s-(h-k)n)$, $k\le h\le\ft$, satisfy
$0\le s-(h-k)n\le3g-3-3h$: use the bound on $s$ when $n\ge3$
and $t\le g-1$ when $n\le2$. We have
\begin{equation}\label{dfbadfna}
n={s\over t-k}={s-n\over t-(k+1)}={s-2n\over t-(k+2)}=\ldots={s-(\ft-k)n\over t-\ft}.
\end{equation}
Thus the paths of the blocks $\cD^{k,s}_{t},\cD^{k+1,s-n}_{t},\cD^{k+2,s-2n}_{t},\cD^{\ft,s-(\ft-k)n}_{t}$
meet at the same point.
If $k\ge1$ and $s+n\le 3g-3-3(k-1)$ then we also have $n={s+n\over t-(k-1)}$.
So, without loss of generality, we may assume that \eqref{dfbadfna} starts with the smallest possible~$k$.
At level $t$, the blocks of the subcategory 
\begin{equation}\label{	kwhrb,whjBG}
\langle\cD^{k,s}_{t},\cD^{k+1,s-n}_{t},\cD^{k+2,s-2n}_{t},\ldots,\cD^{\ft,s-(\ft-k)n}_{t}\rangle
\end{equation}
undergo a mutation resulting in the blocks
\begin{equation}\label{qerhqehqe5he}
\langle \cD^{\ft,s-(\ft-k)n}_{t+\epsilon}, \ldots, \cD^{k+2,s-2n}_{t+\epsilon}, \cD^{k+1,s-n}_{t+\epsilon},
\cD^{k,s}_{t+\epsilon}\rangle.
\end{equation}
Note that the order of the blocks is reversed (according to the slopes of the paths
\eqref{skhjgb,sHEG}). 
The twisting line bundle
$\cO\left(\left\lfloor{s\over t-k}\right\rfloor, s+\left\lfloor{s\over t-k}\right\rfloor(k-1)\right)$
changes precisely when ${s\over t-k}$ decreases and passes through an integer value~$n$. 

Put $r=s-(\lfloor t\rfloor-k)n$.  For the last block we have
$k'=\lfloor t\rfloor$ and hence, by the second line of
\eqref{wonderingblocks},
$
\cD_t^{\lfloor t\rfloor,r}
 =\cD_{\lfloor t\rfloor}^{\lfloor t\rfloor}
   \otimes\cO\bigl(r,r\lfloor t\rfloor\bigr)$.
The same formula holds immediately after the crossing.  (If $t\in\mathbb
Z$, then $r=0$.)  Thus the last block does not change.

\begin{figure}[htbp]
\includegraphics[width=0.8\textwidth]{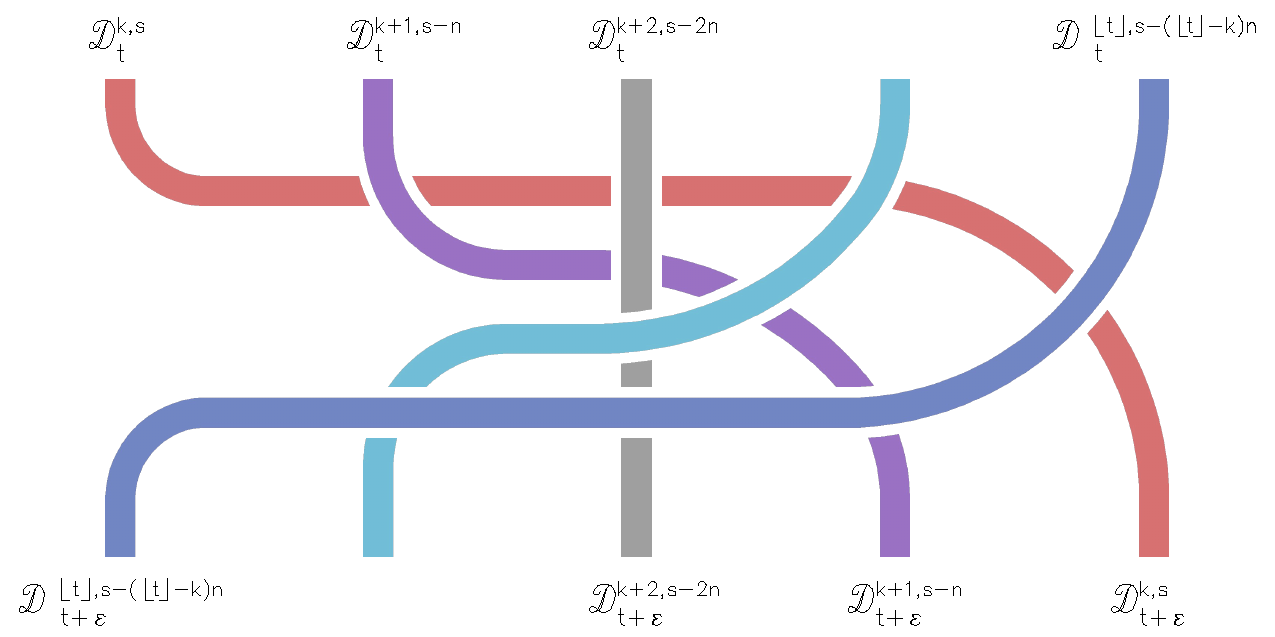}
\caption{Basic Farey Twill Mutation}\label{BasicTwill}
\end{figure}

\begin{lemma}\label{wrgwRHAREHEHR}
For the admissible chains described above with $t\le g-1$, the
mutation in $D^b(M_{\ft})$ from the semiorthogonal decomposition
\eqref{	kwhrb,whjBG} to the decomposition~\eqref{qerhqehqe5he}
is depicted in Figure~\ref{BasicTwill}.
\end{lemma}

If $t=i\in\mathbb Z$, the last fraction in \eqref{dfbadfna} is $0/0$.
By the boundary convention above, the corresponding endpoint block is
$\cD_i^{i,0}$, and it is included in the chain mutation.  Several chain
mutations may occur at level $i$.  Their subsequences
\eqref{	kwhrb,whjBG} are disjoint except for the common endpoint block
$\cD_i^{i,0}$, which participates consecutively, beginning with the
rightmost subsequence.

To prove Theorem~\ref{TwillTheorem}, we use the Farey Twill only
up to level $t=g-1$. First run the procedure for $t<g-1$.
At the last wall, apply the window embedding $\iota_{g-1}$ and add
the endpoint block $\cD^{g-1,0}_{g-1}=\cD^{g-1}$, without yet
performing any chain mutation at this final level.

This gives a semi-orthogonal decomposition of  $D^b(M_{g-1})$ with blocks
$\langle\cD^k_{g-1}\otimes L^{k,s}_{g-1}\rangle$ for 
$0\le k\le g-1$, $0\le s\le 3g-3-3k$, where 
$$L^{k,s}_{g-1}=\begin{cases}
\cO\left(
\left\lfloor{s\over g-1-k}\right\rfloor, s+\left\lfloor{s\over g-1-k}\right\rfloor(k-1)\right)&\hbox{\rm if}\quad k<g-1,\cr
\cO &\hbox{\rm if}\quad k=g-1.\cr
\end{cases}
$$
Note that the block $\cD^{g-1,0}_{g-1}$ has not participated in the Farey Twill yet. At level $t=g-1$, we do a single additional 
chain mutation
of  Lemma~\ref{wrgwRHAREHEHR}, from the sequence of blocks
$\langle \cD^{0,3g-3}_{g-1}, \cD^{1,3g-6}_{g-1}, \ldots, \cD^{g-2,3}_{g-1}, \cD^{g-1,0}_{g-1}\rangle$
to the sequence of blocks 
$\langle \cD^{g-1,0}_{g-1}, \cD^{g-2,3}_{g-1+\epsilon}, \ldots, \cD^{1,3g-6}_{g-1+\epsilon}, \cD^{0,3g-3}_{g-1+\epsilon}\rangle$.
This gives precisely the very  last blocks 
$\bigl\langle \cD^k(2,2g-4+j)\bigr\rangle_{ j+k=g-1\atop j,k\ge0}$
of the third mega-block
of Theorem~\ref{TwillTheorem}, in the correct order. 

It remains to explain how to obtain the remaining blocks in Theorem~\ref{TwillTheorem}.
These blocks are arranged into three mega-blocks
$$\Bigl\langle
\bigl\langle \cD^k(0,j)\bigr\rangle,\quad
\bigl\langle \cD^k(1,g-2+j)\bigr\rangle,\quad
\bigl\langle \cD^k(2,2g-4+j)\bigr\rangle\Bigr\rangle$$
that appear on the left of the semi-orthogonal decomposition.
In each of the three mega-blocks, $j,k\ge0$, $j+k\le g-2$, and the blocks are first arranged  by $j+k$
 and then by $j$ (both in increasing order).

These blocks come from $\langle\cD^{k,s}_{g-1}\rangle$ with
$0\le k\le g-2$ and $0\le s\le3g-4-3k$. Write
\[
 s=q(g-1-k)+j,\qquad q\in\{0,1,2\},\qquad 0\le j\le g-2-k.
\]
For $s>0$, stop the trajectory at $t=t(s,k)$ defined by
\[
 \frac{s}{t-k}=q+\frac{j+k}{g-1}.
\]
Then $k<t(s,k)\le g-1$, and the twisting line bundle is constant
from this time to $g-1$. All intervening crossings are nonintegral
and interchange mutually orthogonal blocks by Lemma~\ref{sgsrhsRHJ}.
Order the stopped blocks by their stopping values, that is, first by
$q$ and then by $j+k$. For fixed $q$ and $j+k$, the stopping time
increases with $j$, giving the required order within each row.
Place each block $\langle\cD^k\rangle$, $0\le k\le g-2$,
first in its row by commuting it past the blocks
$\langle\cD^a(0,j)\rangle$ with $j>0$ and $a+j<k$.
These interchanges are orthogonal: at $t=k$,
$0<j/(k-a)<1$, so the horizontal-endpoint case of
Lemma~\ref{sgsrhsRHJ} applies, and the subsequent windows embeddings
preserve orthogonality.
 No further integral chain mutation is
performed at $t=g-1$ beyond the $n=3$ chain described above.

\section{Cross Warp}\label{CrossWarpSection}

\begin{figure}[htbp]
\includegraphics[width=\textwidth]{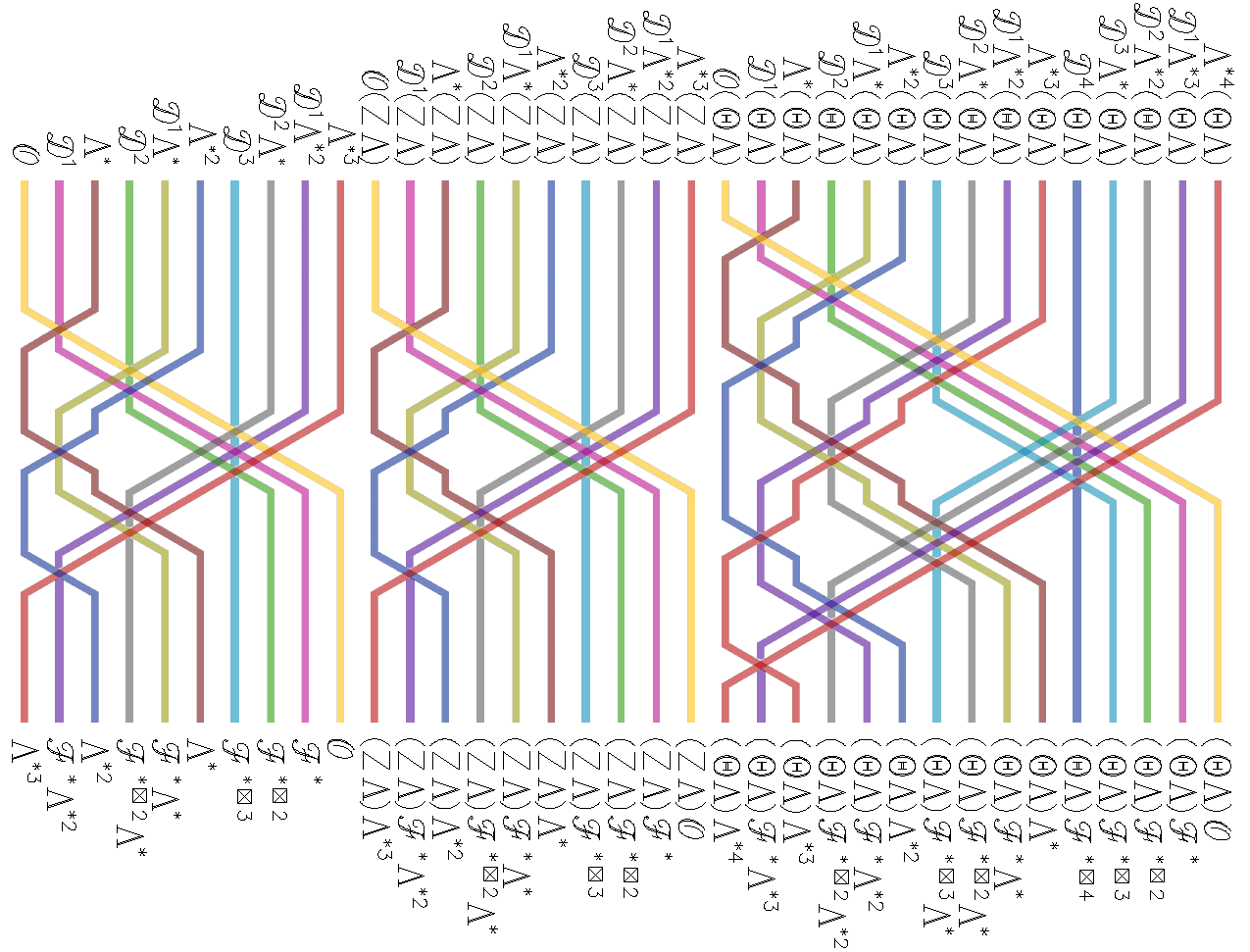}
\caption{Cross Warp in genus $5$}\label{genus5warp}
\end{figure}

In this section we will study a rather striking weaving pattern in $D^b(M)$, which is illustrated
in Figure~\ref{genus5warp} (in genus~$5$). We call it the Cross Warp.

The Cross Warp transforms the decomposition of
Theorem~\ref{TwillTheorem} into one whose kernels are tensor bundles
\(\cF^{*\boxtimes k}\) twisted by line bundles.

\begin{theorem}\label{CrossTheorem}
$D^b(M)$ has a semi-orthogonal decomposition into admissible subcategories
arranged into three mega-blocks, as follows:
$$
\Bigl\langle
\!\!\bigl\langle {\bLambda^*}^k{\cF^{*\boxtimes j}}\bigr\rangle_{ j+k\le g-2\atop j,k\ge0},\  
\bigl\langle (Z\bLambda)(\bLambda^*)^k\cF^{*\boxtimes j}\bigr\rangle_{j+k\le g-2\atop j,k\ge0},\ 
\bigl\langle (\theta\bLambda)(\bLambda^*)^k\cF^{*\boxtimes j}\bigr\rangle_{j+k\le g-1\atop j,k\ge0}
\!\!\Bigr\rangle.
$$
Within each mega-block, the blocks are ordered first by decreasing $k$ and then, for fixed $k$, by decreasing $j$.
\end{theorem}

The following theorem implies Theorem~\ref{CrossTheorem} and the lemmas 
of Section~\ref{TwillSection}.

\begin{theorem}\label{sGSRHSRH}
Let $0\le k\le i\le g-1$.
The Fourier--Mukai functors $\cP_{\cD_i^k}:\,D^b(\Sym^kC)\to D^b(M_i)$
introduced in Section~\ref{TwillSection} have the following properties:
\begin{enumerate}
\item[(a)] $\cP_{\cD_i^k}$ is fully faithful.
\item[(b)] 
$ \iota_i\circ\mathcal P_{\mathcal D_{i-1}^k}
 \simeq
 \mathcal P_{\mathcal D_i^k}$,
where $\iota_i$ 
is the windows
embedding ~\eqref{windowsembedding}, if $k<i$.
\item[(c)] The Fourier--Mukai functors $\cP_{\cD_i^k}, \ \cP_{\cF^{*\boxtimes k}}:\,D^b(\Sym^kC)\to D^b(M_i)$ are related via the Cross Warp mutation
illustrated in Figure~\ref{BasicCross}. 
\begin{figure}[htbp]
\includegraphics[width=0.8\textwidth]{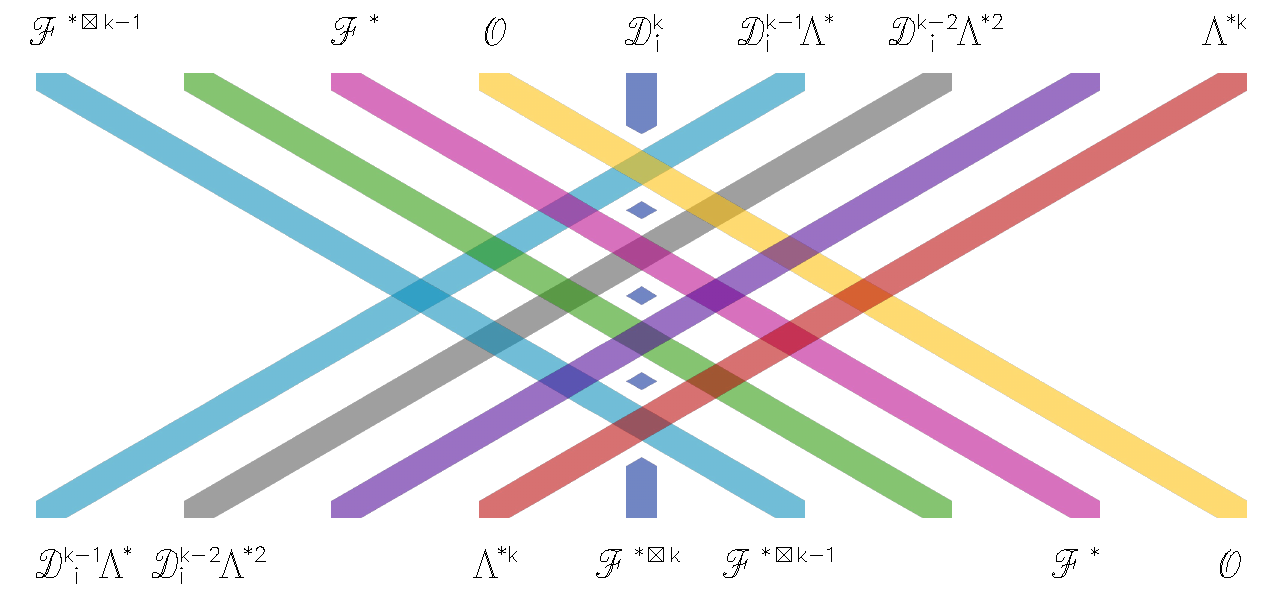}
\caption{Basic Cross Warp Mutation}\label{BasicCross}
\end{figure}
\end{enumerate}
\end{theorem}

\begin{remark}\label{Remark33}
This mutation combines two mutations as in Figure~\ref{ladyvanishes}. 
For $1\le k\le i$, the following subcategories $\mathcal B_k$ and $\mathcal A_k$
 are completely orthogonal:
\[
\mathcal B_k=
\bigl\langle\cF^{*\boxtimes(k-1)},\ldots,\cF^*,\cO\bigr\rangle,
\quad
\mathcal A_k=
\bigl\langle
\cD_i^{k-1}\bLambda^*,\ldots,
\cD_i^1{\bLambda^*}^{k-1},{\bLambda^*}^{k}
\bigr\rangle.
\]
The functorial triangle
\eqref{1srgseGWEG} identifies the right mutation of
$\cP_{\cF^{*\boxtimes k}}(X)[k]$ through $\mathcal B_k$
with $\cP_{\cF^{\bullet\boxtimes k}}(X)$, where the Fourier--Mukai kernel $\cF^{\bullet\boxtimes k}$ 
is defined in Definition~\ref{,aENFVkejhfv}.
The functorial
triangle \eqref{1Xfgsrhgasrha} gives a natural isomorphism
\[
 \mathbf L_{\mathcal A_k}\circ
 \mathcal P_{\mathcal F^{\bullet\boxtimes k}}
 \simeq
 \mathcal P_{\mathcal D_i^k}\circ
 \bigl(-\otimes
 \mathcal O_{\operatorname{Sym}^kC}(-\Delta_k/2)\bigr).
\]
Here $\Delta_k\subset\operatorname{Sym}^kC$ is the diagonal divisor of
nonreduced effective divisors, and
$\mathcal O_{\operatorname{Sym}^kC}(-\Delta_k/2)$ is the  line bundle
 introduced in
Lemma~\ref{amazinggrace}.
\end{remark}

We  will first prove various corollaries of Theorem~\ref{sGSRHSRH} 
and introduce additional notation,
before proving the theorem itself in the next section.

\begin{proof}[Proof of Theorem~\ref{CrossTheorem}]
Each of the three mega-blocks of the semi-orthogonal decomposition of Theorem~\ref{TwillTheorem} mutates
into the corresponding mega-block of the decomposition of Theorem~\ref{CrossTheorem},
without any interaction between the mega-blocks. 
Furthermore, the second and third mega-blocks differ from the first by
tensoring with $Z\otimes\boldsymbol\Lambda$ and
$\theta\otimes\boldsymbol\Lambda$, respectively.
In addition, the third mega-block is larger (it includes the block whose
Fourier--Mukai kernel is $\mathcal D^{g-1}$ whereas the first two mega-blocks continue only up to $\cD^{g-2}$).

So it is enough to describe the mutation for the first mega-block, which is illustrated in 
the left third of Figure~\ref{genus5warp}.
The  mutations of Theorem~\ref{sGSRHSRH}(c) fit together along their
boundary blocks.  
The mutation whose upper central block is $\cD^k$ shares
its upper-left block ${\mathcal F^{*\boxtimes (k-1)}}$ with the preceding mutation and its
lower-left block $\cD^{k-1}\otimes\bLambda^*$ with the next row.  Hence the
basic mutations concatenate as illustrated below to form the Cross Warp pattern.
$$
\begin{matrix}
\cO\cr
&\ldots&\cr
\ldots&&\cD^{g-5}&&&\cr
&\ldots&&\cD^{g-4}&&\cr
\ldots&&\cD^{g-5}\bLambda^*&&\cD^{g-3}&\cr
&\ldots&&\cD^{g-4}\bLambda^*&&\cD^{g-2}\cr
\ldots&&\cD^{g-5}{\bLambda^*}^2&&\cD^{g-3}\bLambda^*&\cr
&\ldots&&\cD^{g-4}{\bLambda^*}^2&&\cr
\ldots&&\cD^{g-5}{\bLambda^*}^3&&&\cr
&\ldots&\cr
{\bLambda^*}^{g-2}\cr
\end{matrix}
$$
Altogether this gives the mutation from the semi-orthogonal decomposition
$$\langle\cO\rangle,\langle\cD^1,\bLambda^*\rangle,\langle \cD^2,\cD^1\bLambda^*,{\bLambda^*}^2\rangle,\ldots,\langle\cD^{g-2},\cD^{g-3}{\bLambda^*},\ldots,\cD^1{\bLambda^*}^{g-3},{\bLambda^*}^{g-2}\rangle$$ of Theorem~\ref{TwillTheorem} to
 the semi-orthogonal decomposition
$$\langle{\bLambda^*}^{g-2}\rangle,\langle\cF^*{\bLambda^*}^{g-3},{\bLambda^*}^{g-3}\rangle,
\ldots,
\langle {\cF^{*\boxtimes g-2}}, {\cF^{*\boxtimes g-3}},\ldots,\cF^*,\cO\rangle$$ 
of Theorem~\ref{CrossTheorem}. The other two mega-blocks differ by  line bundle twists.
The third mega-block continues to $\cD^{g-1}$ rather than $\cD^{g-2}$.
\end{proof}

\begin{proof}[Proof of Lemma~\ref{asgasrhare}]
This is a reformulation of Theorem~\ref{sGSRHSRH}~(a).
\end{proof}

\begin{proof}[Proof of Lemma~\ref{dfbadnsetn}]
The first part is a reformulation of Theorem~\ref{sGSRHSRH}~(b).
The second part follows from Theorem~\ref{sGSRHSRH}~(c) and induction on $k$,
since the weights of $\cF^{*\boxtimes j}$ are $0,\ldots,j$ by \cite[Theorem~3.21]{TT}.
\end{proof}

\begin{notation}\label{notation:ordered-quotient}
Let $\tau_k:C^k\times M_i\to\Sym^kC\times M_i$ be the 
quotient map, and let $\pi_j:C^k\times M_i\to C\times M_i$ be the projection
onto the $j$th curve component and the moduli space factor.  The morphism $\tau_k$ is
finite flat; its formation and invariant direct image commute with
base change on $M_i$ \cite[Lemmas~2.2 and~2.4]{TT}.
A hat will be used for objects on $C^k\times M_i$.  In particular,
\[
 \widehat{\mathcal F}^{\boxtimes k}=
 \bigotimes_{j=1}^k\pi_j^*\mathcal F,\qquad
 \mathcal F^{\boxtimes k}=
 (\tau_{k*}\widehat{\mathcal F}^{\boxtimes k})^{S_k},
\]
with the ordinary permutation action.  For an equivariant complex $\hat K$,
put
\[
 \operatorname{Alt}_k(\hat K)=
 (\tau_{k*}(\hat K\otimes\operatorname{sgn}_k))^{S_k},\qquad
 \overline{\mathcal F}^{\boxtimes k}=
 \operatorname{Alt}_k(\widehat{\mathcal F}^{\boxtimes k}).
\]
\end{notation}

\begin{lemma}\label{amazinggrace}
For every $k\ge0$,
$
 \mathcal F^{*\boxtimes k}
 \simeq
 (\Lambda^{*\boxtimes k}\boxtimes
  \boldsymbol\Lambda^{*k})\otimes
 \mathcal F^{\boxtimes k},
$
and
\[
 (\mathcal F^{*\boxtimes k})^*
 \simeq
 \overline{\mathcal F}^{\boxtimes k}
 \otimes\mathscr L_k^*.
\]
Here $\mathscr L_k$ is the line bundle
$
 \mathscr L_k:=\operatorname{Alt}_k(\mathcal O_{C^k})=
 \mathcal O_{\operatorname{Sym}^kC}(-\Delta_k/2)$,
where  
 $\Delta_k\subset\operatorname{Sym}^kC$ is the diagonal divisor. We have $\mathscr L_0$, $\mathscr L_1=\mathcal O$.
In~particular,
$$\mathcal O(-\Delta_k/2)^{\otimes2}\simeq\mathcal O(-\Delta_k).$$
\end{lemma}

\begin{proof}
The first isomorphism follows from
$\mathcal F^*\simeq\mathcal F\otimes\det\mathcal F^*$,
$\det\mathcal F=\Lambda\boxtimes\boldsymbol\Lambda$, and the
projection formula for $\tau_k$.  The tensor power of the line bundle
$\Lambda^*$ descends with its ordinary permutation action; see
\cite[Lemma~3.13]{TT}.

Put $\widehat\Delta_k=\sum_{a<b}\Delta_{ab}$ on the ordered product $C^k$.
Evaluation of alternating functions gives the equivariant isomorphism
$
 \tau_k^*\mathscr L_k\otimes\operatorname{sgn}_k
 \simeq\mathcal O(-\widehat\Delta_k)
$.

Coherent duality for finite morphisms with smooth source and target gives
\begin{equation}\label{eq:relative-dualizing-linearization}
 \omega_{\tau_k}\simeq\mathcal O(\widehat\Delta_k)
 \simeq\tau_k^*\mathscr L_k^*\otimes\operatorname{sgn}_k.
\end{equation}
Consequently, by equivariant coherent duality and the projection formula,
\begin{align*}
 (\mathcal F^{*\boxtimes k})^*
 &\simeq
 \left(\tau_{k*}\bigl(
 (\widehat{\mathcal F}^{*\boxtimes k})^*
       \otimes\omega_{\tau_k}\bigr)\right)^{S_k}\\
 &\simeq
 \left(\tau_{k*}\bigl(
 \widehat{\mathcal F}^{\boxtimes k}
       \otimes\operatorname{sgn}_k\bigr)\right)^{S_k}
       \otimes\mathscr L_k^*
 =\overline{\mathcal F}^{\boxtimes k}\otimes\mathscr L_k^*.
\end{align*}
This is also the duality calculation in \cite[Lemma~3.13]{TT}.
\end{proof}

\begin{remark}
A general point of 
$\widehat\Delta_k\subset C^k$ lies on exactly one pairwise diagonal $\Delta_{ab}$ and has stabilizer $S_2$.  This stabilizer acts
nontrivially on the corresponding  fiber of $\widehat{\mathcal F}^{\boxtimes k}$; hence
$\mathcal F^{\boxtimes k}$ is not a
descent of $\widehat{\mathcal F}^{\boxtimes k}$.
\end{remark}

\begin{corollary}
The Fourier--Mukai functors $D^b(\Sym^kC)\to D^b(M)$
with kernels 
$(\cF^{*\boxtimes k})^*$ and $(\bar\cF^{\boxtimes k})$
differ by an autoequivalence of $D^b(\Sym^kC)$ given by tensoring with a line bundle $\mathscr L_k$ on $\Sym^kC$. 
Consequently, we have a noncanonical isomorphism
$
(\cF_D^{*\boxtimes k})^*
 \simeq
\bar\cF_D^{\boxtimes k}$ for any point $D\in \Sym^kC$.
\end{corollary}

\begin{notation}\label{notation:other-determinants}
For a line bundle $\Lambda'$ of degree $d>0$ and
$0\le i\le\lfloor{d-1\over2}\rfloor$, write $M_i(\Lambda')$, or $M_i(d)$
when the determinant is not important, for the Thaddeus space of stable pairs with the
stability convention of \cite[Section~3]{TT}.

For $0\le k\le i\le g-1$, the projection $D_i^k\to\Sym^kC$ is smooth
with fiber $M_{i-k}(\Lambda(-2D))$ over $D$. In particular,
$D_k^k\to\Sym^kC$ is a projective bundle. For $i\ge1$ and $k<i$,
$
 \left.\cO(m,n)\right|_{M_{i-k}(\Lambda(-2D))}
 \simeq\cO(m,n-km)$.
For $i=k\ge1$, the fiber is $\bP^{3g-3-2k}$ and the restriction is
$
 \left.\cO(m,n)\right|_{M_0(\Lambda(-2D))}
 \simeq\cO_{\bP^{3g-3-2k}}\bigl(n+m(1-k)\bigr)
$.
These are the  formulas of \cite[Remark~3.7]{TT}.
\end{notation}

We use  vanishing theorems for tensor vector bundles  from \cite{TT},
which we state here for ease of reference.

\begin{theorem}[{\cite[Theorem 7.1]{TT}}]\label{hardvanishing}
Suppose $2<d\leq 2g+1$, and $1\leq j\leq \lfloor {d-1\over 2}\rfloor$. 
Let $D\in\Sym^a C$ and $D'\in\Sym^bC$ with $a,\ b\leq d+g-2j-1$, and let $t$ be an integer satisfying
$
    a-j-1 < t < d+g-2j-1-b
$.
Suppose, further, that 
$t\notin [0,a]$. 
Then we have
\begin{equation}\label{afgaawrhr}
R\Gamma_{M_j(d)}({{\bar\cF^{\boxtimes a}_D}}{}^*\otimes \bar\cF^{\boxtimes b}_{D'}\otimes\bLambda^t)=0
\end{equation}
and, if $D=\sum\alpha_i x_i$,
\begin{equation}\label{afgaasfargarhwrhr}
R\Gamma_{M_j(d)}(
\bigotimes_i{\cF^*_{x_i}}^{\otimes \alpha_i}\otimes
\bar\cF^{\boxtimes b}_{D'}
\otimes\bLambda^t)
=0.
\end{equation}
\end{theorem}

\begin{theorem}[{\cite[Lemma 7.3, Theorem 7.4]{TT}}]\label{TT_Theorem7.4}
Let  $d>2$ and $1\leq j\leq \lfloor {d-1\over 2}\rfloor$. 
Let $D\in\Sym^a C$ and $D'\in\Sym^bC$,
and let $t$ be an integer satisfying
$
    a < t < d+g-2j-1-b$.
Then 
\begin{equation}
R\Gamma_{M_j(d)}({\bar\cF^{\boxtimes a}_D}{}^*\otimes {\bar\cF^{\boxtimes b}_{D'}}\otimes\bLambda^t)=
R\Gamma_{M_j(d)}({\cF^{\boxtimes a}_D}^*\otimes\cF^{\boxtimes b}_{D'} \otimes\bLambda^t)=
0.
\end{equation}
Moreover, the same vanishing holds for $d>0$ and $j=0$.
\end{theorem}

\begin{theorem}[{\cite[Theorems~9.3 and~9.7]{TT}}]\label{TT_Theorem9.6}
Suppose $2<d\le 2g-1$ and
$1\le j\le\lfloor(d-1)/2\rfloor$.  Let
$D\in\Sym^aC$ and $D'\in\Sym^bC$, with
$a\le j$ and $b<d+g-2j-1$.  If $D\not\le D'$ (in particular, if
$a>b$), then
\[
R\Gamma_{M_j(d)}
\bigl((\bar\cF_D^{\boxtimes a})^*
      \otimes\bar\cF_{D'}^{\boxtimes b}\bigr)=0.
\]
\end{theorem}

\begin{remark}
Here Theorem~9.7 of \cite{TT} uses full faithfulness in smaller tensor
degrees, supplied in the stated range by Theorem~9.3 of \cite{TT}.
\end{remark}


\begin{proof}[Proof of Lemma~\ref{sgsrhsRHJ}]
Put $i=\ft$. For an ordinary crossing, let
$r=s/(t-k)=s'/(t-k')\notin\bZ$ and $k'>k$.
For a crossing with the horizontal endpoint at $t=i$, set
$k'=i$, $s'=0$, and $r=s/(i-k)\notin\bZ$.
In both cases, $r(k'-k)=s-s'\in\bZ$, so
$\{r\}(k'-k)\in\bZ$ and $k'\ge k+2$.
$\cD^{k,s}_t=
\cD^k_i\left(\fr, s+\fr(k-1)\right)$ 
and
$\cD^{k',s'}_t=\cD^{k'}_{i}\otimes L$,
where 
$$L=\begin{cases}
\cO\left(
\fr, s'+\fr(k'-1)\right)&\hbox{\rm if}\quad k'<i,\cr
\cO(s',s'k') &\hbox{\rm if}\quad k'=i.\cr
\end{cases}
$$
We claim that $\RHom(\cP_{\cD^{k,s}_t}(X),\cP_{\cD^{k',s'}_t}(X'))=0$ for any $X\in D^b(\Sym^kC)$, $X'\in D^b(\Sym^{k'}C)$.
By Theorem~\ref{sGSRHSRH} (c), it suffices to show that 
$$\RHom(\cP_{{\bLambda^*}^l{\cF^{*\boxtimes m}}}(X)\left(\fr, s+\fr(k-1)\right),\cP_{\cD^{k',s'}_t}(X'))=0,$$
where $X\in D^b(\Sym^mC)$, $X'\in D^b(\Sym^{k'}C)$
and  $l+m\le k$, $l,m\ge0$.
Furthermore, we can assume that $X$ (resp., $X'$)
is the skyscraper sheaf of a point $D\in \Sym^mC$ (resp.,~$D'\in\Sym^{k'}C$), and so we need to check that 
$$\RHom({\bLambda^*}^l{\cF_D^{*\boxtimes m}}\left(\fr, s+\fr(k-1)\right),\cO_{M_{i-k'}(\Lambda(-2D'))}\otimes L)=0.$$
We first consider the case when $\boxed{k'=i}$, so that $M_{i-k'}(\Lambda(-2D'))=\bP^{3g-3-2k'}$ and we need to show that 
\begin{equation}\label{Sgasrhasrhasrha}
R\Gamma(\bP^{3g-3-2k'},\ 
{\bLambda}^l({\cF_D^{*\boxtimes m}})^*\cO(s'-s-\fr k+\fr k'))=
\end{equation}
$$R\Gamma(\bP^{3g-3-2k'},\ 
{({\cF_D^{\boxtimes m}})^*\bLambda}^{l+m-(s'-s-\fr k+\fr k')})=0.$$
But $s'-s-\fr k+\fr k'=\{r\}(k-k')<0$ since $r$ is not an integer.
On the other hand, $l+m+\{r\}(k'-k)\le k'\le 3g-3-2k'$. So we have \eqref{Sgasrhasrhasrha}
by Theorem~\ref{TT_Theorem7.4} with $j=0$, $d=2g-1-2k'$, $a=m$, $b=0$, $t=l+m+\{r\}(k'-k)$.

Now suppose $\boxed{k'<i}$. Arguing as above, we need to show that
$$R\Gamma(M_{i-k'}(\Lambda(-2D')),\ 
{\cF_D^{\boxtimes m}}^*{\bLambda}^{t})=0,$$
where $t=l+m+\{r\}(k'-k)$. 
This  follows from Theorem~\ref{TT_Theorem7.4} with parameters
$a=m$, $b=0$, $j=i-k'$, $d=2g-1-2k'$, $t=l+m+\{r\}(k'-k)$.
Indeed, put $\delta=\{r\}(k'-k)$.  Since $r\notin\mathbb Z$ and
$\delta\in\mathbb Z$, one has $1\le\delta\le k'-k-1$.  Therefore
$
m<t=l+m+\delta\le k'-1$.
Moreover,
\[
1\le j=i-k'\le g-1-k'
 =\left\lfloor\frac{d-1}{2}\right\rfloor,
\qquad d=2g-1-2k'\ge3,
\]
and
\[
t\le k'-1<i\le g-1<3g-2i-2=d+g-2j-1.
\]
Thus every hypothesis of Theorem~\ref{TT_Theorem7.4} is satisfied.

The opposite vanishing follows from the order of the semi-orthogonal decomposition immediately before the crossing, which has already been constructed. The two kernels do not change at this crossing. Hence the two blocks are mutually orthogonal and may be interchanged.
\end{proof}

\begin{proof}[Proof of Lemma~\ref{wrgwRHAREHEHR}]
Write $p:D_i^k\to\Sym^kC$ and let
\[
 E_i^{k,1}=\{(D,F,s)\in D_i^k:\deg Z(s)\ge k+1\},
 \qquad \mathcal E_i^{k,1}=\mathcal O_{E_i^{k,1}}.
\]
We use its reduced relative divisor structure; for $k=i$ it is empty.
On the fiber over $D$, this is the divisor of pairs with a zero in
$M_{i-k}(\Lambda(-2D))$.  

By \cite[Remark~3.7 and Definition~6.4]{TT},
$
 \mathcal O_{D_i^k}(1,k-1)(-E_i^{k,1})
$
is trivial on every fiber of $p$, including the projective-space fibers
when $k=i$.  Put
$
 L=p_*\bigl(\mathcal O_{D_i^k}(1,k-1)(-E_i^{k,1})\bigr).
$
The fibers of $p$ are connected smooth projective Thaddeus spaces, and
$h^0(\mathcal O)=1$ on each fiber.  Cohomology and base change therefore
make $L$ invertible and the evaluation morphism an isomorphism
$
 p^*L\simeq\mathcal O_{D_i^k}(1,k-1)(-E_i^{k,1})
$.
The divisor sequence, pushed into $\Sym^kC\times M_i$, is thus
\[
 0\longrightarrow
 \mathcal D_i^k(-1,-(k-1))\otimes p^*L
 \longrightarrow\mathcal D_i^k
 \longrightarrow\mathcal E_i^{k,1}\longrightarrow0.
\]
Here the tensor product with $p^*L$ is taken on $D_i^k$ before the
closed pushforward.  For $X\in D^b(\Sym^kC)$ this gives
\begin{equation}\label{afbasfbab}
 \mathcal P_{\mathcal D_i^k}(X\otimes L)(-1,-(k-1))
 \longrightarrow\mathcal P_{\mathcal D_i^k}(X)
 \longrightarrow\mathcal P_{\mathcal E_i^{k,1}}(X)
 \longrightarrow .
\end{equation}
Tensoring by $L$ is an autoequivalence of the source.  It must be
retained in this functorial triangle, but it does not change the
essential-image subcategory used in the mutation.

\begin{claim}\label{qdfberb}
$\cP_{\cE_i^{k,1}}(X)\in\langle\cD_i^{k+1},\ldots,\cD_i^i\rangle$.
\end{claim}

We will prove Claim~\ref{qdfberb} by induction. 

\begin{definition}\label{SGsrhsrha}
For $k_1,\ldots,k_r\ge0$ with $k_1+\cdots+k_r\le i$, define
\[
 D_i^{k_1,\ldots,k_r}=
 \{(D_1,\ldots,D_r,F,s):s|_{D_1+\cdots+D_r}=0\}
 \subset\prod_{a=1}^r\Sym^{k_a}C\times M_i.
\]
For $k_r\ge1$ and $k_1+\cdots+k_{r-1}\le i$, let
\[
 E_i^{k_1,\ldots,k_r}=
 \{(D_1,\ldots,D_{r-1},F,s)\in D_i^{k_1,\ldots,k_{r-1}}:
        \deg Z(s)\ge k_1+\cdots+k_r\},
\]
with its reduced structure.  It is empty when
$k_1+\cdots+k_r>i$.  Zero-degree factors are points and may be omitted;
$D_i^{\varnothing}=M_i$.  The corresponding structure sheaves are
denoted by $\mathcal D_i^{k_1,\ldots,k_r}$ and
$\mathcal E_i^{k_1,\ldots,k_r}$.
\end{definition}

For $k_r\ge1$ and $k_1+\cdots+k_r\le i$, the locus $D_i^{k_1,\ldots,k_r}$ is smooth and 
\begin{equation}\label{normalization of Ei}
\nu:\, D_i^{k_1,\ldots,k_r}\to E_i^{k_1,\ldots,k_r},\quad (D_1,\ldots,D_r,F,s)\mapsto  (D_1,\ldots,D_{r-1},F,s)
\end{equation}
is a normalization morphism. 

\begin{claim}
The relative form of \cite[Lemma~6.5]{TT} holds: 
we have a commutative diagram
\begin{equation}\label{conductor}
    \begin{tikzcd}[column sep = 0.3cm]
    0 \arrow[r] & \nu_*\cD_i^{k_1,\ldots,k_r}(-E_i^{k_1,\ldots,k_r,1}) \arrow[r] &\nu_*\cD_i^{k_1,\ldots,k_r}\arrow[r] &\nu_*\cE_i^{k_1,\ldots,k_r,1} \arrow[r] &0\\
    0\arrow[r] &\mathcal{I}_{E_i^{k_1,\ldots,k_{r-1},k_r+1}}\arrow{u}[sloped]{\sim}\arrow[r] &\cE_i^{k_1,\ldots,k_r}\arrow[u,hook]\arrow[r] &\cE_i^{k_1,\ldots,k_{r-1},k_r+1}\arrow[u,hook]\arrow[r] &0
    \end{tikzcd}
\end{equation}
where $\mathcal{I}_{E_i^{k_1,\ldots,k_{r-1},k_r+1}}\simeq \nu_*\cD_i^{k_1,\ldots,k_r}(-E_i^{k_1,\ldots,k_r,1})$ is the conductor sheaf of the normalization (\ref{normalization of Ei}) and $E_i^{k_1,\ldots,k_r,1}$ (resp.~$E_i^{k_1,\ldots,k_{r-1},k_r+1}$)
is a conductor subscheme in $D_i^{k_1,\ldots,k_r}$ (resp.~$E_i^{k_1,\ldots,k_r}$).
\end{claim}

\begin{proof}
Fix $k_1,\ldots,k_{r-1}$ and put $m=\sum_{b<r}k_b$.
For $1\le a\le i-m$, write
\[
 T=D_i^{k_1,\ldots,k_{r-1}},\qquad
 D_a=D_i^{k_1,\ldots,k_{r-1},a},\qquad
 E_a=E_i^{k_1,\ldots,k_{r-1},a}.
\]
After dividing the universal section by $D_1+\cdots+D_{r-1}$,
$D_a$ parametrizes degree-$a$ sub-divisors of its zero scheme, and
$E_a$ is its reduced image in $T$. The schemes $T$ and $D_a$ are
smooth, with $\dim D_a=\dim T-a$. The forgetful map
$\nu_a:D_a\to E_a$ is finite and is an isomorphism over
$E_a\setminus E_{a+1}$, where the residual zero divisor has degree
exactly $a$. Thus $\nu_a$ is the normalization.

Put $E_{i-m+1}=\varnothing$ and
$R_a=E_i^{k_1,\ldots,k_{r-1},a,1}\subset D_a$.
This is the reduced inverse image of $E_{a+1}$ and is a Cartier
divisor, or is empty. After fixing all selected divisors, the fibers
of $D_a$ are the residual Thaddeus spaces, and $R_a$ restricts to
their reduced zero-divisor loci. Indeed, the restriction is a
Cartier divisor, generically reduced by the simple-zero incidence
description. This also identifies the restriction of
$\mathcal O_{D_a}(-R_a)$ used in~\eqref{afbasfbab}.

At a general point of each component of $E_{a+1}$, the residual
zero divisor consists of $a+1$ distinct simple points.
\'Etale-locally, persistence of these zeros is given by functions
$h_1,\ldots,h_{a+1}$ on $T$. Their common zero locus is smooth of
codimension $a+1$, since the corresponding incidence scheme is
smooth. Thus they are part of a regular system of parameters and
\[
 \mathcal I_{E_a/T}
 =\bigcap_{b=1}^{a+1}
   (h_1,\ldots,\widehat h_b,\ldots,h_{a+1}).
\]
Transversally to $E_{a+1}$, this is the union of the coordinate
axes in $\mathbb A^{a+1}$, as in \cite[Claim~6.6]{TT}.
We now induct on $a$. The reduced divisor $E_1\subset T$ is
Cohen--Macaulay. If $E_a$ is Cohen--Macaulay, it satisfies $S_2$;
the preceding codimension-one description therefore implies that
$E_a$ is seminormal. Its conductor is consequently reduced.
Its support is $E_{a+1}$ downstairs and $R_a$ on the normalization,
so
$
 \mathcal I_{E_{a+1}/E_a}
 =\nu_{a*}\mathcal O_{D_a}(-R_a)$.
This gives the two rows of~\eqref{conductor}. The sheaf on the right
is Cohen--Macaulay, since $\nu_a$ is finite and $D_a$ is smooth.
The exact sequence
\[
 0\longrightarrow\nu_{a*}\mathcal O_{D_a}(-R_a)
 \longrightarrow\mathcal O_{E_a}
 \longrightarrow\mathcal O_{E_{a+1}}\longrightarrow0
\]
and the depth lemma show that $E_{a+1}$ is Cohen--Macaulay.
This proves the relative form of \cite[Lemma~6.5]{TT}.
\end{proof}

Instead of  Claim~\ref{qdfberb}, we will prove a slightly stronger 
\begin{claim}
For every  $X\in D^b(\Sym^{k_1}C\times\ldots\times \Sym^{k_{r-1}}C)$, we have
$$\cP_{\cE_i^{k_1,\ldots,k_r}}(X)\in\langle\cD_i^{k_1+\ldots+k_r},\ldots,\cD_i^i\rangle.$$
\end{claim}

\begin{proof}[Proof of the claim]
We argue by  downward induction on $k_1+\ldots+k_r$.
By inductive assumption and \eqref{conductor}, it suffices to prove the same statement for 
the Fourier--Mukai functors with kernels
$\nu_*\cD_i^{k_1,\ldots,k_r}$ and $\nu_*\cE_i^{k_1,\ldots,k_r,1}$.
By~projection formula and inductive assumption, it suffices to prove that 
$\cP_{\cD_i^{k_1,\ldots,k_r}}(X)\in\langle\cD_i^{k_1+\ldots+k_r},\ldots,\cD_i^i\rangle$
for  $X\in D^b(\Sym^{k_1}C\times\ldots\times \Sym^{k_{r}}C)$.

In fact, consider projections $\pi_1:\,D_i^{k_1,\ldots,k_r}\to \Sym^{k_1}C\times\ldots\times \Sym^{k_{r}}C$
and  $\pi_2:\,D_i^{k_1,\ldots,k_r}\to M_i$. It suffices to show  that $R{\pi_2}_*(L\pi_1^*X)\in \langle\cD_i^{k_1+\ldots+k_r}\rangle$.
We have a Cartesian diagram
$$    \begin{tikzcd}
    D_i^{k_1,\ldots,k_r} \arrow["\pi_1",d]\arrow[r,"\mu"] & D_i^{k_1+\ldots+k_r}\arrow[d,"\pi_1"]\\
    \Sym^{k_1}C\times\ldots\times\Sym^{k_r}C\arrow[r,"\mu"] &\Sym^{k_1+\ldots+k_r}C
     \end{tikzcd}
$$
where $\mu$ is the addition map $(D_1,\ldots,D_r)\mapsto D_1+\ldots+D_r$.
Furthermore, 
$R{\pi_2}_*(L\pi_1^*X)=R{\pi_2}_*(R\mu_*L\pi_1^*X)=R{\pi_2}_*(L\pi_1^*(R\mu_*X))$
by cohomology and base change. But the latter is contained in $\langle\cD_i^{k_1+\ldots+k_r}\rangle$
by definition.
\end{proof}

To finish the proof that $\langle\cD_i^k(-1,-(k-1))\rangle$
is the mutation of $\langle\cD_i^k\rangle$, we need to show that
$\langle\cD_i^k(-1,-(k-1))\rangle\subset{}^\perp\langle\cD_i^{k+1},\ldots,\cD_i^i\rangle$.
It suffices to show that, for $X\in D^b(\Sym^kC)$,  $X'\in D^b(\Sym^\alpha C)$, $k<\alpha\leq i$, 
we have
$\RHom(\cP_{\cD_i^k}(X)(-1,-(k-1)), \cP_{\cD_i^\alpha}(X'))=0$.

Apply~\eqref{afbasfbab} with $k=\alpha$ and use downward induction
on $\alpha$, together with Claim~\ref{qdfberb}.
As the argument ranges over every $X'\in D^b(\Sym^\alpha C)$, the source
line bundle in that triangle can be absorbed by replacing $X'$ with its
tensor product with that line bundle.  It therefore suffices to prove
\[
 \RHom\bigl(
 \mathcal P_{\mathcal D_i^k}(X)(-1,-(k-1)),
 \mathcal P_{\mathcal D_i^\alpha}(X')(-1,-(\alpha-1))
 \bigr)=0.
\]

By Theorem~\ref{sGSRHSRH}(c), it is enough instead to prove that
$$\RHom(\cP_{{\bLambda^*}^l{\cF^{*\boxtimes m}}}(X)(-1,-(k-1)),\cP_{\cD_i^\alpha}(X')(-1,-(\alpha-1)))=0$$
for $l+m\le k$, $l,m\ge0$. Furthermore, 
we can assume that $X$ (resp., $X'$)
is a skyscraper sheaf of a point $D\in \Sym^mC$ (resp.,~$D'\in\Sym^{\alpha}C$).

After restricting to $M_{i-\alpha}(\Lambda(-2D'))$, the difference
of the two ambient twists is $\mathcal O(0,k-\alpha)$.
Up to tensoring with a fixed one-dimensional vector space, the
required group is therefore
\[
 R\Gamma\!\left(M_{i-\alpha}(\Lambda(-2D')),
       (\mathcal F_D^{\boxtimes m})^*
       \otimes\boldsymbol\Lambda^{l+m+\alpha-k}\right).
\]
Apply Theorem~\ref{TT_Theorem7.4} with
$
 j=i-\alpha$, $d=2g-1-2\alpha$,
 $a=m$, $b=0$, $t=l+m+\alpha-k$.
If $j=0$, its projective space case applies since $d\ge1$.
If~$j>0$, then $d\ge3$ and
$1\le j\le\lfloor(d-1)/2\rfloor$.
Furthermore $\alpha>k$ gives $t>m$, and
$
 t=l+m+\alpha-k\le\alpha\le i
       <3g-2-2i=d+g-2j-1
$.
These are the required inequalities, so the cohomology vanishes.
\end{proof}

\section{Analysis of the basic Cross Warp mutation}\label{sgbsdgbsdst}

The goal of this section is to prove Theorem~\ref{sGSRHSRH}.
We use the notation of Section~\ref{CrossWarpSection}. Throughout this
section $0\le k\le i\le g-1$; in particular,
\begin{equation}\label{eq:uniform-chamber-bound}
 i<3g-2i-2.
\end{equation}
We first establish the semiorthogonality required for the basic Cross Warp mutation in
Figure~\ref{BasicCross}.
To prove orthogonality of Fourier--Mukai functors, it suffices to
check their values on skyscraper sheaves, since these form spanning
classes and the functors have both adjoints
\cite[Propositions~3.17 and~5.9]{huybrechts}.

\begin{lemma}\label{wRGwrhgwrHR}
For $1\le k\le i$, $D^b(M_i)$ admits admissible subcategories
$$\langle {\cF^{*\boxtimes k-1}},\ldots,{\cF}^*,\cO,{\bLambda^*}^k\rangle\quad\hbox{\rm and}\quad 
\langle {\bLambda^*}^k,{\cF^{*\boxtimes k}},{\cF^{*\boxtimes k-1}},\ldots,{\cF}^*,\cO\rangle.$$
\end{lemma}

\begin{proof}
Full faithfulness follows from \cite[Theorem~9.3]{TT}. By
Lemma~\ref{amazinggrace}, the required orthogonality reduces to
\[
\begin{aligned}
 R\Gamma\bigl((\bar\cF_D^{\boxtimes a})^*
               \otimes\bar\cF_{D'}^{\boxtimes b}\bigr)&=0
       &&(0\le b<a\le k),\\
 R\Gamma\bigl(\bar\cF_D^{\boxtimes l}\otimes\bLambda^{-k}\bigr)&=0
       &&(0\le l\le k),\\
 R\Gamma\bigl(\cF_D^{*\boxtimes l}\otimes\bLambda^k\bigr)&=0
       &&(0\le l<k).
\end{aligned}
\]
Apply Theorems~\ref{TT_Theorem9.6}, \ref{hardvanishing}, and
\ref{TT_Theorem7.4}, respectively, with $j=i$ and $d=2g-1$.
Their hypotheses follow from $1\le k\le i<3g-2i-2$;
in the first line $D\not\le D'$ because $a>b$.
\end{proof}

We prove Theorem~\ref{sGSRHSRH} by induction on $k$, starting with
$k=0$. The mutation Lemmas~\ref{1adfbadfnbadn} and~\ref{1swgasrgaerh}
apply for every $k\ge1$ and invoke the induction hypothesis for
Theorem~\ref{sGSRHSRH} only in smaller degree. They also use the
cohomology calculation of Lemma~\ref{sRGwgwrG} in the current degree,
proved in Section~\ref{asrgarsgargaerhaehae};
Lemma~\ref{1swgasrgaerh} uses Lemma~\ref{1adfbadfnbadn} in the current
degree as well. We first define the complex that will occur in these
mutations.

The universal section $\Sigma$ of the universal stable pair $(\cF,\Sigma)$ on $C\times M_i$
vanishes on the locus $D_i^1$, which has codimension~$2$. 
It follows that  we have an exact Koszul complex on $C\times M_i$,
\begin{equation}\label{sRARSHADTNADT}
0\to \Lambda^*\boxtimes\bLambda^*\to\cF^*\arrow^\Sigma\cO\to\cD_i^1\to0.
\end{equation}

\begin{definition}\label{def:universal-two-term}
We define the complex $\cF^\bullet=[\cF^*\arrow^\Sigma\cO]$ in $D^b(C\times M_i)$, normalized so that $\cO$ is placed in cohomological degree~$0$. 
\end{definition}

By exactness of \eqref{sRARSHADTNADT}, we compute 
\begin{equation}\label{sbasrhasrha}
\cH^{0}(\cF^\bullet)=\cD_i^1\quad\hbox{\rm and}\quad \cH^{-1}(\cF^\bullet)=\Lambda^*\boxtimes\bLambda^*.
\end{equation}
This gives two exact triangles in $D^b(C\times M_i)$,
\begin{equation}\label{dfghsdrgadrghadrgha}
\Lambda^*\boxtimes\bLambda^*[1]\to\cF^\bullet\to\cD_i^1\to
\qquad\hbox{\rm and}\qquad 
\cO\to\cF^\bullet\to\cF^*[1]\to.
\end{equation}

For arbitrary $k$, we use the following ordered incidence scheme.

\begin{definition}\label{def:ordered-incidence}
For $0\le k\le i$, put
$\widehat D_i^k=D_i^k\times_{\Sym^kC}C^k$ and
$\widehat{\mathcal D}_i^k=\mathcal O_{\widehat D_i^k}$.
These are $D_i^{1,\ldots,1}$ and $\cD_i^{1,\ldots,1}$, respectively,
in notation of Definition~\ref{SGsrhsrha}.
\end{definition}

\begin{remark}
The morphism $D_i^k\to\operatorname{Sym}^kC$ is smooth, and its
embedding in $\operatorname{Sym}^kC\times M_i$ is regular of codimension
$2k$.  Flat base change gives the analogous assertions for
$\widehat D_i^k\to C^k$.
The scheme $\widehat D_i^k$ is not, in general, the scheme-theoretic
intersection $\bigcap_a\pi_a^{-1}(D_i^1)$.  They agree when the ordered
points are distinct.  At a repeated point the scheme $\widehat D_i^k$ imposes vanishing
to the sum of the multiplicities, whereas $\bigcap_a\pi_a^{-1}(D_i^1)$ imposes only the
vanishing at the individual points.  This scheme can
therefore have additional components along the diagonals, and the full
Koszul complex of $(\pi_1^*\Sigma,\ldots,\pi_k^*\Sigma)$ need not be
exact.  Rather than working with it directly, we will instead compute the cohomology of the tensor product of
the two-term complexes $[\pi_a^*\mathcal F^*\to\mathcal O]$.
\end{remark}

\begin{definition}\label{,aENFVkejhfv}
We  define
$${\cF^{\bullet\boxtimes k}}=\operatorname{Alt}_k({\hat\cF^{\bullet\boxtimes k}}),\quad\hbox{\rm where}\quad {\hat\cF^{\bullet\boxtimes k}}=L\pi_1^*\cF^\bullet\otimes^L\ldots\otimes^LL\pi_k^*\cF^\bullet.$$
Recall that the symmetric group acts on ${\hat\cF^{\bullet\boxtimes k}}$ in such a way that
an adjacent transposition $(t,t+1)$, in addition to permuting the corresponding factors of a homogeneous tensor
$g_1\otimes\ldots\otimes g_k$, 
also multiplies it by $(-1)^{\deg(g_t)\deg(g_{t+1})}$.
Since $\mathcal F^\bullet$ is a complex of vector bundles,
$\mathcal F^{\bullet\boxtimes k}$ is represented by a locally free
complex in degrees $[-k,0]$: finite flat pushforward preserves local
freeness and taking invariants is a direct summand in characteristic
zero.  
\end{definition}

\begin{definition}\label{qetgathaehehe}
Let $0\le l\le k$, put $r=k-l$, and let
\begin{equation}\label{asfvasfgasgarga}
W_{k,l}:=
 \{(D,D')\in\operatorname{Sym}^kC\times\operatorname{Sym}^lC
       \mid D\ge D'\}.
\end{equation}
Write
$\rho_{k,l}:W_{k,l}\to\operatorname{Sym}^rC$ for the difference map
$(D,D')\mapsto D-D'$.
\end{definition}

Lemma~\ref{sRARHAEJ} concerns a brutal truncation of the 
complex $\mathcal F^{\bullet\boxtimes k}$. A rather delicate
Lemma~\ref{sRGwgwrG}, proved  in
Section~\ref{asrgarsgargaerhaehae}, computes its cohomology sheaves.
Lemma~\ref{sRHwrhjeRJAETJ} then follows from canonical truncation.

\begin{lemma}\label{sRARHAEJ}
Let
$
 G_k:=\sigma_{\ge-k+1}
       (\mathcal F^{\bullet\boxtimes k})
$
be the brutal truncation of the
 complex ${\cF^{\bullet\boxtimes k}}$ obtained by deleting the term in degree $-k$.
We have the following exact triangle in $D^b(\Sym^kC\times M_i)$:
\begin{equation}\label{adtqerghqerhq}
G_k\to{\cF^{\bullet\boxtimes k}}\to{\cF^{*\boxtimes k}}[k]\to.
\end{equation}
Furthermore, 
$\cP_{G_k}(X)\in\langle{\cF^{*\boxtimes k-1}},\ldots,\cF^*,\cO\rangle$
for every $X\in D^b(\Sym^kC)$.
\end{lemma}

\begin{proof}
The degree $-k$ term of $\mathcal F^{\bullet\boxtimes k}$ is $\mathcal F^{*\boxtimes k}$: the Koszul character on the
tensor product of the $k$ degree-$(-1)$ factors cancels the external sign
character. Hence the short exact sequence of complexes
\[
0\longrightarrow G_k\longrightarrow
\mathcal F^{\bullet\boxtimes k}\longrightarrow
\mathcal F^{*\boxtimes k}[k]\longrightarrow0
\]
induces the distinguished triangle~\eqref{adtqerghqerhq}.

 \begin{claim}\label{afsgarbar}
 Let $0\le l<k$ and let $r=k-l$.
 Let $p_{12}$, $p_{13}$ and $p_{23}$ be the projections
from
$\operatorname{Sym}^kC\times\operatorname{Sym}^lC\times M_i$.
The degree-$(-l)$ term of
$\mathcal F^{\bullet\boxtimes k}$ is
\begin{equation}\label{sergtergwerhqe5}
 \left(\cF^{\bullet\boxtimes k}\right)^{-l}=Rp_{13*}\!\left(
  Lp_{12}^* \mathcal O_{W_{k,l}}
  \otimes\rho_{k,l}^*\mathscr L_r
  \otimes Lp_{23}^*\mathcal F^{*\boxtimes l}
 \right),
 \end{equation}
 where $\mathscr L_r=\mathcal O_{\operatorname{Sym}^rC}(-\Delta_r/2)$
 (as defined in  Lemma~\ref{amazinggrace}.)
 \end{claim}
 
\begin{proof}[Proof of Claim~\ref{afsgarbar}]
Since $\mathcal F^\bullet=[\mathcal F^*\to\mathcal O]$, with
$\mathcal O$ in degree $0$, the degree-$(-l)$ term of ${\widehat\cF^{\bullet\boxtimes k}}$ is
\begin{equation}\label{asdvasfvasfgasrg}
\left(\widehat\cF^{\bullet\boxtimes k}\right)^{-l}=
 \bigoplus_{\substack{U\subset\{1,\ldots,k\}\\ |U|=l}}
 \bigotimes_{j\in U}\pi_j^*\mathcal F^*.
\end{equation}
The group $S_k$ permutes the summands.  Choose
$U_0=\{r+1,\ldots,k\}$.  Its stabilizer is the subgroup
$S_r\times S_l\subset S_k$ that permutes the first $r$ and the last $l$
coordinates separately.  On the chosen summand, permutation of the $l$
degree-$(-1)$ factors contributes the Koszul character
$\operatorname{sgn}_l$, whereas the external sign character restricts as
$
 \operatorname{sgn}_k|_{S_r\times S_l}
 =\operatorname{sgn}_r\boxtimes\operatorname{sgn}_l$.
Consequently the total character is
$ \operatorname{sgn}_l\,
 (\operatorname{sgn}_r\boxtimes\operatorname{sgn}_l)
 =\operatorname{sgn}_r$.
Equivalently, after tensoring by the external sign representation, the
entire degree-$(-l)$ term is the induced $S_k$-equivariant sheaf
\begin{equation}\label{eq:claim48-induced-term}
 \left(\widehat\cF^{\bullet\boxtimes k}\right)^{-l}\otimes\operatorname{sgn}_k
 \simeq
 \operatorname{Ind}_{S_r\times S_l}^{S_k}
 \left(
 ( \mathcal O_{C^r}\otimes\operatorname{sgn}_r)
  \boxtimes\widehat{\mathcal F}^{*\boxtimes l}
 \right).
\end{equation}
Frobenius reciprocity 
shows that
the first $r$ coordinates acquire the sign character and the last
$l$ coordinates do not. The equivariant pushforwards  
of $\mathcal O_{C^r}\otimes\operatorname{sgn}_r$ and $\widehat{\mathcal F}^{*\boxtimes l}$
are
$\mathscr L_r=\mathcal O_{\operatorname{Sym}^rC}(-\Delta_r/2)$ and
$\mathcal F^{*\boxtimes l}$, respectively.

The quotient of $C^k$ that remembers both the total divisor and the divisor
formed by the last $l$ coordinates is the incidence correspondence \eqref{asfvasfgasgarga}.
On this correspondence the complementary divisor is
$\rho_{k,l}(D,D')=D-D'$.  Applying the preceding descent and then pushing
along $p_{13}$ gives \eqref{sergtergwerhqe5}.
\end{proof}

For $0\le l<k$, Claim~\ref{afsgarbar} identifies the degree $-l$ term as a
convolution whose second Fourier--Mukai factor is
$\mathcal P_{\mathcal F^{*\boxtimes l}}$.  The finite complex $G_k$ has a
filtration by brutal truncations with these degree terms as successive
quotients.  Applying $\mathcal P_{(-)}(X)$ therefore gives
$
 \mathcal P_{G_k}(X)\in
 \langle\mathcal F^{*\boxtimes(k-1)},\ldots,
          \mathcal F^*,\mathcal O\rangle
$.
\end{proof}

\begin{lemma}\label{sRGwgwrG}
Let $0\le l\le k\le i$ and put $r=k-l$.
Let
$$
 \Gamma^W_{k,r}=\{(D,E)\in W_{k,r}:E\cap(D-E)\ne\varnothing\}
$$
with its effective Cartier divisor structure.
Set
\[
 \mathcal A_{k,r}=
 \mathcal O_{W_{k,r}}(\Gamma^W_{k,r})\otimes
 \left(\Lambda^{*\boxtimes k}\boxtimes
            (\mathscr L_r\otimes\Lambda^{\boxtimes r})\right).
\]
Here a line bundle on $W_{k,r}$ is also viewed as its pushforward to
$\Sym^kC\times\Sym^rC$.  There is an isomorphism of Fourier--Mukai
functors
\[
 \mathcal P_{\mathcal H^{-l}(\mathcal F^{\bullet\boxtimes k})}
 \simeq
 \mathcal P_{\mathcal D_i^r\otimes(\boldsymbol\Lambda^*)^{\otimes l}}
       \circ\mathcal P_{\mathcal A_{k,r}}.
\]
In particular, for every $X\in D^b(\Sym^kC)$ its image belongs to
$\langle\mathcal D_i^r\otimes(\boldsymbol\Lambda^*)^{\otimes l}\rangle$.
\end{lemma}

The lemma is proved in
Section~\ref{asrgarsgargaerhaehae}.
\begin{lemma}\label{sRHwrhjeRJAETJ}
We have the following  exact triangle in $D^b(\Sym^kC\times M_i)$:
\begin{equation}\label{adfbadfb}
 H_k\longrightarrow\mathcal F^{\bullet\boxtimes k}
 \longrightarrow\mathcal D_i^k(-\Delta_k/2)\longrightarrow
\end{equation}
with
\[
 \mathcal P_{H_k}(X)\in
 \bigl\langle
 \mathcal D_i^{k-1}\boldsymbol\Lambda^*,\ldots,
 \mathcal D_i^1(\boldsymbol\Lambda^*)^{k-1},
 (\boldsymbol\Lambda^*)^k
 \bigr\rangle
\]
for every $X\in D^b(\operatorname{Sym}^kC)$, where
$\Delta_k\subset\operatorname{Sym}^kC$ is the diagonal divisor.
\end{lemma}

\begin{proof}
Take $H_k=\tau_{\le-1}(\mathcal F^{\bullet\boxtimes k})$.
Lemma~\ref{sRGwgwrG} with $l=0$ identifies the last term of its
truncation triangle with $\mathcal D_i^k(-\Delta_k/2)$.
The successive quotients of the canonical filtration of $H_k$ are
$\mathcal H^{-l}(\mathcal F^{\bullet\boxtimes k})[l]$, $1\le l\le k$.
Their Fourier--Mukai images belong to
$\langle\mathcal D_i^{k-l}(\boldsymbol\Lambda^*)^l\rangle$
by the same lemma, proving the assertion.
\end{proof}

\begin{corollary}
Applying 
the Fourier--Mukai functors to the exact triangles \eqref{adtqerghqerhq}
and \eqref{adfbadfb}
gives exact triangles
\begin{equation}\label{1srgseGWEG}
\cP_{G_k}(X)\to\cP_{{\cF^{\bullet\boxtimes k}}}(X)\to\cP_{{\cF^{*\boxtimes k}}}(X)[k]\to
\end{equation}
and 
\begin{equation}\label{1Xfgsrhgasrha}
 \mathcal P_{H_k}(X)
 \longrightarrow\mathcal P_{\mathcal F^{\bullet\boxtimes k}}(X)
 \longrightarrow
 \mathcal P_{\mathcal D_i^k}\bigl(X(-\Delta_k/2)\bigr)
 \longrightarrow
\end{equation}
for every $X\in D^b(\Sym^kC)$.
Here $\cP_{G_k}(X)\in\langle{\cF^{*\boxtimes k-1}},\ldots,\cF^*,\cO\rangle$
and $\cP_{H_k}(X)\in\langle\cD_i^{k-1}\bLambda^*,\ldots,\cD_i^1{\bLambda^*}^{k-1},{\bLambda^*}^k\rangle$.
\end{corollary}

\medskip

\begin{figure}[htbp]
\includegraphics[width=\textwidth]{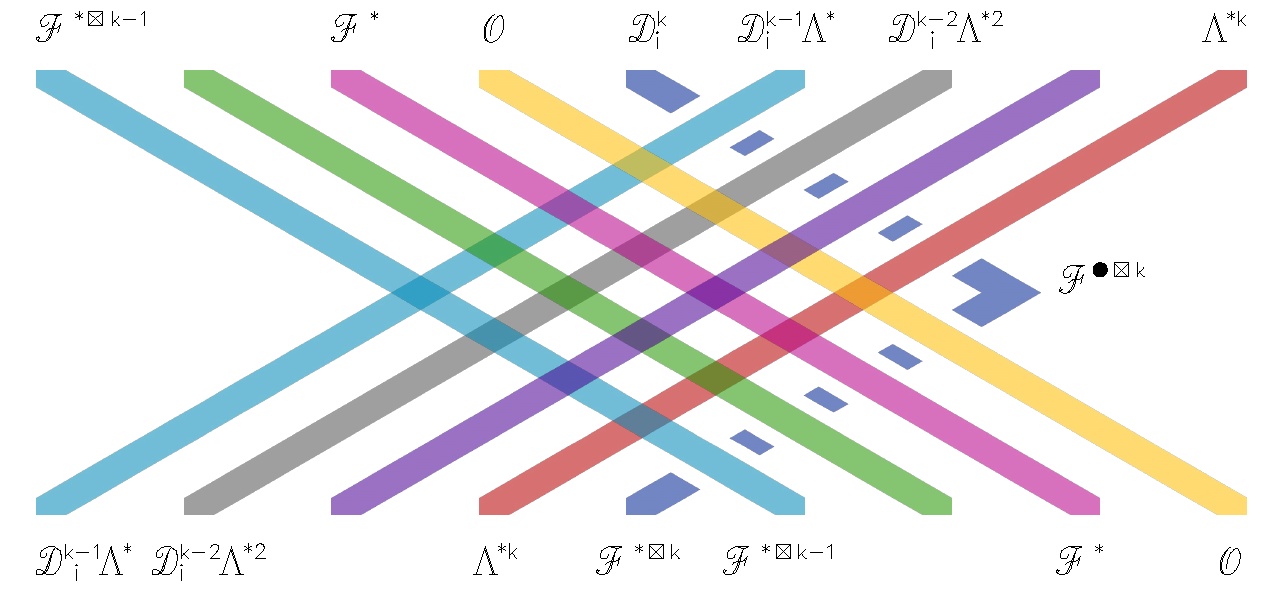}
\caption{The intermediate subcategory $\langle\cF^{\bullet\boxtimes k}\rangle$.}\label{ladyvanishes}
\end{figure}

To finish the proof of Theorem~\ref{sGSRHSRH}, we need Lemma~\ref{1adfbadfnbadn} and Lemma~\ref{1swgasrgaerh},
which break the required mutation in two steps, as illustrated in Figure~\ref{ladyvanishes}.

\begin{lemma}\label{1adfbadfnbadn}
$\cP_{{\cF^{\bullet\boxtimes k}}}$ is fully faithful and there is a mutation 
$$\langle {\cF^{*\boxtimes k}},\cF^{*\boxtimes k-1},\ldots,\cO\rangle\to\langle\cF^{*\boxtimes k-1},\ldots,\cO,{\cF^{\bullet\boxtimes k}}\rangle.$$
\end{lemma}

\begin{proof}
Use $\mathcal B_k$ from Remark~\ref{Remark33} and the convention \(\cF^{*\boxtimes0}=\cO\). We first prove that
\begin{equation}\label{1sdasbsRHr}
\RHom\bigl(
\cP_{\cF^{\bullet\boxtimes k}}(X),
\cP_{\cF^{*\boxtimes m}}(Y_m)
\bigr)=0
\end{equation}
for \(X\in D^b(\Sym^kC)\), \(Y_m\in D^b(\Sym^mC)\), and
\(0\le m<k\).

By~\eqref{1Xfgsrhgasrha} and the description of \(\cP_{H_k}(X)\), it is
enough to prove
\begin{equation}\label{eq:lemma414-block-vanishing}
\RHom\bigl(
\cP_{\cD_i^l{\bLambda^*}^{k-l}}(X_l),
\cP_{\cF^{*\boxtimes m}}(Y_m)
\bigr)=0
\end{equation}
for \(0\le l\le k\), \(X_l\in D^b(\Sym^lC)\), and
\(Y_m\in D^b(\Sym^mC)\). 
It suffices to take $X_l=\cO_D$ and $Y_m=\cO_{D'}$.
Put $Y=M_{i-l}(\Lambda(-2D))\subset M_i$. Adjunction gives
$\det N_{Y/M_i}\simeq\bLambda^l$, up to a fixed one-dimensional
factor. Coherent duality for regular immersion therefore gives
\[
 \RHom_{M_i}\bigl(\cO_Y\otimes\bLambda^{l-k},
                         \cF_{D'}^{*\boxtimes m}\bigr)
 \simeq
 R\Gamma\bigl(Y,\cF_{D'}^{*\boxtimes m}|_Y
                    \otimes\bLambda^k\bigr)[-2l],
\]
with the same convention on constant factors.
Apply Theorem~\ref{TT_Theorem7.4} with
$j=i-l$, $d_0=2g-1-2l$, $a=m$, $b=0$, and $t=k$:
\[
 m<k\le i<3g-2i-2=d_0+g-2j-1.
\]
For $j=0$ use its zeroth-chamber clause, which allows $d_0=1$.
This proves~\eqref{eq:lemma414-block-vanishing}
and~\eqref{1sdasbsRHr}.

 Rotate \eqref{1srgseGWEG} to obtain
$
 \mathcal P_{\mathcal F^{\bullet\boxtimes k}}(X)
 \longrightarrow
 \mathcal P_{\mathcal F^{*\boxtimes k}}(X)[k]
 \longrightarrow
 \mathcal P_{G_k}(X)[1]
 \longrightarrow
$.
Here $\mathcal P_{G_k}(X)[1]\in\mathcal B_k$, whereas
\eqref{1sdasbsRHr} gives
$\mathcal P_{\mathcal F^{\bullet\boxtimes k}}(X)
 \in{}^\perp\mathcal B_k$.  Therefore this is the canonical
right-mutation triangle and
\[
 \mathbf R_{\mathcal B_k}
 \bigl(\mathcal P_{\mathcal F^{*\boxtimes k}}(X)[k]\bigr)
 \simeq
 \mathcal P_{\mathcal F^{\bullet\boxtimes k}}(X).
\]

Right mutation is an equivalence on the component being mutated
\cite[Theorem~2.10]{KuzHH}. Since $\cP_{\cF^{*\boxtimes k}}$ is fully
faithful by Lemma~\ref{wRGwrhgwrHR}, the displayed natural isomorphism
also proves that $\cP_{\cF^{\bullet\boxtimes k}}$ is fully faithful.
\end{proof}

\begin{lemma}\label{1swgasrgaerh}
$\cP_{\cD_i^k}$ is fully faithful. There is a mutation 
$$\langle \cD_i^{k-1}\otimes\bLambda^*,\ldots,{\bLambda^*}^k,\cF^{\bullet\boxtimes k}\rangle
\to\langle\cD_i^k, \cD_i^{k-1}\otimes\bLambda^*,\ldots,{\bLambda^*}^k\rangle.$$
\end{lemma}

\begin{proof}
We argue as part of the simultaneous induction on $k$ used in the proof of
Theorem~\ref{sGSRHSRH}. Thus all parts of that theorem are available in
indices $<k$, and Lemma~\ref{1adfbadfnbadn} has already been proved in
index $k$. Use $\mathcal A_k$ from Remark~\ref{Remark33} and the
convention $\cD_i^0=\cO_{M_i}$.
By the inductive Cross Warp in degree $k-1$, followed by tensoring with
$\bLambda^*$, $\mathcal A_k$ is an admissible 
subcategory of $D^b(M_i)$.
We first establish an auxiliary vanishing. For
$
0\le r\le k$,
$0\le \ell<k$,
$X_r\in D^b(\Sym^rC)$,
$Y_\ell\in D^b(\Sym^\ell C)$,
we claim that
\[
\RHom\Bigl(
 \cP_{\cF^{*\boxtimes r}}(X_r),
 \cP_{\cD_i^{\ell}\otimes{\bLambda^*}^{k-\ell}}(Y_\ell)
\Bigr)=0.
\]
Both essential images in the preceding display are admissible: this follows
from Lemma~\ref{wRGwrhgwrHR} for the first one and from the inductive
hypothesis for the second one. Hence, by
\cite[Lemma~10.2]{TT}, it is enough to use the spanning classes of
skyscraper sheaves. Take
$X_r=\cO_D$ and $Y_\ell=\cO_{D'}$, where
$D\in\Sym^rC$ and $D'\in\Sym^\ell C$. Then
\begin{align*}
&\RHom_{M_i}\Bigl(
 \cF_D^{*\boxtimes r},
 \cO_{M_{i-\ell}(\Lambda(-2D'))}
       \otimes\bLambda^{\ell-k}
\Bigr)\\
&\qquad\simeq
R\Gamma\Bigl(
 M_{i-\ell}(\Lambda(-2D')),
 (\cF_D^{*\boxtimes r})^*\otimes\bLambda^{\ell-k}
\Bigr)\\
&\qquad\simeq
R\Gamma\Bigl(
 M_{i-\ell}(\Lambda(-2D')),
 \bar\cF_D^{\boxtimes r}\otimes\bLambda^{\ell-k}
\Bigr).
\end{align*}
In the last line we suppress the one-dimensional factor obtained by
restricting
$\mathcal O_{\operatorname{Sym}^rC}(\Delta_r/2)$ to the point
$D\in\operatorname{Sym}^rC$; it has no effect on the vanishing.
Apply Theorem~\ref{hardvanishing} with
$
 j=i-\ell$,
 $d_0=2g-1-2\ell$,
 $a=0$,
 $b=r$,
 $t=\ell-k$.
Here $j\ge1$, $d_0\ge3$, and
\[
 -j-1=\ell-i-1<\ell-k<0<3g-2i-2-r
                  =d_0+g-2j-1-r.
\]
Thus Theorem~\ref{hardvanishing} gives the auxiliary vanishing.

We now prove
\begin{equation}\label{1rgarharreh}
\RHom\Bigl(
 \cP_{\cF^{\bullet\boxtimes k}}(X),
 \cP_{\cD_i^{\ell}\otimes{\bLambda^*}^{k-\ell}}(Y)
\Bigr)=0
\end{equation}
for $0\le\ell<k$, $X\in D^b(\Sym^kC)$, and
$Y\in D^b(\Sym^\ell C)$. Apply the functor 
$\RHom(-,\cP_{\cD_i^{\ell}\otimes{\bLambda^*}^{k-\ell}}(Y))$
to the triangle~\eqref{1srgseGWEG}. Its first term lies in
$\langle\cF^{*\boxtimes k-1},\ldots,\cF^*,\cO\rangle$, while its third
term is $\cP_{\cF^{*\boxtimes k}}(X)[k]$. Both outer terms vanish by the auxiliary vanishing; hence~\eqref{1rgarharreh} follows.

We next prove the semi-orthogonality in the other direction. The
inductive Cross Warp implies
\[
\langle\cD_i^\ell\rangle
\subset
\left\langle
 \cF^{*\boxtimes b}\otimes{\bLambda^*}^{q}
 \ \middle|\
 0\le b\le\ell,\ 0\le q\le\ell-b
\right\rangle
\qquad (\ell<k).
\]

After tensoring this containment by
${\bLambda^*}^{k-\ell}$, every generator on the right has the form
$
\cF^{*\boxtimes b}\otimes{\bLambda^*}^{t}$,
$0\le b<k$,
$0<t\le k-b$,
where $t=q+k-\ell$. Let $D\in\Sym^kC$ and
$D'\in\Sym^bC$. We have
\begin{align*}
&\RHom_{M_i}\Bigl(
 \cF_{D'}^{*\boxtimes b}\otimes{\bLambda^*}^{t},
 \cO_{M_{i-k}(\Lambda(-2D))}
\Bigr)\\
&\qquad\simeq
R\Gamma\Bigl(
 M_{i-k}(\Lambda(-2D)),
 (\cF_{D'}^{*\boxtimes b})^*\otimes\bLambda^t
\Bigr)\\
&\qquad\simeq
R\Gamma\Bigl(
 M_{i-k}(\Lambda(-2D)),
 \bar\cF_{D'}^{\boxtimes b}\otimes\bLambda^t
\Bigr).
\end{align*}
Here the one-dimensional factor $(\mathscr L_b^*)_{D'}$ from
Lemma~\ref{amazinggrace} is suppressed.

 Apply
Theorem~\ref{TT_Theorem7.4} with
$
 j_1=i-k$,
 $d_1=2g-1-2k$,
 $a=0$,
 $b=b$,
 $t=t$.
Its hypotheses follow from
\[
 0<t\le k-b<3g-2i-2-b=d_1+g-2j_1-1-b,
\]
using the zeroth-chamber clause when $j_1=0$.

The spanning-class argument at the beginning of this section therefore
extends the pointwise vanishing to
\begin{equation}\label{1sdgwetwethw}
 \RHom\bigl(\mathcal A_k,\cP_{\cD_i^k}(X)\bigr)=0
 \qquad\text{for every }X\in D^b(\Sym^kC).
\end{equation}
This step uses the adjoints of the Fourier--Mukai functors, not full
faithfulness of $\cP_{\cD_i^k}$.
It remains to identify the mutation functorially. The triangle
\eqref{1Xfgsrhgasrha} is natural in $X$; its first term belongs to
$\mathcal A_k$, while its third term belongs to $\mathcal A_k^\perp$ by
\eqref{1sdgwetwethw}. Hence it is the canonical left-mutation triangle and
gives a natural isomorphism
\[
 \mathbf L_{\mathcal A_k}\circ
 \mathcal P_{\mathcal F^{\bullet\boxtimes k}}
 \simeq
 \mathcal P_{\mathcal D_i^k}\circ
 \bigl(-\otimes
 \mathcal O_{\operatorname{Sym}^kC}(-\Delta_k/2)\bigr).
\]

By~\eqref{1rgarharreh},
$\langle\mathcal A_k,\cF^{\bullet\boxtimes k}\rangle$ is a
semiorthogonal sequence of admissible subcategories. Left mutation is
an equivalence on its second component \cite[Theorem~2.10]{KuzHH}.
The displayed natural isomorphism and
Lemma~\ref{1adfbadfnbadn} therefore prove full faithfulness of
$\cP_{\cD_i^k}$. The source twist is an autoequivalence, so its
essential image is $\langle\cD_i^k\rangle$.
Equation
\eqref{1rgarharreh} gives the semi-orthogonality of the collection before
mutation, and~\eqref{1sdgwetwethw} gives that of the collection after
mutation. Therefore
$
\bigl\langle\mathcal A_k,\cF^{\bullet\boxtimes k}\bigr\rangle
\longrightarrow
\bigl\langle\cD_i^k,\mathcal A_k\bigr\rangle
$
is the claimed mutation.
\end{proof}

\begin{lemma}
The subcategories $\mathcal A_k$ and $\mathcal B_k$ of
Remark~\ref{Remark33} are completely orthogonal:
$
\RHom(\mathcal A_k,\mathcal B_k)=0$ and 
$
\RHom(\mathcal B_k,\mathcal A_k)=0$.
\end{lemma}

\begin{proof}
Indeed, the first vanishing is exactly
\eqref{eq:lemma414-block-vanishing} with $0\le l<k$ and $0\le m<k$.
The second is the auxiliary vanishing proved at the beginning of
Lemma~\ref{1swgasrgaerh}, with $0\le r<k$ and $0\le\ell<k$.
\end{proof}

\begin{proof}[Proof of Theorem~\ref{sGSRHSRH}]
Part~(a) is Lemma~\ref{1swgasrgaerh}.  Part~(c) follows from
Lemmas~\ref{wRGwrhgwrHR}, \ref{1adfbadfnbadn}, and
\ref{1swgasrgaerh}, together with the induction hypothesis; the
orthogonality proved immediately above permits the indicated
reordering without changing the embeddings.
For part~(b), argue by induction on $k$, including $k=0$, where the
kernel is $\mathcal O$.  Strengthen the induction assertion by requiring
the common lift of $\mathcal D^k$ to have weights $0,\ldots,k$.
Fix $k\le i-1$.  Use the window and restrictions of
\eqref{windowsembedding}, denoted locally by
\[
 r_-:\mathcal G_{[0,i)}\xrightarrow{\sim}D^b(M_{i-1}),
 \qquad
 r_+:\mathcal G_{[0,i)}\hookrightarrow D^b(M_i).
\]
By \cite[Theorem~3.21]{TT}, $\mathcal F^*$ has wall weights $0,1$,
$\mathcal O$ has weight $0$, and $\boldsymbol\Lambda^*$ has weight $1$.
Thus the same equivariant kernel $\mathcal F^{\bullet\boxtimes k}$
defines a functor into this window, with weights $0,\ldots,k$.
For $0\le r<k$, the induction hypothesis gives a common lift of
$\mathcal D^r(\boldsymbol\Lambda^*)^{k-r}$ with weights
$k-r,\ldots,k$.  These, and only these lower-degree incidence kernels,
are used to lift $\mathcal A_k$ from Remark~\ref{Remark33}.

Form the left mutation of the lifted
$\mathcal P_{\mathcal F^{\bullet\boxtimes k}}$ through that lifted
subcategory, inside the window.  The projection term and its cone have
weights in $0,\ldots,k$, since this weight range is closed under
extensions.  Full faithfulness of both restriction functors identifies
the evaluation morphism with the corresponding evaluation morphism in
each chamber.  Consequently the two restrictions of this mutation
are, by~\eqref{1Xfgsrhgasrha},
$
 \mathcal P_{\mathcal D_{i-1}^k}\circ(-\otimes\mathscr L_k)$,
$\mathcal P_{\mathcal D_i^k}\circ(-\otimes\mathscr L_k)$.
The source twist is the same in both chambers and has wall weight zero.
Cancel it.  Since $\iota_i=r_+r_-^{-1}$, this gives the natural
isomorphism
$
 \iota_i\circ\mathcal P_{\mathcal D_{i-1}^k}
       \simeq\mathcal P_{\mathcal D_i^k}$,
and also constructs the required common lift with weights
$0,\ldots,k$.  
\end{proof}

\section{Cohomology of the complex ${\cF^{\bullet\boxtimes k}}$}\label{asrgarsgargaerhaehae}

Fix $r,l\ge0$ with $r+l=k\le i$. Write
\[
 j=j_{r,l}:Z_{r,l}=D_i^r\times\Sym^lC
 \hookrightarrow Y_{r,l}=\Sym^rC\times\Sym^lC\times M_i
\]
for the closed immersion, and let
\[
 a=a_{r,l}:Y_{r,l}\longrightarrow\Sym^kC\times M_i,
 \qquad(E,U,F,s)\longmapsto(E+U,F,s).
\]
Let $\Gamma_{r,l}\subset Y_{r,l}$ be the effective Cartier divisor
where $E$ and $U$ have a common point. We take $\Gamma_{r,l}=0$
when $r=0$ or $l=0$.
The relative dualizing sheaf of $a$ is
$\omega_a\simeq\mathcal O_{Y_{r,l}}(\Gamma_{r,l})$, as proved in
Lemma~\ref{lem:addition-duality}.
Define the following line bundle on $Z_{r,l}$ and the cohomology sheaf on $\Sym^kC\times M_i$,
\begin{equation}\label{eq:corrected-incidence-lines}
 \mathscr M_{r,l}
 =j^*\bigl(\mathscr L_r\boxtimes\Lambda^{*\boxtimes l}
                     \boxtimes(\boldsymbol\Lambda^*)^{\otimes l}\bigr),
 \qquad
 M_{r,l}=\mathcal H^{-l}(\mathcal F^{\bullet\boxtimes k}).
\end{equation}
Let $\beta=\beta_{r,l}=a\circ j$.
Let $\Sigma_r=\beta_{r,l}(Z_{r,l})\subset {\Sym^kC\times M_i}$ with its reduced structure.
Its points are the pairs $(D,F,s)$ such that $D$ contains a divisor
$E$ of degree $r$ with $s|_E=0$.  We refer to $E$ as the {\em selected}
divisor when discussing a point of $Z_{r,l}$
and to $U=D-E$ as the {\em residual}
divisor.

The smooth projective morphism $D_i^r\to\Sym^rC$ has nonempty
irreducible fibers $M_{i-r}(\Lambda(-2E))$
\cite[Remark~3.7]{TT}, \cite[(3.1), (3.6)]{thaddeus}.
Thus $D_i^r$ and $Z_{r,l}=D_i^r\times\Sym^lC$ are smooth and
irreducible, and $\Sigma_r$ is integral.

The goal of this section is to prove the following theorem.

\begin{theorem}\label{sdfgnsfghsftjt}
For $r+l=k\le i$, there is a natural isomorphism
\begin{equation}\label{eq:corrected-cohomology}
 \Phi^+_{r,l}:\beta_{r,l*}\bigl(\mathscr M_{r,l}\otimes j^*\omega_{a_{r,l}}\bigr)
 \xrightarrow{\sim}
 \mathcal H^{-l}(\mathcal F^{\bullet\boxtimes k})=M_{r,l}
\end{equation}
of sheaves on $\Sym^kC\times M_i$.
Thus $M_{r,l}$ is a Cohen--Macaulay sheaf of pure codimension $2r$, and its
scheme-theoretic support is the reduced subscheme~$\Sigma_r$.
The finite morphism
$\beta_{r,l}:Z_{r,l}\to\Sigma_r$ is the
normalization.  As an $\mathcal O_{\Sigma_r}$-module, $M_{r,l}$ is
torsion-free of generic rank one.
\end{theorem}

\begin{proof}[Proof of Lemma~\ref{sRGwgwrG}]
The subscheme $W_{k,r}$ appearing in 
Lemma~\ref{sRGwgwrG} is given by $W_{k,r}:=
 \{(D,E)\in\operatorname{Sym}^kC\times\operatorname{Sym}^rC
       \mid D\ge E\}$.
Use the isomorphism
\[
 W_{k,r}\simeq\Sym^rC\times\Sym^lC,\qquad
 (D,E)\longmapsto(E,D-E),
\]
so that its second projection is the difference map
$\rho_{k,r}$ of Definition~\ref{qetgathaehehe}.
The projection $W_{k,r}\to \Sym^rC$ is smooth.  Its fiber product with
$D_i^r\to\Sym^rC$ is $Z_{r,l}$ and there are no higher Tor sheaves
in the convolution of the two kernels.  Tensoring with the line
bundle $\mathcal O_{W_{k,r}}(\Gamma^W_{k,r})$ does not change this
Tor-independence.  Multiplicativity of the tensor line bundles gives
$\rho_{k,r}^*(\Lambda^{*\boxtimes l})\simeq
 \left(\Lambda^{*\boxtimes k}\boxtimes\Lambda^{\boxtimes r}\right)
       |_{W_{k,r}}$.
The pullback of $\Gamma^W_{k,r}$ to $Z_{r,l}$ is $j^*\Gamma_{r,l}$.
Thus the tensor product of the two pulled-back kernels in Lemma~\ref{sRGwgwrG} is
$\mathscr M_{r,l}\otimes j^*\omega_a$.
The remaining projection is the finite map $\beta_{r,l}$, and the
convolution is isomorphic to 
$\beta_{r,l*}(\mathscr M_{r,l}\otimes j^*\omega_a)$, which by Theorem~\ref{sdfgnsfghsftjt} is isomorphic
 to $M_{r,l}$.
\end{proof}

\begin{lemma}\label{lem:CMcandidate}
The subscheme $\Sigma_r$ has pure codimension $2r$ in ${\Sym^kC\times M_i}$.
For every line bundle $L$ on $Z_{r,l}$,
the sheaf $\beta_{r,l*}L$
is Cohen--Macaulay of pure codimension $2r$, with scheme-theoretic
support $\Sigma_r$.
\end{lemma}

\begin{proof}
Since $\dim Z_{r,l}=3g-3+k-2r$ and $\beta_{r,l}$ is finite,
$\Sigma_r$ has codimension $2r$. Finite pushforward of a line
bundle on the smooth variety $Z_{r,l}$ is Cohen--Macaulay.
Its annihilator is
$\ker(\mathcal O_{\Sym^kC\times M_i}\to\beta_{r,l*}\mathcal O_{Z_{r,l}})$,
which is the ideal sheaf of the reduced image $\Sigma_r$.
\end{proof}

\begin{lemma}\label{lem:altbound}\label{prop:support}
\begin{equation}\label{eq:support}
 \Supp\cH^{-l}({\cF^{\bullet\boxtimes k}})\subseteq\Sigma_r.
\end{equation}
\end{lemma}

\begin{proof}
Fix $x=(D,F,s)$ and write
\[
 D=\sum_\alpha m_\alpha p_\alpha,\qquad
 b_\alpha=\operatorname{ord}_{p_\alpha}(s),\qquad
 \nu_x=\sum_\alpha\min(m_\alpha,b_\alpha).
\]
Then $x\in\Sigma_r$ exactly when $\nu_x\ge r$.
We compute first on the fiber over $(F,s)$. Separate the support
points \'etale-locally and choose a uniformizer $z_\alpha$ at each
$p_\alpha$. In a suitable frame $s=(z_\alpha^{b_\alpha},0)$, so
its evaluation complex splits as
$[\mathcal O\xrightarrow{z_\alpha^{b_\alpha}}\mathcal O]
       \oplus\mathcal O[1]$.
For the $m_\alpha$ labels near $p_\alpha$, choose $q_\alpha$
factors of the first type. The K\"unneth formula places the
cohomology of this summand in degree $-k+\sum_\alpha q_\alpha$.
If $b_\alpha=0$ and $q_\alpha>0$, the summand is acyclic.

The stabilizer of the choice of factors is
$\prod_\alpha(S_{q_\alpha}\times S_{m_\alpha-q_\alpha})$.
The permutation sign on the factors $\mathcal O[1]$ cancels the
corresponding part of the external sign character. Frobenius
reciprocity therefore leaves ordinary invariants on those factors
and alternating invariants on the $q_\alpha$ Koszul factors.
The latter contribute
\[
 \left(\bigl(\mathbb C[z]/(z^{b_\alpha})\bigr)^{\otimes q_\alpha}
               \otimes\operatorname{sgn}_{q_\alpha}\right)^{S_{q_\alpha}}
 \simeq\bigwedge^{q_\alpha}\bigl(\mathbb C[z]/(z^{b_\alpha})\bigr).
\]
Thus a nonzero contribution requires
$q_\alpha\le\min(m_\alpha,b_\alpha)$ for every $\alpha$.
The restricted complex has no cohomology above $-k+\nu_x$.
Represent the universal complex at $x$ by a finite free complex
$P^\bullet$ over $A=\mathcal O_{\Sym^kC\times M_i,x}$, and put
$I=\mathfrak m_{M_i,(F,s)}A$. If $n$ is the highest nonzero
cohomological degree of $P^\bullet$, then
$
 H^n(P^\bullet\otimes_A A/I)=H^n(P^\bullet)\otimes_A A/I
$.
The Nakayama lemma and the fiber calculation give $n\le-k+\nu_x$.
For $x\notin\Sigma_r$ we have $\nu_x<r$, so the cohomology in
degree $-l=-k+r$ vanishes, as required.
\end{proof}

\begin{definition}
Let $B_{\mathrm{mult}}\subset Z_{r,l}$ be the locus where $E+U$ has a
point of multiplicity at least three or two distinct multiple points.
Let $B_{\mathrm{extra}}$ be the locus where $s|_{E+p}=0$ for some
$p\in\operatorname{Supp}U$, and put
\[
 B_{r,l}=B_{\mathrm{mult}}\cup B_{\mathrm{extra}},\quad
 U_{r,l}=(\Sym^kC\times M_i)\setminus\beta_{r,l}(B_{r,l}),\quad
 V_{r,l}=U_{r,l}\cap\Sigma_r.
\]
Write $V_{r,l}^{\mathrm{red}}$ for the open subset where the total
divisor is reduced.
\end{definition}

Both excluded loci are closed of codimension at least two in $Z_{r,l}$.
For $B_{\mathrm{mult}}$, this follows from the smooth projection
$Z_{r,l}\to\Sym^rC\times\Sym^lC$ and the finite addition map: the
locus of two collisions has codimension at least two in $\Sym^kC$.
For $l>0$, the locus $B_{\mathrm{extra}}$ is the finite image of
$
 \left(D_i^{r+1}\times_{\Sym^{r+1}C}
       (\Sym^rC\times C)\right)\times\Sym^{l-1}C$,
under $(E,p,F,s,U')\mapsto(E,p+U',F,s)$; its source has dimension
$\dim Z_{r,l}-2$. Here $r+1\le k\le i$. For $l=0$ this locus is
empty. Finiteness of $\beta_{r,l}$ now gives
\begin{equation}\label{eq:bigopen}
 \operatorname{codim}_{\Sigma_r}(\Sigma_r\setminus U_{r,l})\ge2.
\end{equation}
The locus $V_{r,l}^{\mathrm{red}}$ is dense in $\Sigma_r$, by
smoothness over the product of symmetric powers and the
codimension of the excluded loci. Every point of $V_{r,l}$ has a unique selected divisor:
if $E'\ne E$ were another one, then $s|_{E+p}=0$ for a point
$p\in\Supp(D-E)$, contrary to the definition of $B_{\mathrm{extra}}$.

\begin{lemma}\label{lem:addition-duality}
The morphism $a:Y_{r,l}\longrightarrow {\Sym^kC\times M_i}$ is finite flat and Gorenstein, and is \'etale away
from $\Gamma_{r,l}$.  There is a canonical isomorphism
\begin{equation}\label{eq:omega-addition}
 \omega_a\simeq\mathcal O_{Y_{r,l}}(\Gamma_{r,l})
\end{equation}
under which the trace section
$\tau_a:\mathcal O_{Y_{r,l}}\to\omega_a$ is the natural inclusion
$\mathcal O_{Y_{r,l}}\hookrightarrow\mathcal O_{Y_{r,l}}(\Gamma_{r,l})$.
Moreover,
\begin{equation}\label{eq:discriminant-factor}
 a^*\Delta_k=\operatorname{pr}_1^*\Delta_r+
             \operatorname{pr}_2^*\Delta_l+2\Gamma_{r,l},
\end{equation}
and
\begin{equation}\label{eq:signline-addition}
 a^*\scrL_k\simeq\operatorname{pr}_1^*\scrL_r\otimes
                    \operatorname{pr}_2^*\scrL_l\otimes
                    \mathcal I_{\Gamma_{r,l}}.
\end{equation}
\end{lemma}

\begin{proof}
The addition map is a finite morphism between smooth varieties of the same
dimension, hence flat and Gorenstein. It is \'etale when the two
divisors are disjoint. Near a general point of $\Gamma_{r,l}$,
it is the quadratic cover ordering the common double point.
Write its algebra as $B=A\oplus A\xi$, with involution
$\xi\mapsto-\xi$. The functional
\[
 \lambda(b_0+b_1\xi)=2b_1
\]
generates $\operatorname{Hom}_A(B,A)$, and
$\operatorname{Tr}_{B/A}=\xi\lambda$.
Thus the trace section vanishes to order one on $\Gamma_{r,l}$,
proving~\eqref{eq:omega-addition} with its stated normalization.

On the ordered product, the sum of all diagonals is the sum of the
selected, residual, and mixed diagonals. Taking alternating
summands for $S_r\times S_l$ gives~\eqref{eq:signline-addition};
squaring gives~\eqref{eq:discriminant-factor}.
\end{proof}

\begin{lemma}\label{lem:selected-vandermonde-map}
Let $\widehat\Delta_r=\sum\limits_{1\le a<a'\le r}
\Delta_{aa'}$.
There is an $S_r$-equivariant morphism
\begin{equation}\label{eq:selected-vandermonde-map}
 \vartheta:\widehat{\mathcal D}_i^r(-\widehat\Delta_r)
 \longrightarrow
 \mathcal H^0(\widehat{\mathcal F}^{\bullet\boxtimes r}),
\end{equation}
where the source has its natural ideal-sheaf $S_r$-linearization.
Moreover, the morphism $\vartheta$ becomes an isomorphism after alternating direct
image:
\begin{equation}\label{asfarfsfweADJFHBAJDHRFB}
 \mathcal D_i^r\otimes\mathscr L_r
 \xrightarrow{\sim}\mathcal H^0(\mathcal F^{\bullet\boxtimes r}).
\end{equation}
\end{lemma}

\begin{proof}
The assertion is clear for $r=0$. For $r>0$, let
$\mathcal I=\mathcal I_{\widehat D_i^r}$ and let $\mathcal J$ be the
image of the differential into degree zero.
Locally choose a coordinate $t$ separating the support points, so
that $1,t,\ldots,t^{r-1}$ trivialize the algebra of the universal
divisor. In a frame of $\mathcal F$ along this divisor, write the
section as
\[
 \left(\sum_{h=0}^{r-1}f_ht^h,\ \sum_{h=0}^{r-1}g_ht^h\right).
\]
The coefficients $f_h,g_h$ generate $\mathcal I$, whereas their evaluations
at the ordered points generate $\mathcal J$. For the Vandermonde
matrix $V=(t_a^h)$ and
$\delta=\det V=\prod_{a<a'}(t_{a'}-t_a)$, the adjoint matrix identity gives
$
 \mathcal I_{\widehat\Delta_r}\mathcal I
       =\delta\mathcal I\subseteq\mathcal J\subseteq\mathcal I
$.
Hence the inclusion $\mathcal I_{\widehat\Delta_r}\hookrightarrow\mathcal O$
induces~\eqref{eq:selected-vandermonde-map}. This construction is
independent of the coordinate and frame and is equivariant for
the natural ideal-sheaf linearization.
The map
$\operatorname{Alt}_r(\mathcal I_{\widehat\Delta_r})
 \to\operatorname{Alt}_r(\mathcal O)$ is the identity on
$\mathscr L_r$. Since $\mathcal I$ is pulled back from
$\mathcal I_{D_i^r}$, the projection formula identifies the
alternating images of the two outer ideals in the displayed
inclusions with $\mathcal I_{D_i^r}\otimes\mathscr L_r$.
Exactness of alternating direct image now gives
\eqref{asfarfsfweADJFHBAJDHRFB}.
\end{proof}

Write
$Y_U=a^{-1}(U_{r,l})\subset Y_{r,l}$,
$Z_U=\beta_{r,l}^{-1}(U_{r,l})$,
$W=Y_U\setminus\Gamma_{r,l}$.

\begin{lemma}
There is a natural morphism
\begin{equation}\label{eq:kappa-fixed-splitting}
 \alpha^W_{r,l}:j_*\mathscr M_{r,l}|_W\longrightarrow a^*M_{r,l}|_W.
\end{equation}
\end{lemma}

\begin{proof}
Let $p:C^r\times C^l\times M_i\to Y_{r,l}$ be the quotient by
$H=S_r\times S_l$.
On $W$, the selected and residual divisors have disjoint supports.
The pullback of the full ordered quotient is the disjoint union of the
translates of $p^{-1}(W)$, indexed by $S_k/H$.  Finite base change and
Frobenius reciprocity give the canonical identification
\begin{equation}\label{eq:fixed-splitting-target}
 a^*M_{r,l}|_W\simeq
 \left(p_*\bigl(\mathcal H^{-l}
       (\widehat{\mathcal F}^{\bullet\boxtimes k})
             \otimes\operatorname{sgn}_k\bigr)\right)^H\big|_W.
\end{equation}

Write $K_b=\pi_b^*\mathcal F^\bullet$ and
$D_b=\det(\pi_b^*\mathcal F^*)$.
Contraction with the section gives a map of complexes
\begin{equation}\label{eq:prop512-kernel-map}
 D_b[1]\longrightarrow K_b,\qquad
 v\wedge w\longmapsto\Sigma(p_b)(v)w-\Sigma(p_b)(w)v.
\end{equation}
For $b=r+1,\ldots,k$, tensor these maps with
$K_1\otimes\cdots\otimes K_r$ and take cohomology in degree $-l$.
Precomposing with~\eqref{eq:selected-vandermonde-map} gives the map
\[
 \widehat{\mathcal D}_i^r(-\widehat\Delta_r)
       \otimes\bigotimes_{b=r+1}^kD_b
 \longrightarrow
 \mathcal H^{-l}(\widehat{\mathcal F}^{\bullet\boxtimes k}).
\]

Tensor with $\operatorname{sgn}_k$ and take $H$-invariants.
The permutations of the $l$ shifted determinant factors supply
$\operatorname{sgn}_l$, which cancels the external residual sign.
On the selected factor the inclusion
$\mathcal O(-\widehat\Delta_r)\subset\mathcal O_{C^r}$ induces an
isomorphism on alternating summands: every alternating regular function
vanishes on every transposition diagonal.  Its alternating descent is
therefore $\scrL_r$.  Projection formula gives the same statement after
pullback to $D_i^r$.  Thus the source descends to $j_*\mathscr M_{r,l}$.
\end{proof}

Next, we apply coherent duality, which gives
\begin{align}
 \Hom_{Y_U}(j_*\mathscr M_{r,l},a^*M_{r,l})
 &\simeq
 \Hom_{Y_U}(j_*\mathscr M_{r,l}\otimes\omega_a,a^!M_{r,l})\notag\\
 &\simeq
 \Hom_{U_{r,l}}\!\left(
 \beta_{r,l*}(\mathscr M_{r,l}\otimes j^*\omega_a)|_{U_{r,l}},
 M_{r,l}|_{U_{r,l}}\right).
 \label{eq:finite-duality-Hom}
\end{align}

\begin{proposition}\label{prop:split-comparison}
Suppose $M_{r,l}=\mathcal H^{-l}(\mathcal F^{\bullet\boxtimes k})$
is pure of codimension~$2r$.  The morphism
\eqref{eq:kappa-fixed-splitting} has a unique extension
\begin{equation}\label{eq:split-comparison}
 \widetilde\Phi^U_{r,l}:j_*\mathscr M_{r,l}
                      \longrightarrow a^*M_{r,l}
 \quad\text{on }Y_U.
 \end{equation}
Under coherent duality ~\eqref{eq:finite-duality-Hom} it corresponds to an isomorphism
\begin{equation}\label{eq:Phi-U}
 \Phi^{+,U}_{r,l}:
 \beta_{r,l*}(\mathscr M_{r,l}\otimes j^*\omega_a)|_{U_{r,l}}
                       \xrightarrow{\sim}M_{r,l}|_{U_{r,l}}.
\end{equation}
\end{proposition}

\begin{proof}
Since $a$ is finite flat, $a^*M_{r,l}$ is pure of dimension
$\dim Z_{r,l}$. The difference of two extensions has
image supported on $Z_U\cap\Gamma_{r,l}$, of smaller dimension,
and is therefore zero. The same argument proves uniqueness after
\'etale base change.

\smallskip
\noindent\textit{The mixed double point.}
Fix a closed point $y=(E,U,F,s)\in Z_U\cap\Gamma_{r,l}$
and put $x=a(y)$.
By the definition of $U_{r,l}$,
$
 E+U=2p+p_3+\cdots+p_k$,
where the displayed  points are distinct, and $p$ has
multiplicity one in each of $E,U$. The zero of $s$ at $p$ is
simple, since otherwise $s|_{E+p}=0$.

Let $V\subset\Sym^2C\times C^{k-2}\times M_i$ be the open subset
on which the labeled points $q_3,\ldots,q_k$ are pairwise distinct
and disjoint from the unordered degree-two divisor $D_2$.
The morphism
\[
 \nu:V\longrightarrow\Sym^kC\times M_i,\quad
 (D_2,q_3,\ldots,q_k,F,s)\longmapsto
 (D_2+q_3+\cdots+q_k,F,s)
\]
is \'etale. Set $v_0=(2p,p_3,\ldots,p_k,F,s)$ and
$A=\mathcal O_{V,v_0}$. Fix the sets $E',U'\subset\{3,\ldots,k\}$
of the other selected and residual labels.
The base change of $a$ to $V$ has an open-and-closed part where
$E$ has degree one on $D_2$ and contains precisely the labels $E'$.
It is $V\times_{\Sym^2C}(C\times C)$. Over $\Spec A$, write this part
as $\Spec B$, so
$
 \Spec B=(C\times C)\times_{\Sym^2C}\Spec A$.
On this cover write $D_2=q_1+q_2$, with $q_1$ selected and $q_2$
residual. Thus $B$ is finite free of rank two over $A$, with involution
$\sigma$ interchanging $q_1,q_2$. Its closed fiber is supported at
$(p,p)$, so $B$ is local.
By uniqueness of $E$, the support of $j_*\mathscr M_{r,l}$ on every
other open-and-closed part has empty fiber over the closed point of
$\Spec A$. Its image is closed by finiteness, so this support is
empty. Those parts consequently make no contribution to the adjoint.

Choose a regular parameter $t$ at $p$ and a frame of the universal
bundle near $(p,F,s)$. In the $A$-basis $1,t$ of $B$, write the section as
\[
 f=f_0+f_1t,\qquad g=g_0+g_1t,\qquad f_h,g_h\in A.
\]
The simple-zero condition makes one of $f_1,g_1$ invertible.
After interchanging the components, assume $f_1$ is a unit.
Replace $g$ by $g-(g_1/f_1)f$ in both evaluated factors and set
$x_1=f$, $x_2=\sigma(f)$,
$z=g_0-(g_1/f_1)f_0\in A$.
The~two evaluations are now $(x_1,z)$ and $(x_2,z)$, with the same
$z\in A$. Since $1,x_1$ is an $A$-basis of $B$, we have
\begin{equation}\label{eq:prop512-quadratic-algebra}
 \begin{gathered}
 e_1=x_1+x_2,\qquad e_2=x_1x_2,\qquad \xi=x_2-x_1,\qquad
 \xi^2=e_1^2-4e_2,\\
 B\simeq A[T]/(T^2-e_1T+e_2),\quad T\longmapsto x_1,
 \qquad B=A\oplus A\xi.
 \end{gathered}
\end{equation}
Here $\xi=f_1(\sigma(t)-t)$ generates the pullback of
$\mathcal I_{\Gamma_{r,l}}$.

Use the order $q_1,E',q_2,U'$ in the tensor construction, with
increasing orders within $E',U'$, and group the two factors at
$q_1,q_2$ by the graded tensor symmetry. We use the induced sign
identification throughout. Set
\[
 L_j=[Ba_j\oplus Bb_j\xrightarrow{(x_j\ \ z)}B],\qquad
 C^\bullet=(L_1\otimes_B L_2\otimes\operatorname{sgn}_2)^{S_2}.
\]
The terms of $C^\bullet$ have degrees $-2,-1,0$.
 Its degree-zero term is $A\xi$. In degree $-1$, an $A$-basis is
\[
 \begin{aligned}
 u&=a_1-a_2,& v&=x_2a_1-x_1a_2,\\
 w&=z(a_1-a_2)+x_2b_1-x_1b_2,&
 \varepsilon&=b_1-b_2.
 \end{aligned}
\]
These are the anti-invariant vectors; they form a basis because
$1,x_2$ is an $A$-basis of $B$. The differential sends
$u$ to $-\xi$ and sends $v,w,\varepsilon$ to zero.
In degree $-2$, the Koszul sign cancels the external sign, so an
$A$-basis is
\[
 \begin{aligned}
 c_1&=a_1\otimes a_2,\quad
 &c_2&=a_1\otimes b_2+b_1\otimes a_2,\\
 c_3&=x_1a_1\otimes b_2+x_2b_1\otimes a_2,\quad
 &c_4&=b_1\otimes b_2.
 \end{aligned}
\]
The tensor differential gives
\[
 dc_1=-v,\qquad dc_2=-w,\qquad
 dc_3=zv-e_1w+e_2\varepsilon,\qquad dc_4=-z\varepsilon.
\]
It follows that
\begin{equation}\label{eq:mixed-descended}
 H^0(C^\bullet)=0,\qquad
 H^{-1}(C^\bullet)=A/(z,e_2)[\varepsilon].
\end{equation}

For $j\in E'$, write the evaluation of the section as $(f_j,g_j)$ with
$f_j,g_j\in A$, and put $J=(f_j,g_j:j\in E')$.
Let $P$ be the tensor product of these evaluation complexes;
then $P\in D^{\le0}(A)$ and $H^0(P)=A/J$.
At every label in $U'$ the section is nonzero, so the
maps~\eqref{eq:prop512-kernel-map} identify the remaining complexes
with their determinant lines in degree $-1$.
Let $Q$ be the tensor product of these lines on
$\operatorname{Spec}A$.
Since $C^\bullet\in D^{\le-1}(A)$, its tensor product with
$P$ and $Q[l-1]$ has highest cohomology
$H^{-1}(C^\bullet)\otimes_A H^0(P)\otimes_A Q$ in degree $-l$.
Thus
\begin{equation}\label{eq:prop512-full-local-modules}
 \begin{aligned}
 M_{r,l}|_{\operatorname{Spec}A}
   &\simeq A/(J,z,e_2)[\varepsilon]\otimes_A Q,\\
 a^*M_{r,l}|_{\operatorname{Spec}B}
   &\simeq B/(J,z,x_1x_2)[\varepsilon]\otimes_A Q.
 \end{aligned}
\end{equation}

The selected divisor is reduced here. Order its labels with $q_1$
first, and evaluate its sign line at this ordering to choose the
basis with value $1$. Its map~\eqref{eq:selected-vandermonde-map}
then sends $\overline1$ to $[1]$ in
$H^0(L_1\otimes_A P)=B/(J,x_1,z)$.
Put $\theta_2=a_2\wedge b_2$ and use these same orders in
\eqref{eq:fixed-splitting-target}. The source is
$B/(J,x_1,z)\theta_2\otimes_A Q$.
Define
\begin{equation}\label{eq:mixed-split-map}
 \begin{aligned}
 B/(J,x_1,z)\theta_2\otimes_A Q
 &\longrightarrow B/(J,z,x_1x_2)[\varepsilon]\otimes_A Q,\\
 \overline b\theta_2\otimes q&\longmapsto-bx_2[\varepsilon]\otimes q.
 \end{aligned}
\end{equation}
The relations $J,x_1,z$ annihilate its image, so the map is well defined.
To  compare its restriction to $D(\xi)$ with \eqref{eq:kappa-fixed-splitting}, recall that
\eqref{eq:prop512-kernel-map} sends $\theta_2$ to
$\eta_2=-za_2+x_2b_2$. In $L_1\otimes_B L_2$,
$
 \eta_2+x_2\varepsilon=-d(b_1\otimes a_2)$.
Consequently, the map ~\eqref{eq:kappa-fixed-splitting}
sends $\overline1\theta_2\otimes q$ to
$[1\otimes\eta_2]\otimes q=-x_2[\varepsilon]\otimes q$.
Moving selected classes  past
$\eta_2$ introduces no Koszul sign. Thus this is exactly
\eqref{eq:mixed-split-map}.

For the duality adjoint, use the generator $\lambda$ of
$\omega_{B/A}$ from Lemma~\ref{lem:addition-duality}.
Since $x_2=(e_1+\xi)/2$, we have $\lambda(x_2)=1$.
Tensoring~\eqref{eq:mixed-split-map} over $B$ with $\omega_{B/A}$
and using coherent duality sends $(\theta_2\otimes\lambda)\otimes q$
to the functional
$
 b\longmapsto-\lambda(x_2b)[\varepsilon]\otimes q$,
where coefficients in the target are taken modulo $(J,z,e_2)$.
The counit of adjunction evaluates this functional at $1$ (\StacksTag{08YP}).
Moreover, \eqref{eq:prop512-quadratic-algebra} gives
$B/(J,x_1,z)\simeq A/(J,z,e_2)$. Hence the adjoint is an isomorphism
\[
 \begin{aligned}
 A/(J,z,e_2)(\theta_2\otimes\lambda)\otimes_A Q
 &\longrightarrow A/(J,z,e_2)[\varepsilon]\otimes_A Q,\\
 (\theta_2\otimes\lambda)\otimes q
 &\longmapsto-[\varepsilon]\otimes q.
 \end{aligned}
\]

These local maps agree on overlaps by uniqueness and descend to
$Y_U$, giving~\eqref{eq:split-comparison}. Their adjoints under
\eqref{eq:finite-duality-Hom} are therefore isomorphisms at the
mixed double points.

\smallskip
\noindent\textit{Disjoint selected and residual divisors.}
Let $y=(E,U,F,s)\in Z_U\cap W$. The section does not vanish on
$\Supp U$. The contraction maps \eqref{eq:prop512-kernel-map}
identify the residual factors with their determinant lines in degree $-1$.
Their permutation sign cancels the external residual sign.
The K\"unneth formula, together with~\eqref{eq:fixed-splitting-target},
identifies $\alpha^W_{r,l}$ with~\eqref{asfarfsfweADJFHBAJDHRFB}
tensored with the residual line. Thus it is an isomorphism near $y$.
The map $a$ is \'etale at the unique selected splitting. After an
\'etale base change it is the identity on the only component meeting
$Z_U$, and the trace identifies the duality counit with the identity.
The adjoint is consequently an isomorphism here as well.
Outside $\Sigma_r$, both sides of~\eqref{eq:Phi-U} vanish.
The two cases prove the proposition.
\end{proof}

\begin{lemma}\label{lem:S2}
Let $(R,\mathfrak m)$ be a regular local ring and let $P^\bullet$ be a finite free complex in degrees $[-n,0]$.  Put
$M_j=H^{-j}(P^\bullet)$, $c_j=2(n-j)$.

Fix $0\le q\le n$.  Assume
\begin{enumerate}
\item $\codim\Supp M_q\ge c_q$;
\item for every $j>q$, the module $M_j$ is either zero or Cohen--Macaulay of codimension $c_j$.
\end{enumerate}
Then $M_q$ satisfies Serre's condition $S_2$, and every minimal prime of $M_q$ has height exactly $c_q$.
\end{lemma}

\begin{proof}
Localize at $\mathfrak p\in\Supp M_q$ and put $S=R_{\mathfrak p}$,
$h=\dim S$. The canonical-truncation spectral sequence is
\[
 E_2^{a,j}=\Ext_S^a((M_j)_{\mathfrak p},S)
 \Longrightarrow H^{a+j}(\RHom_S(P^\bullet_{\mathfrak p},S)),
 \quad d_s:E_s^{a,j}\to E_s^{a+s,j-s+1}.
\]
The spectral sequence is bounded, and the abutment vanishes above degree~$n$.
For $j>q$, Cohen--Macaulayness gives
\begin{equation}\label{eq:CMext}
 E_2^{a,j}=0\quad\text{unless }a=c_j
\end{equation}
(\StacksTag{0AWX}, \StacksTag{0B5A}).

If $h\ge c_q+1$, the term $E_2^{h,q}$ has no outgoing differential,
since all targets have Ext-degree greater than $h$.
A nonzero incoming $d_s$ would require
$h-s=c_{q+s-1}=c_q-2s+2$, hence $s=c_q+2-h\le1$.
It therefore survives, but its total degree $h+q$ is greater than $n$.
Thus
\begin{equation}\label{eq:exttop}
 \Ext_S^h((M_q)_{\mathfrak p},S)=0\qquad(h\ge c_q+1).
\end{equation}
The same argument for $E_2^{h-1,q}$, now with $h\ge c_q+2$, gives
\begin{equation}\label{eq:extnext}
 \Ext_S^{h-1}((M_q)_{\mathfrak p},S)=0\qquad(h\ge c_q+2):
\end{equation}
its only possible incoming differential would require
$s=c_q+3-h\le1$, and again its total degree exceeds $n$.

At a minimal prime of $M_q$, the top Ext is nonzero. Thus
\eqref{eq:exttop} and the codimension bound force its height to be
$c_q$, and $\dim_S(M_q)_{\mathfrak p}=h-c_q$.
The Auslander--Buchsbaum formula and
\eqref{eq:exttop}--\eqref{eq:extnext} now give
\[
 \depth_S(M_q)_{\mathfrak p}\ge\min\{2,h-c_q\}.
\]
This is $S_2$; in particular, there are no embedded associated primes.
\end{proof}

\begin{remark}
While we only need Lemma~\ref{lem:S2} for the stated function $c_j=2(n-j)$,
the same proof works for other codimensions $c_j$, provided that
$c_j+j<c_q+q$ for every $j>q$ and $c_q+q\ge n$.
\end{remark}

\begin{proof}[Proof of Theorem~\ref{sdfgnsfghsftjt}]\label{sec:arbitrary-k}
Fix $k$ and descend through $l=k,k-1,\ldots,0$. For the values of $l$
already treated, Lemma~\ref{lem:CMcandidate} gives Cohen--Macaulay
cohomology sheaves with the asserted codimensions.
For the current value of $l$, put $r=k-l$.
Lemma~\ref{prop:support} and Lemma~\ref{lem:S2} imply that
$M_{r,l}$ is pure of codimension $2r$ and satisfies $S_2$.

Proposition~\ref{prop:split-comparison} gives~\eqref{eq:Phi-U}
on $U_{r,l}$. Both sheaves are $S_2$ and have support $\Sigma_r$,
since they agree on the dense open subset $V_{r,l}$.
Its complement in the support has codimension at least two by
\eqref{eq:bigopen}. Hence~\eqref{eq:Phi-U} extends uniquely to
\eqref{eq:corrected-cohomology} by \StacksTag{0E9I}.
The resulting isomorphism makes $M_{r,l}$ Cohen--Macaulay and
identifies its scheme-theoretic support with $\Sigma_r$, completing
the induction.

Over $V_{r,l}^{\mathrm{red}}$, the selected splitting identifies
$\beta_{r,l}$ with an isomorphism. Since $Z_{r,l}$ is smooth and
$\beta_{r,l}$ is finite birational, it is the normalization of
$\Sigma_r$. Its pushforward of a line bundle is torsion-free of
generic rank one, giving the remaining assertions.
\end{proof}

\section{Broken Loom}\label{BrokenLoomSection}

This section transforms the three-mega-block decomposition adapted to the
wall-crossing map $\psi:M\dashrightarrow\mathbb P^{3g-3}$ into a
four-mega-block decomposition adapted to the forgetful  morphism
$
\zeta:M\longrightarrow N$,
$(F,s)\longmapsto F$.
The resulting arrangement is the input for the final comparison with
$D^b(N)$.

\begin{lemma}\label{hgc,cmnvcmvc}
$D^b(M)$ has a semi-orthogonal decomposition into admissible subcategories
arranged into three mega-blocks, as follows:
$$\Bigl\langle
\bigl\langle     Z^{3+j+2k-g}{\theta^*}^{k+j+1}{\cE^{\boxtimes j}}\bigr\rangle_{ j+k\le g-2\atop j,k\ge0},\qquad\qquad\qquad\qquad\qquad$$
\begin{equation}\label{asfvRGqrWRH}
\bigl\langle Z^{2+j+2k-g}{\theta^*}^{k+j}{\cE^{\boxtimes j}}\bigr\rangle_{ j+k\le g-2\atop j,k\ge0},\ 
\end{equation}
$$\qquad\qquad\qquad\qquad\qquad\bigl\langle Z^{1+j+2k-g}{\theta^*}^{k+j-1}{\cE^{\boxtimes j}}\bigr\rangle_{ j+k\le g-1\atop j,k\ge0}
\Bigr\rangle.$$

Within each mega-block, the blocks are ordered first by decreasing $k$ and then, for fixed $k$, by decreasing $j$.
\end{lemma}

\begin{proof}
Apply the common twist $\omega_M^{3-g}$ to the decomposition of
Theorem~\ref{CrossTheorem}.  The corresponding full helix rotations differ from this twist only
by cohomological shifts, which do not change the admissible subcategories.

\begin{figure}[htbp]
\includegraphics[width=\textwidth]{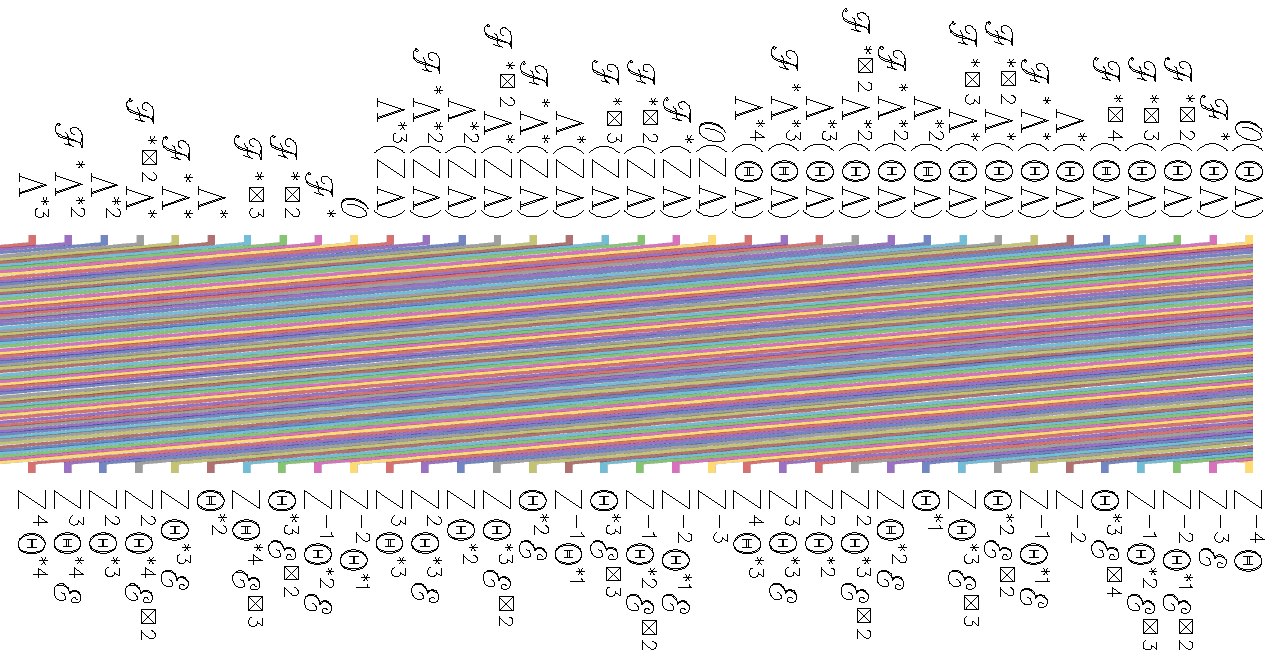}
\caption{Repeated helix mutation in genus $5$}\label{spiralbraid}
\end{figure}

The additional common twist $\theta^{5-2g}$ is pulled back from $N$.
Since $\omega_M=Z\theta^{-2}$, their product is
$
 \omega_M^{3-g}\theta^{5-2g}=Z^{3-g}\theta^{-1}$.
For a block of the first mega-block, use
$\mathcal F=\mathcal E(-Z)$ and $\boldsymbol\Lambda=\theta Z^{-2}$:
\[
 (\boldsymbol\Lambda^*)^k\mathcal F^{*\boxtimes j}
 \simeq \mathcal E^{\boxtimes j}\otimes
      \bigl(\Lambda^{*\boxtimes j}\boxtimes
                  Z^{-j}\boldsymbol\Lambda^{-k-j}\bigr)
 \simeq \mathcal E^{\boxtimes j}\otimes
      \bigl(\Lambda^{*\boxtimes j}\boxtimes
                  Z^{2k+j}\theta^{-k-j}\bigr).
\]
Tensoring the source by $\Lambda^{*\boxtimes j}$ is an
autoequivalence, so this source factor does not affect the essential
image.  Multiplication by the common twist gives
$Z^{3+j+2k-g}\theta^{-k-j-1}\mathcal E^{\boxtimes j}$.
The second and third original mega-blocks have the additional
factors $Z\boldsymbol\Lambda=\theta Z^{-1}$ and
$\theta\boldsymbol\Lambda=\theta^2 Z^{-2}$, respectively.
They therefore give the other two expressions
in~\eqref{asfvRGqrWRH}.  The common twists do not change the order
of blocks inside a mega-block.
\end{proof}

\begin{remark}
The twist by the line bundle $\omega_M^{3-g}$ is required to create blocks of the semi-orthogonal decomposition 
compatible with the contraction $\zeta$.
In~contrast, the twist by the line bundle ${\theta^*}^{2g-5}$, which is pulled back from~$N$,
is unnecessary; we use it only to simplify the formulas.
In the illustrations of the braid, we ignore this twist even though we still use it to label the blocks, for consistency.
\end{remark}

\begin{figure}[htbp]
\includegraphics[width=\textwidth]{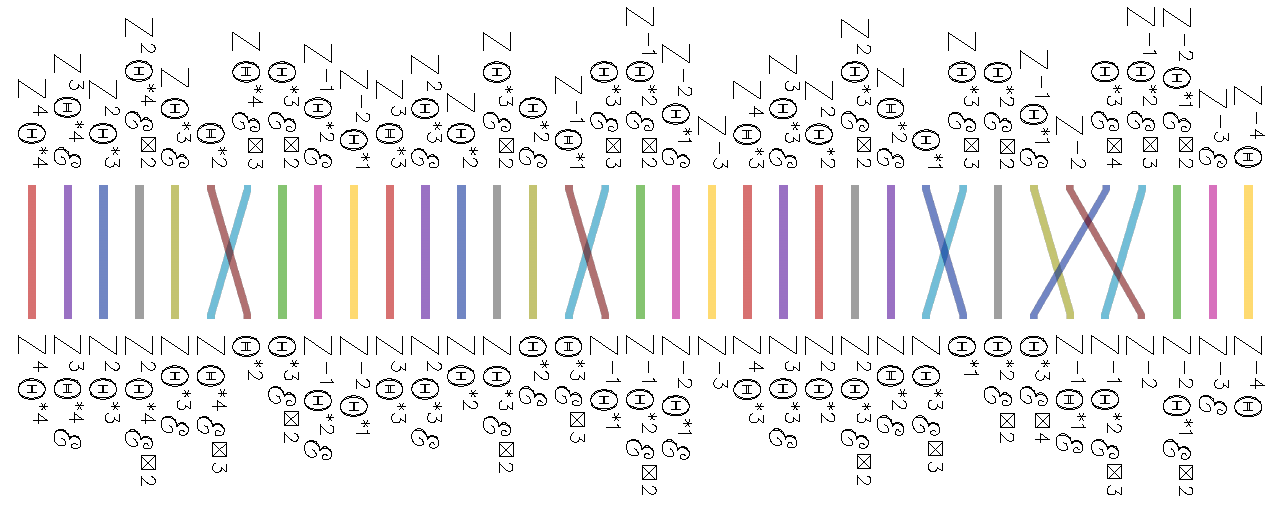}
\caption{Reordering of the blocks in genus $5$}\label{brokdfbqerhenloom}
\end{figure}

The next step is crucial. As it turns out, the blocks within each mega-block of  
\eqref{asfvRGqrWRH} can be reordered differently. This is illustrated in Figure~\ref{brokdfbqerhenloom}.

\begin{theorem}\label{wrgasrharharh}
$D^b(M)$ has a semi-orthogonal decomposition into admissible subcategories
arranged into three mega-blocks, as follows:
$$\Bigl\langle
\bigl\langle     Z^{3+\lambda-g}{\theta^*}^{\lambda-k+1}{\cE^{\boxtimes \lambda-2k}}\bigr\rangle_{\lambda-k\le g-2\atop \lambda-2k,k\ge0},\qquad\qquad\qquad\qquad\qquad$$
\begin{equation}\label{asfvRGdfvagqerqrWRH}
\bigl\langle Z^{2+\lambda-g}{\theta^*}^{\lambda-k}{\cE^{\boxtimes \lambda-2k}}\bigr\rangle_{\lambda-k\le g-2\atop \lambda-2k,k\ge0},\ 
\end{equation}
$$\qquad\qquad\qquad\qquad\qquad\bigl\langle Z^{1+\lambda-g}{\theta^*}^{\lambda-k-1}{\cE^{\boxtimes \lambda-2k}}\bigr\rangle_{ \lambda-k\le g-1\atop \lambda-2k,k\ge0}
\Bigr\rangle.$$
Within each of the three mega-blocks, the blocks are arranged 
in  decreasing order of $\lambda$. Blocks with the same $\lambda$
are arranged in  decreasing order of $k$.
\end{theorem}

We start with a lemma.

\begin{lemma}\label{SRBADRHADHATD}
Suppose $2<d\leq 2g+1$ and $1\leq j\leq \lfloor {d-1\over 2}\rfloor$. 
Let $D\in\Sym^a C$ and $D'\in\Sym^bC$ with $a,\ b\leq \min(d+g-2j-1, j)$,
and let $t$ be an integer satisfying
$$
    a-j-1 < t < d+g-2j-1-b.
$$
Suppose, further, that 
$2t<a-b$. 
Then 
\begin{equation}\label{afgawrhr}
R\Gamma_{M_j(d)}({\cF^{\boxtimes a}_D}^*\otimes\cF^{\boxtimes b}_{D'}\otimes\bLambda^t)=0
\end{equation}
\end{lemma}

\begin{proof}
Write $D=\sum_i\alpha_i x_i$ and $D'=\sum_\mu\beta_\mu y_\mu$,
with distinct $x_i$ and distinct $y_\mu$, and put
$H=d+g-2j-1$. By \cite[Corollary~2.9]{TT} and semicontinuity,
it suffices to prove the vanishing for the point tensors
\begin{equation}\label{afasrgarggawrhr}
 \bigotimes_i(\cF_{x_i}^*)^{\otimes\alpha_i}
 \otimes\bigotimes_\mu\cF_{y_\mu}^{\otimes\beta_\mu}
 \otimes\bLambda^t.
\end{equation}
Both tensor bundles in the statement occur as the isomorphic
nonzero fibers of their respective deformations, so vanishing at the
special fiber suffices.

Rank-two Schur--Weyl duality decomposes these tensors into symmetric
powers, with determinant factors. A summand has parameters
\[
 a'=a-2u,\qquad b'=b-2v,\qquad t'=t-u+v,
 \qquad u,v\ge0.
\]
These remain in the required range: $a',b'\le\min(H,j)$, and
\[
\begin{aligned}
 t'-(a'-j-1)&=t-(a-j-1)+u+v>0,\\
 H-b'-t'&=H-b-t+u+v>0,\\
 2t'-(a'-b')&=2t-(a-b)<0.
\end{aligned}
\]
Thus it suffices to prove, throughout this range,
\begin{equation}\label{afasrghr}
 R\Gamma_{M_j(d)}\!\left(
 \bigotimes_i\Sym^{\alpha_i}\cF_{x_i}^*
 \otimes\bigotimes_\mu\Sym^{\beta_\mu}\cF_{y_\mu}
 \otimes\bLambda^t\right)=0.
\end{equation}
We suppress the constant one-dimensional determinant fibers on $C$,
which do not affect vanishing. Rank-two duality interchanges the two
families and replaces $(a,b,t)$ by $(b,a,t-a+b)$, preserving the same
three inequalities. We may therefore assume $b\ge a$, so $t<0$.
We prove~\eqref{afasrghr} by induction on $a+b$.
The case $a=b=0$ follows from Theorem~\ref{hardvanishing}.
For the inductive step, use the degree filtration of the ordered fiber
of the symmetric quotient, as in \cite[Lemma~7.7]{TT}.
More explicitly, put
\[
 A_\beta=\mathbb C[t_1,\ldots,t_\beta]/(\sigma_1,\ldots,\sigma_\beta),
 \quad A=\bigotimes_\mu A_{\beta_\mu},
 \quad G=\prod_\mu S_{\beta_\mu},
\]
where $\sigma_i$ are the elementary symmetric polynomials. This describes
the ordered fiber over $D'$, after choosing a parameter at each $y_\mu$.
Let $\mathcal T$ be the tensor product of the pullbacks of $\cF$ to
$\Spec A\times M_j(d)$, and let $\tau$ be the projection to $M_j(d)$.
Then
$ \bar\cF_{D'}^{\boxtimes b}
 =\bigl(\tau_*\mathcal T\otimes\operatorname{sgn}_G\bigr)^G$.
The degree filtration $A_{\ge m}$ induces a filtration by subbundles
of $\tau_*\mathcal T$. Its quotients are
\begin{equation}\label{eq:claim56-ver3f-associated-graded}
 \operatorname{gr}_m(\tau_*\mathcal T)
 \simeq A_m\otimes\bigotimes_\mu\cF_{y_\mu}^{\otimes\beta_\mu}.
\end{equation}
Indeed, modulo $A_{\ge m+1}$ the positive-degree terms of the transition
matrices act trivially, leaving their values at $y_\mu$.

The top degree $N=\sum_\mu\binom{\beta_\mu}{2}$ of $A$ is the
one-dimensional sign representation, generated by the product of the
Vandermonde determinants. Applying the exact sign projector to the
inclusion of this top subbundle gives
\begin{equation}\label{adrgarhaerhaetha}
 0\longrightarrow\bigotimes_\mu\Sym^{\beta_\mu}\cF_{y_\mu}
 \longrightarrow\bar\cF_{D'}^{\boxtimes b}
 \longrightarrow Q\longrightarrow0.
\end{equation}
The identification of the first term is up to a constant nonzero scalar;
the subbundle and its quotient are independent of that choice.

Schur--Weyl duality applied to~\eqref{eq:claim56-ver3f-associated-graded}
shows that $Q$ has a finite filtration whose quotients are sums of
\begin{equation}\label{eq:claim56-ver3f-graded-Q}
 \bLambda^v\otimes\bigotimes_\mu
       \Sym^{\beta_\mu-2v_\mu}\cF_{y_\mu},
 \qquad v=\sum_\mu v_\mu>0.
\end{equation}
To justify the strict inequality, the summand with all $v_\mu=0$
has trivial $G$-representation on its tensor-power factor. It can
survive the sign projector only from a sign summand of $A$.
Every alternating polynomial is divisible by the corresponding
Vandermonde, so $A$ contains the sign representation only in its top
degree. That degree has been removed in forming $Q$.

For a quotient in~\eqref{eq:claim56-ver3f-graded-Q}, the new parameters
are $(a,b-2v,t+v)$. They satisfy
\[
 \begin{gathered}
 a-j-1<t+v,\qquad 2(t+v)<a-(b-2v),\\
 H-(b-2v)-(t+v)=H-b-t+v>0.
 \end{gathered}
\]
Their total degree is smaller. The induction hypothesis therefore
annihilates the cohomology of
$\bigotimes_i\Sym^{\alpha_i}\cF_{x_i}^*\otimes Q\otimes\bLambda^t$.
The same tensor product with the middle term
$\bar\cF_{D'}^{\boxtimes b}$ of~\eqref{adrgarhaerhaetha} has zero
cohomology by~\eqref{afgaasfargarhwrhr}, because $t<0$ and the dual
symmetric-power factors are direct summands of the corresponding
point tensors. The long exact sequence of~\eqref{adrgarhaerhaetha}
now proves~\eqref{afasrghr}, completing the induction.
\end{proof}

\begin{remark}\label{srGASRHARHARH}
The proof of Lemma~\ref{SRBADRHADHATD} also shows that \eqref{afgawrhr}  holds with 
$\cF^{\boxtimes a}_D$ replaced with
$\bar\cF^{\boxtimes a}_D$
and with 
$\cF^{\boxtimes b}_{D'}$
replaced with 
$\bar\cF^{\boxtimes b}_{D'}$.
\end{remark}

\begin{proof}[Proof of Theorem~\ref{wrgasrharharh}]
Substitute $\lambda=2k+j$ in Lemma~\ref{hgc,cmnvcmvc}.  It is
enough to check the new order in the third mega-block, since the
other two differ by common line-bundle twists and have smaller
index ranges.  We must prove the vanishing from the block indexed
by $(\lambda_1,k_1)$ to the block indexed by $(\lambda_2,k_2)$ when
$\lambda_1<\lambda_2$, or when $\lambda_1=\lambda_2$ and $k_1<k_2$.
Put
\[
 j_i=\lambda_i-2k_i,\qquad
 t=(k_1+j_1)-(k_2+j_2)
      =(\lambda_1-k_1)-(\lambda_2-k_2).
\]
Using $\mathcal E=Z\mathcal F$ and
$\boldsymbol\Lambda=\theta Z^{-2}$, the required group for
$D\in\operatorname{Sym}^{j_1}C$ and
$D'\in\operatorname{Sym}^{j_2}C$ becomes
$
 R\Gamma_M\bigl(
 \boldsymbol\Lambda^t\otimes
 (\mathcal F_D^{\boxtimes j_1})^*
       \otimes\mathcal F_{D'}^{\boxtimes j_2}\bigr)
$.
Suppose first that $\lambda_1<\lambda_2$.  Apply
Lemma~\ref{SRBADRHADHATD} with degree $2g-1$, wall index $g-1$,
$a=j_1$, and $b=j_2$.  The interval condition is
\begin{equation}\label{asfasgasrgh}
 j_1-g<t<g-j_2.
\end{equation}
Indeed,
\[
 \begin{aligned}
 t-(j_1-g)&=g-(\lambda_2-k_2)+k_1>0,\\
 (g-j_2)-t&=g-(\lambda_1-k_1)+k_2>0,
 \end{aligned}
\]
using $\lambda_i-k_i\le g-1$ and $k_i\ge0$.
The separate strict inequality is
\begin{equation}\label{asfgaaasvasbsrsrgh}
 2t<j_1-j_2.
\end{equation}
It follows not from these interval bounds but directly from
\[
 2t-(j_1-j_2)
     =(j_1+2k_1)-(j_2+2k_2)=\lambda_1-\lambda_2<0.
\]
This proves the required vanishing when $\lambda_1<\lambda_2$.

Now suppose $\lambda_1=\lambda_2$ and $k_1<k_2$.  Put
$\delta=k_2-k_1>0$.  Then
$t=\delta$,
$j_1-j_2=2\delta$,
$q:=t+j_2-j_1=-\delta<0$.
By Lemma~\ref{amazinggrace}, up to one-dimensional line factors coming from
the source divisors,
\[
\bLambda^t(\cF_D^{\boxtimes j_1})^*
 \otimes\cF_{D'}^{\boxtimes j_2}
\simeq
\bLambda^q\bar\cF_D^{\boxtimes j_1}
 \otimes\cF_{D'}^{*\boxtimes j_2}.
\]
Apply Theorem~\ref{hardvanishing} with wall index $g-1$, degree $2g-1$,
$a=j_2$, $b=j_1$, and twist $q$.  We have $j_1,j_2\le g-1$,
$q\notin[0,j_2]$, and
$
j_2-g<q<g-j_1$.
Indeed, the right inequality follows from $j_1\le g-1$, while the left is
equivalent to $j_2+\delta<g$, which follows from
$\lambda-k_2\le g-1$ and $k_1<k_2$.  By 
\cite[Corollary~2.9]{TT}, we can pass the vanishing from point tensors to
$\cF_{D'}^{*\boxtimes j_2}$.
\end{proof}

Finally, we prove the following theorem:

\begin{theorem}\label{LoomTheorem}
$D^b(M)$ has a semiorthogonal decomposition whose blocks are arranged into
the following four mega-blocks:
$$
\bigl\langle 
Z^{\lambda+2-g}{\theta^*}^{\lambda-k+1}\cE^{\boxtimes \lambda-2k}
\bigr\rangle_{0\le \lambda\le g-2 \atop 0\le k\le \lfloor{\lambda\over2}\rfloor},\bigl\langle 
Z^{\lambda+3-g}{\theta^*}^{\lambda-k+1}\cE^{\boxtimes \lambda-2k}
\bigr\rangle_{0\le \lambda\le 2(g-2) \atop 0\le k\le  \lfloor{\lambda\over2}\rfloor, \lambda-k\le g-2},$$
$$\bigl\langle 
Z^{\lambda+2-g}{\theta^*}^{\lambda-k}\cE^{\boxtimes \lambda-2k}
\bigr\rangle\!\!\!{}_{0\le \lambda\le 2(g-2) \atop 0\le k\le \lfloor{\lambda\over2}\rfloor, \lambda-k\le g-2}\!\!\!\!,
\bigl\langle 
Z^{\lambda+1-g}{\theta^*}^{\lambda-k-1}\cE^{\boxtimes \lambda-2k}
\bigr\rangle\!\!\!\!{}_{g-1\le \lambda\le 2(g-1) \atop 0\le k\le  \lfloor{\lambda\over2}\rfloor, \lambda-k\le g-1}\!\!\!\!.
$$ 
Within each mega-block, the blocks are arranged in decreasing order of $\lambda$
and those with identical $\lambda$ are further arranged by decreasing $k$.
\end{theorem}

\begin{figure}[htbp]
\includegraphics[width=\textwidth]{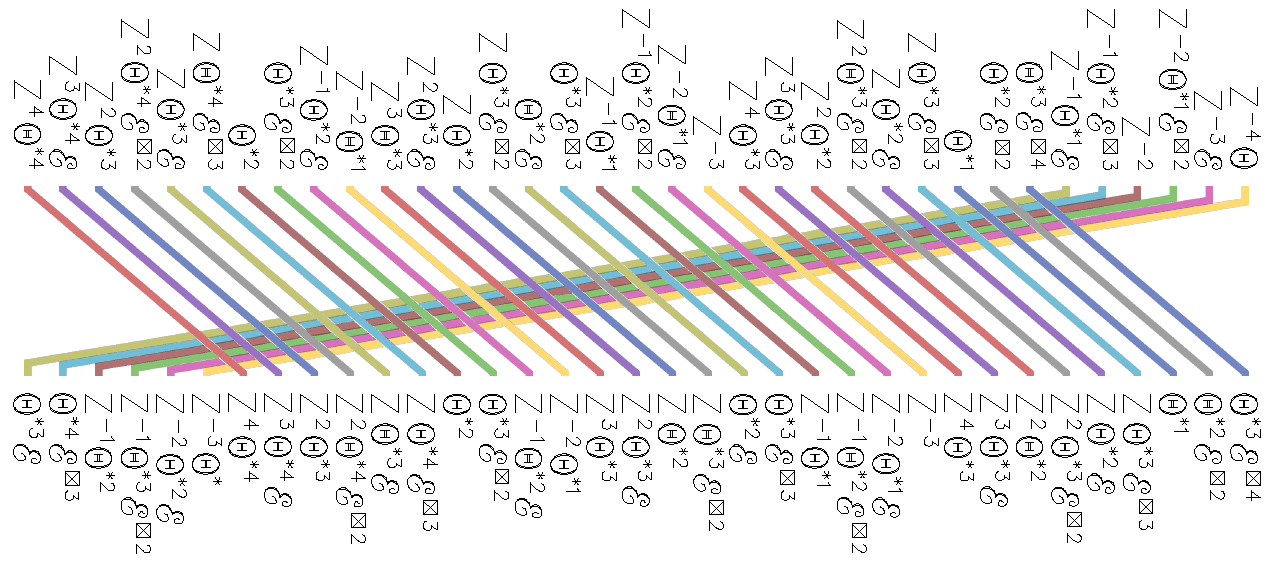}
\caption{Creation of the fourth mega-block in genus $5$}\label{lasttwist}
\end{figure}

\begin{proof}
Split the third mega-block of Theorem~\ref{wrgasrharharh} into the
part with $g-1\le\lambda\le2(g-1)$ and the part with
$0\le\lambda\le g-2$.  The latter is at the end because the
blocks are ordered by decreasing $\lambda$.  Move this terminal
part to the front by a helix mutation.  Up to the cohomological
shift of the Serre functor, which does not change an essential
image, this tensors its kernels by $\omega_M=Z\theta^{-2}$.
For each moved block the calculation is
\[
 (Z\theta^{-2})\,
 Z^{1+\lambda-g}\theta^{-(\lambda-k-1)}
             \mathcal E^{\boxtimes\lambda-2k}
 =Z^{2+\lambda-g}\theta^{-(\lambda-k+1)}
             \mathcal E^{\boxtimes\lambda-2k}.
\]
These are the kernels of the first mega-block in the statement.
On $0\le\lambda\le g-2$, the old inequality
$\lambda-k\le g-1$ is automatic.  The other three index sets and
the orders within them are unchanged.  Figure~\ref{lasttwist}
illustrates this mutation.
\end{proof}

\section{Plain Weave}\label{PlainWeaveSection}

\begin{figure}[htbp]
\includegraphics[width=0.85\textwidth]{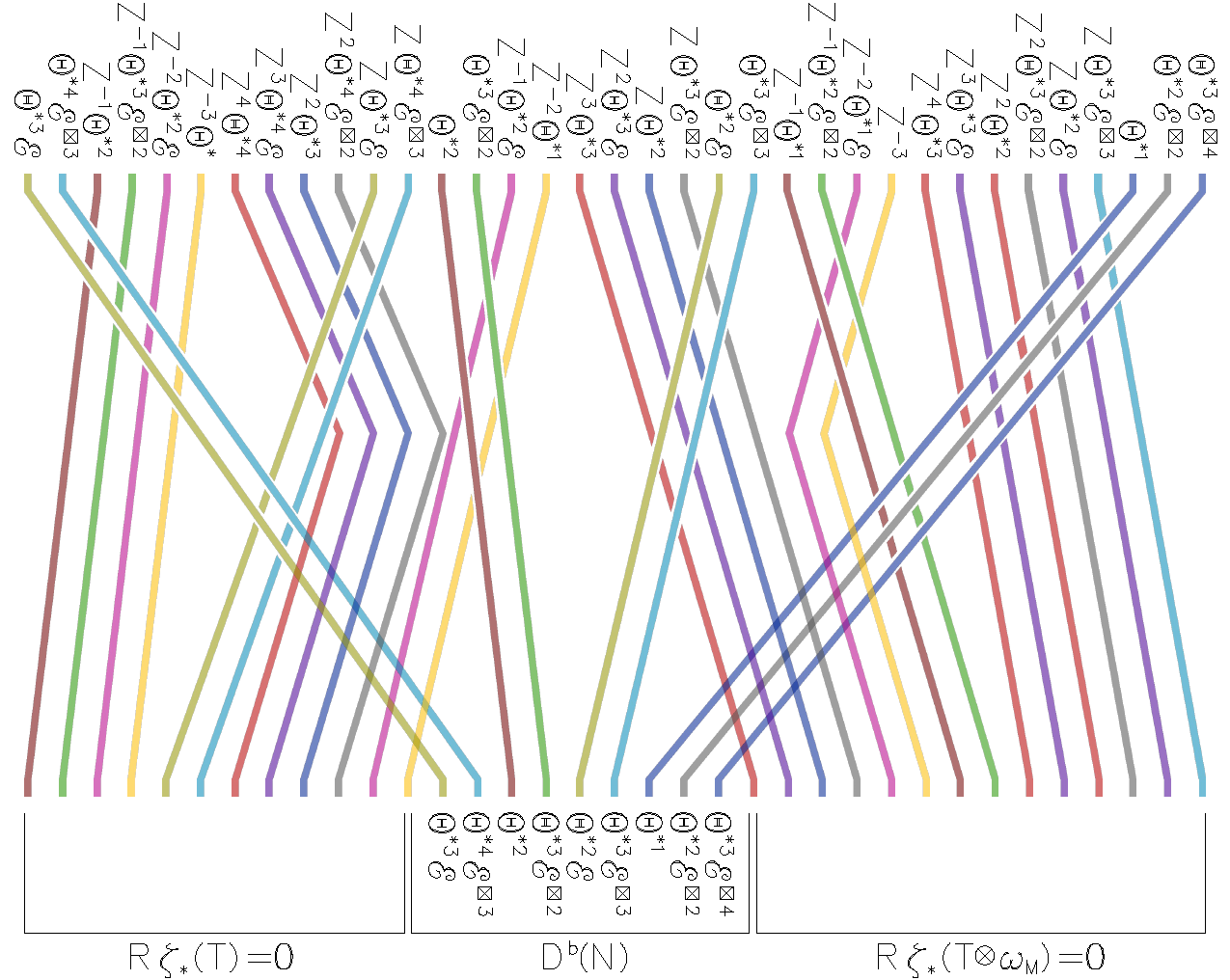}
\caption{Plain Weave in genus $5$. 
}\label{plainweave}
\end{figure}

In this final section we prove Theorem~\ref{MainTheorem} by constructing the Plain Weave mutation starting from
the decomposition of Theorem~\ref{LoomTheorem}.  
The blocks with trivial
power of $Z$ occur at $\lambda=g-2,g-3,g-2,g-1$ in mega-blocks
I, II, III, IV, respectively.  They are the pullbacks of the blocks
in Theorem~\ref{MainTheorem}.  We retain these four families, with their
orders unchanged, and mutate all other blocks to the two sides as in
Figure~\ref{plainweave}.  

\begin{notation}
Let $\mathcal A=\ker R\zeta_*$,
$\mathcal A'=\ker R\zeta_*(-\otimes\omega_M)
             =\ker R\zeta_*(-\otimes Z)$.
Since $R\zeta_*\cO_M=\cO_N$ (see Lemma~\ref{pw:direct-images}),
\[
 D^b(M)=\langle\mathcal A,L\zeta^*D^b(N)\rangle
       =\langle L\zeta^*D^b(N),\mathcal A'\rangle.
\]
We will put the nonretained blocks from the mega-blocks I and II on the left in
$\mathcal A$, and those from the mega-blocks  III and IV on the right in $\mathcal A'$.
A~final helix mutation tensors the latter by $\omega_M$, up to a
cohomological shift, and moves them to the left into $\mathcal A$.
\end{notation}

We first prove several statements needed for
Plain Weave.  

\begin{lemma}\label{pw:direct-images}
Let $p_q:\Sym^qC\times N\to N$ be the projection.  For every $q\geq0$,
$$ R\zeta_*Z^{-q}\simeq
       \bigl(Rp_{q*}\cE^{\boxtimes q}\bigr)^\vee,
 \quad
 R\zeta_*Z^{q+1}\simeq Rp_{q*}\cE^{\boxtimes q}.
$$In particular, $R\zeta_*\cO_M\simeq R\zeta_*Z\simeq\cO_N$.
\end{lemma}

\begin{proof}
Choose an effective divisor $D_0$ with $R^1p_{1*}\cE(D_0)=0$,
and represent
\[
 Rp_{1*}\cE\simeq[\mathcal U\xrightarrow{u}\mathcal V],\qquad
 \mathcal U=p_{1*}\cE(D_0),\quad
 \mathcal V=p_{1*}(\cE(D_0)|_{D_0\times N}),
\]
in degrees $0,1$. Write $\rk\mathcal U=a$; then
$\rk\mathcal V=a-1$ because $\chi(C,E)=1$.
In the projective bundle of lines $\pi:\mathbb P(\mathcal U)\to N$,
$M$ is the zero scheme of the induced section of
$\pi^*\mathcal V(1)$, and $\mathcal O(-1)|_M=Z$.
The section is regular because $M$ has the expected codimension.

Tensor its Koszul resolution by $\mathcal O(q)$ and apply
$R\pi_*$. For $q\ge0$, the projective bundle formula identifies the
result with
$
 R\zeta_*Z^{-q}\simeq
       \bigl(\operatorname{Sym}^q[\mathcal U\to\mathcal V]\bigr)^\vee
$,
where the symmetric power uses the graded tensor symmetry.
The relative K\"unneth formula gives
\[
 Rp_{q*}\cE^{\boxtimes q}
 \simeq\bigl((Rp_{1*}\cE)^{\otimes^{\mathbf L}q}\bigr)^{S_q}
 \simeq\operatorname{Sym}^q[\mathcal U\to\mathcal V].
\]
These identifications are induced by tautological evaluation and
the Koszul augmentation. The second formula follows by relative
duality, since $\zeta^!\mathcal O_N=Z$. Taking $q=0$ proves the
last assertion.
\end{proof}

\begin{notation}\label{pw:block-families}
For $0\leq h\leq g-1$ and $c\in\mathbb Z$, write
\begin{equation}\label{pw:family-definition}
 \mathcal B_h(c)=
 \left\langle\theta^{c+j}\cE^{\boxtimes(h-2j)}
 \right\rangle_{j=\lfloor h/2\rfloor,\ldots,0}
 \subset D^b(N).
\end{equation}
These are admissible subcategories.  Indeed, choose
$H\in\{g-2,g-1\}$ with $h=H-2m$ for some $m\geq0$.
The terms with $j\geq m$ in the $Z^0$ part of mega-block III or IV of
Theorem~\ref{LoomTheorem} form a semiorthogonal sequence on $M$.
Since $L\zeta^*$ is fully faithful by Lemma~\ref{pw:direct-images},
these terms also form a semiorthogonal sequence on $N$.
After a common twist by a power of $\theta$, their generated category
is \eqref{pw:family-definition}.  In particular,
\begin{equation}\label{pw:prefix}
 \mathcal B_{H-2m}(c+m)=
 \left\langle\theta^{c+j}\cE^{\boxtimes(H-2j)}
 \right\rangle_{j=\lfloor H/2\rfloor,\ldots,m}.
\end{equation}
\end{notation}

\begin{lemma}\label{pw:pieri}
Fix $x\in C$.  Consider the maps given by addition of $x$:
\[
 a_x:\Sym^dC\longrightarrow\Sym^{d+1}C,\qquad
 i_x:\Sym^{d-1}C\longrightarrow\Sym^dC.
\]
For $d\geq1$ there is a natural exact sequence on $\Sym^dC\times N$:
\begin{equation}\label{pw:pieri-sequence}
 0\rightarrow (a_x\times1)^*\cE^{\boxtimes(d+1)}
 \rightarrow \cE_x\otimes\cE^{\boxtimes d}
 \rightarrow (i_x\times1)_*
       \bigl(\det\cE_x\otimes\cE^{\boxtimes(d-1)}\bigr)
 \rightarrow0.
\end{equation}
Consequently, for $1\leq d\leq g-2$,
\begin{equation}\label{pw:pieri-images}
 \cE_x\otimes\langle\cE^{\boxtimes d}\rangle
 \subset
 \langle\cE^{\boxtimes(d+1)},\theta\cE^{\boxtimes(d-1)}\rangle.
\end{equation}
For $d=0$, the corresponding assertion is
$\cE_x\in\langle\cE\rangle$.
\end{lemma}

\begin{proof}
The inverse image of $a_x(\Sym^dC)$ in $C^{d+1}$ is the reduced
simple-normal-crossing divisor $Y=\bigcup_jY_j$, where
$Y_j=\{p_j=x\}$. Tensor its augmented \v Cech resolution by
$\bigotimes_jp_j^*\cE$, take finite direct image to $\Sym^dC\times N$,
and then take $S_{d+1}$-invariants.

For an intersection of $b$ of the $Y_j$, the stabilizer is
$S_b\times S_{d+1-b}$. The orientation sign in the \v Cech resolution
gives $\bigwedge^b\cE_x$ for the first factor, so the corresponding
term is the pushforward of
$
 \bigwedge^b\cE_x\otimes\cE^{\boxtimes(d+1-b)}
$
under addition of $(b-1)x$. The terms with $b\ge3$ vanish because
$\rk\cE_x=2$. The augmentation term is
$(a_x\times1)^*\cE^{\boxtimes(d+1)}$ by finite base change.
The surviving terms $b=1,2$ therefore give
\eqref{pw:pieri-sequence}.

The projection formula gives~\eqref{pw:pieri-images}, since
$\det\cE_x\simeq\theta$.
\end{proof}

\begin{lemma}\label{pw:admissible-deformation}
Let $N$ be a smooth projective variety over $\mathbb C$, and let
$\mathcal S\subset D^b(N)$ be an admissible subcategory.
For $t\in\mathbb A^1(\mathbb C)$, let
$i_t:N\hookrightarrow N\times\mathbb A^1$ be the fiber inclusion.
Suppose that $\mathcal K$ is a perfect complex on
$N\times\mathbb A^1$ and that its derived fibers
$
 \mathcal K_t:=Li_t^*\mathcal K
$
are pairwise isomorphic for $t\ne0$. 

If $\mathcal K_0\in\mathcal S$,
then $\mathcal K_t\in\mathcal S$ for every $t$.
\end{lemma}

\begin{proof}
Let $\Pi$ be the projection to $\mathcal S^\perp$, with its
Fourier--Mukai kernel
\cite[Lemma~2.9 and Theorem~7.1]{kluznet}.
Applying this kernel to the family gives a perfect complex
$\mathcal Q$ on $N\times\mathbb A^1$ with fibers
$\mathcal Q_t\simeq\Pi(\mathcal K_t)$.
Since $\mathcal Q_0=0$, the derived Nakayama lemma implies that
$\Supp\mathcal Q$ misses $N\times\{0\}$. Its image in
$\mathbb A^1$ is closed by properness of $N$.
The fibers over $t\ne0$ are isomorphic. If one were nonzero, this
closed image would contain $\mathbb A^1\setminus\{0\}$ and hence
$0$, a contradiction. Thus all $\mathcal Q_t$ vanish.
\end{proof}

\begin{proposition}\label{pw:tensor-products}
If $a,b\geq0$ and $a+b\leq g-1$, then
\begin{equation}\label{pw:tensor-inclusions}
 \langle\cE^{\boxtimes a}\rangle\otimes
       \langle\cE^{\boxtimes b}\rangle \subset\mathcal B_{a+b}(0),\quad 
 \langle\cE^{\boxtimes a}\rangle\otimes
       \langle\bar\cE^{\boxtimes b}\rangle
       \subset\mathcal B_{a+b}(0).
\end{equation}
The left sides mean the categories generated by tensor products of
objects of the indicated essential images.
\end{proposition}

\begin{proof}
Iterating \eqref{pw:pieri-images} first shows that, for arbitrary
points $x_1,\ldots,x_h$ with repetitions allowed,
\begin{equation}\label{pw:point-products}
 \cE_{x_1}\otimes\cdots\otimes\cE_{x_h}\in\mathcal B_h(0)
       \qquad(0\leq h\leq g-1).
\end{equation}
Indeed, tensoring $\theta^j\langle\cE^{\boxtimes(h-2j)}\rangle$
by another $\cE_x$ gives terms with indices $j$ and $j+1$ in
$\mathcal B_{h+1}(0)$.  

For every divisor $D=\sum x_j$, both tensor bundles
$\cE_D^{\boxtimes h}$ and $\bar\cE_D^{\boxtimes h}$ admit the
isotrivial deformation to $\bigotimes_j\cE_{x_j}$ of
\cite[Corollary~2.9]{TT}.  Tensor two such families.
Their special fiber satisfies \eqref{pw:point-products},
so their other fibers belong to $\mathcal B_{a+b}(0)$ by Lemma~\ref{pw:admissible-deformation}. 
 This proves the desired assertion for every product of
evaluation bundles, including nonreduced divisors.

Finally, compose the Fourier--Mukai functor on
$\Sym^aC\times\Sym^bC$ given by the tensor product of the two
kernels with the projection to $\mathcal B_{a+b}(0)^\perp$.
Its kernel has zero derived fiber at every $(D,D')$, so vanishes
by the derived Nakayama lemma. This proves
\eqref{pw:tensor-inclusions} for the entire essential images.
\end{proof}

\begin{lemma}\label{pw:pushforward-bounds}
Let $0\le n\le g-1$ be an integer and let $c,s\in\mathbb Z$. Then
\[
 R\zeta_*\langle Z^s\theta^c\cE^{\boxtimes n}\rangle
 \subset\mathcal B_{n-s}(c+s)
 \qquad\text{if }n-g+1\le s\le0,
\]
and
\begin{equation}\label{pw:right-bound}
 R\zeta_*\langle Z^{s+1}\theta^c\cE^{\boxtimes n}\rangle
 \subset\mathcal B_{n+s}(c)
 \qquad\text{if }0\le s\le g-1-n.
\end{equation}
\end{lemma}

\begin{proof}
Put $V_q=Rp_{q*}\cE^{\boxtimes q}$. For $0\le q\le g-1$,
rank-two duality, Lemma~\ref{amazinggrace}, and relative duality give
\[
 V_q\in\langle\cE^{\boxtimes q}\rangle,\quad
 V_q^\vee\simeq\theta^{-q}\otimes
 \cP_{\bar\cE^{\boxtimes q}}\bigl(
 \Lambda^{*\boxtimes q}\otimes\mathscr L_q^*
 \otimes\omega_{\Sym^qC}[q]\bigr)
 \in\theta^{-q}\langle\bar\cE^{\boxtimes q}\rangle.
\]
For $F\in\langle\cE^{\boxtimes n}\rangle$, the projection formula,
Lemma~\ref{pw:direct-images}, and Proposition~\ref{pw:tensor-products}
give
\[
\begin{aligned}
 R\zeta_*(Z^{-q}\otimes L\zeta^*(\theta^cF))
  &\simeq\theta^cF\otimes V_q^\vee
    \in\mathcal B_{n+q}(c-q),\\
 R\zeta_*(Z^{q+1}\otimes L\zeta^*(\theta^cF))
  &\simeq\theta^cF\otimes V_q
    \in\mathcal B_{n+q}(c),
\end{aligned}
\]
provided $n+q\le g-1$. Substitute $q=-s$ and $q=s$, respectively.
\end{proof}

\begin{figure}[htbp]
\includegraphics[width=\textwidth]{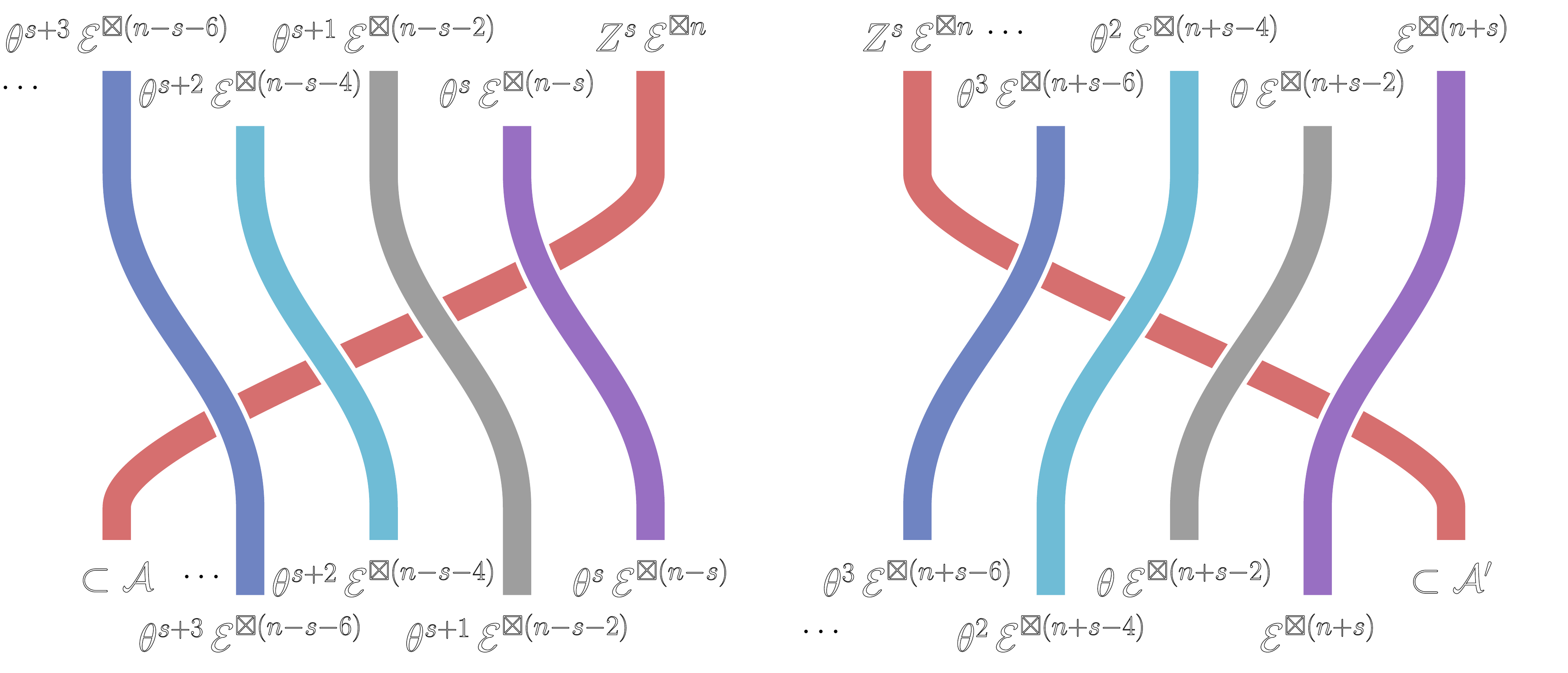}
\caption{Plain Weave mutations. Only strands with nonnegative tensor exponent
are included.}
\label{wdfkljgnaldkrjgbalerkjbge}
\end{figure}

\begin{proposition}[Plain Weave mutations]\label{pw:enlarged-mutations}
Let $\mathcal B\subset D^b(N)$ be admissible.  For $T\in D^b(M)$,
\begin{align}
 R\zeta_*T\in\mathcal B
 &\ \Longrightarrow\ 
       \mathbf L_{L\zeta^*\mathcal B}(T)\in\mathcal A,
       \label{pw:left-criterion}\\
 R\zeta_*(T\otimes Z)\in\mathcal B
 &\ \Longrightarrow\ 
       \mathbf R_{L\zeta^*\mathcal B}(T)\in\mathcal A'.
       \label{pw:right-criterion}
\end{align}
Consequently, for $n\in\mathbb Z_{\ge0}$ and $c,s\in\mathbb Z$, one has
\begin{align}
 \mathbf L_{L\zeta^*\mathcal B_{n-s}(c+s)}
       \langle Z^s\theta^c\cE^{\boxtimes n}\rangle
       &\subset\mathcal A
       &&(n-g+1\le s\le0),\label{pw:left-mutation}\\
 \mathbf R_{L\zeta^*\mathcal B_{n+s}(c)}
       \langle Z^s\theta^c\cE^{\boxtimes n}\rangle
       &\subset\mathcal A'
       &&(0\le s\le g-1-n).\label{pw:right-mutation}
\end{align}
These mutations are illustrated in Figure~\ref{wdfkljgnaldkrjgbalerkjbge} (where $c=0$).
\end{proposition}

\begin{proof}
The right and left adjoints of $L\zeta^*$ are $R\zeta_*$ and
$R\zeta_*(-\otimes Z)$, respectively. Under the two hypotheses,
the defining mutation triangles are
\[
\begin{gathered}
 L\zeta^*R\zeta_*T\longrightarrow T
       \longrightarrow\mathbf L_{L\zeta^*\mathcal B}T,\\
 \mathbf R_{L\zeta^*\mathcal B}T\longrightarrow T
       \longrightarrow L\zeta^*R\zeta_*(T\otimes Z).
\end{gathered}
\]
Applying the respective adjoints proves both criteria.
The remaining assertions follow from Lemma~\ref{pw:pushforward-bounds}.
\end{proof}

We will also use the common twists of these statements.  Conjugating by
$(-)\otimes Z^t$ replaces $L\zeta^*\mathcal B$ by
$Z^tL\zeta^*\mathcal B$, and replaces $\mathcal A,\mathcal A'$ by
$Z^t\mathcal A,Z^t\mathcal A'$.  In~particular,
$Z\mathcal A'=\mathcal A$,
$Z^{-1}\mathcal A=\mathcal A'$.

\begin{corollary}\label{pw:available-prefix}
Let $t\in\mathbb Z$ and let $\mathcal B,\mathcal D\subset D^b(N)$ be
admissible. Suppose $Z^tL\zeta^*\mathcal B$ and
$Z^tL\zeta^*\mathcal D$ occur consecutively, in this order, in a
semiorthogonal decomposition of $D^b(M)$.
\begin{enumerate}
\item If an adjacent block $\mathcal C$ follows them and
$R\zeta_*(Z^{-t}\mathcal C)\subset\mathcal B$, then
\[
 \mathbf L_{Z^tL\zeta^*\langle\mathcal B,\mathcal D\rangle}\mathcal C
 =\mathbf L_{Z^tL\zeta^*\mathcal B}\mathcal C
 \subset Z^t\mathcal A.
\]
\item If an adjacent block $\mathcal C$ precedes them and
$R\zeta_*(Z^{1-t}\mathcal C)\subset\mathcal B$, then
\[
 \mathbf R_{Z^tL\zeta^*\langle\mathcal B,\mathcal D\rangle}\mathcal C
 =\mathbf R_{Z^tL\zeta^*\mathcal B}\mathcal C
 \subset Z^t\mathcal A'.
\]
\end{enumerate}
In both cases, the blocks in the family retain their embeddings and order.
\end{corollary}

\begin{proof}
Tensoring by $Z^{-t}$ reduces the proof to $t=0$.
In the first case, adjunction gives
$
 \operatorname{RHom}(L\zeta^*\mathcal D,\mathcal C)
 \simeq\operatorname{RHom}(\mathcal D,R\zeta_*\mathcal C)=0$.
The reverse vanishing follows from semiorthogonality. Thus $\mathcal C$
commutes unchanged past $L\zeta^*\mathcal D$, after which
\eqref{pw:left-criterion} applies to its left mutation through
$L\zeta^*\mathcal B$.

In the second case, put
$\mathcal C'=\mathbf R_{L\zeta^*\mathcal B}\mathcal C$.
By \eqref{pw:right-criterion}, $\mathcal C'\subset\mathcal A'$, so
adjunction gives
$\operatorname{RHom}(\mathcal C',L\zeta^*\mathcal D)=0$.
The reverse vanishing follows from the semiorthogonal order after this
mutation. Hence $\mathcal C'$ commutes unchanged past
$L\zeta^*\mathcal D$, proving the second equality and containment.
\end{proof}
\begin{proof}[Proof of Theorem~\ref{MainTheorem}]
Start with the  decomposition of Theorem~\ref{LoomTheorem}
ordered by decreasing $\lambda$ and then by decreasing $k$.  
Split the megablock II at
$\lambda=g-2$ into IIa and IIb and the megablock III at $\lambda=g-2$ into IIIa and IIIb: the mega-blocks
IIa and IIIa have $\lambda\geq g-2$, while IIb and IIIb have
$\lambda\leq g-3$.  We process the mega-blocks in the order
IV, I, IIa, IIb, IIIb, IIIa.
For $g=2$, the parts IIb and IIIb are empty and their steps below are omitted.

In every leftward procedure take the leftmost block not yet processed;
in every rightward procedure take the rightmost such block.
With this order, each block to which Corollary~\ref{pw:available-prefix}
is applied is adjacent to the retained family. The embeddings and
orders of the retained blocks remain unchanged.

\smallskip
\noindent\textit{Mega-block IV.}
Its blocks are
\begin{equation}\label{sRGSRGDRH}
 \left\langle Z^{\lambda+1-g}\theta^{m+1-\lambda}
       \cE^{\boxtimes(\lambda-2m)}\right\rangle_{{g-1\leq\lambda\leq2(g-1)}\atop
       {0\leq m\leq\lfloor\lambda/2\rfloor,\ \lambda-m\leq g-1}}.
\end{equation}
Retain the  family of rightmost blocks with $\lambda=g-1$:
\begin{equation}\label{dthehwe5jwjw6}
 \left\langle\theta^{2-g+k}\cE^{\boxtimes(g-1-2k)}
 \right\rangle_{k=\lfloor(g-1)/2\rfloor,\ldots,0}
       =L\zeta^*\mathcal B_{g-1}(2-g).
\end{equation}

\begin{lemma}\label{megablockiv}
Let $g\leq\lambda\leq2(g-1)$,
$0\leq m\leq\lfloor\lambda/2\rfloor$, and $\lambda-m\leq g-1$.
Put $k_0=g-1-\lambda+m$ and let $\mathcal B\subset D^b(M)$ be the
initial segment of \eqref{dthehwe5jwjw6} with
$
 k=\lfloor(g-1)/2\rfloor,\ldots,k_0$.
Then
\[
 \mathbf R_{\mathcal B}
 \left\langle Z^{\lambda+1-g}\theta^{m+1-\lambda}
       \cE^{\boxtimes(\lambda-2m)}\right\rangle
       \subset\mathcal A'.
\]
\end{lemma}

\begin{proof}
Set $n=\lambda-2m$, $s=\lambda+1-g$, and $c=m+1-\lambda$.
Then
\[
 n+s=2\lambda-2m+1-g=g-1-2k_0,\qquad c=2-g+k_0.
\]
Here $s\geq1$, and $0\leq n+s\leq g-1$ follows from $n\geq0$
and $\lambda-m\leq g-1$.  Thus $0\leq k_0\leq\lfloor(g-1)/2\rfloor$.
The assertion is \eqref{pw:right-mutation}, and \eqref{pw:prefix}
identifies its mutation category with the stated existing blocks.
\end{proof}

Process the blocks with $\lambda\geq g$ using this lemma and
Corollary~\ref{pw:available-prefix}
(b).  Each is right-mutated through
the indicated initial segment and then commuted past the rest of
\eqref{dthehwe5jwjw6}.  At the end, the family
\eqref{dthehwe5jwjw6} is at the left of mega-block IV, unchanged,
and the other blocks lie to its right in $\mathcal A'$.
For example, in genus $5$ the block with $(\lambda,m)=(6,2)$ is
mutated through all three blocks
$
 \langle\theta^{-1},\theta^{-2}\cE^{\boxtimes2},
                \theta^{-3}\cE^{\boxtimes4}\rangle$.

\smallskip
\noindent\textit{Mega-block I.}
Its blocks are
\begin{equation}\label{sEGsrhSRhjdt}
 \left\langle Z^{\lambda+2-g}\theta^{m-\lambda-1}
           \cE^{\boxtimes(\lambda-2m)}\right\rangle_{{0\leq\lambda\leq g-2}\atop
           {0\leq m\leq\lfloor\lambda/2\rfloor}}.
\end{equation}
Retain the family of leftmost blocks with $\lambda=g-2$:
\begin{equation}\label{sEGsrhSzxfafdbRhjdt}
 \left\langle\theta^{1-g+k}\cE^{\boxtimes(g-2-2k)}
 \right\rangle_{k=\lfloor(g-2)/2\rfloor,\ldots,0}
       =L\zeta^*\mathcal B_{g-2}(1-g).
\end{equation}

\begin{lemma}\label{megablocki}
For $0\leq\lambda\leq g-3$ and
$0\leq m\leq\lfloor\lambda/2\rfloor$, let $\mathcal B$ be the
initial segment of \eqref{sEGsrhSzxfafdbRhjdt} indexed by $k\geq m$.
Then
\[
 \mathbf L_{\mathcal B}
 \left\langle Z^{\lambda+2-g}\theta^{m-\lambda-1}
              \cE^{\boxtimes(\lambda-2m)}\right\rangle
       \subset\mathcal A.
\]
\end{lemma}

\begin{proof}
Set $n=\lambda-2m$, $s=\lambda+2-g<0$, and $c=m-\lambda-1$.
Then $n-s=g-2-2m$ and $c+s=m+1-g$.
Equations~\eqref{pw:left-mutation} and \eqref{pw:prefix} give the
claim.  The bound $0\leq m\leq\lfloor(g-2)/2\rfloor$ follows
from $2m\leq\lambda\leq g-3$.
\end{proof}

Apply 
Corollary~\ref{pw:available-prefix}
(a).  Before each left mutation,
the unprocessed block commutes past the part of
\eqref{sEGsrhSzxfafdbRhjdt} indexed by $k<m$; no retained block
is changed.  Thus the nonretained part of I moves to the left into
$\mathcal A$, and \eqref{sEGsrhSzxfafdbRhjdt} remains at its right.

\smallskip
\noindent\textit{Mega-block IIa.}
This part consists of
\begin{equation}\label{sBGASFNADNM}
 \left\langle Z^{\lambda+3-g}\theta^{m-\lambda-1}
              \cE^{\boxtimes(\lambda-2m)}\right\rangle_{{g-2\leq\lambda\leq2(g-2)}\atop
        {0\leq m\leq\lfloor\lambda/2\rfloor,\ \lambda-m\leq g-2}}.
\end{equation}
Initially retain the  family of blocks
\begin{equation}\label{AWRHQERHERJQEJ}
 \left\langle Z\theta^{1-g+k}\cE^{\boxtimes(g-2-2k)}
 \right\rangle_{k=\lfloor(g-2)/2\rfloor,\ldots,0}
       =ZL\zeta^*\mathcal B_{g-2}(1-g).
\end{equation}
For a block with $\lambda\geq g-1$, put
$k_0=g-2-\lambda+m$.  Divide the block and the retained family by
$Z$ and apply \eqref{pw:right-mutation}.  Indeed,
\[
 \begin{gathered}
 n=\lambda-2m,\quad s=\lambda+2-g,\quad c=m-\lambda-1,\\
 n+s=g-2-2k_0,\qquad c=1-g+k_0.
 \end{gathered}
\]
Thus right mutation through the initial segment $k\geq k_0$ of
\eqref{AWRHQERHERJQEJ} places this block in
$Z\mathcal A'=\mathcal A$.  The cutoff is in the displayed range:
$k_0\geq0$ follows from $\lambda-m\leq g-2$, while
$g-2-2k_0=n+s\geq0$.  
Corollary~\ref{pw:available-prefix} (b), with
$t=1$, moves each processed block after the whole retained family.

We now have the two consecutive families
\eqref{sEGsrhSzxfafdbRhjdt}, \eqref{AWRHQERHERJQEJ}.
Their boundary mutation is  given by
\begin{equation}\label{pw:positive-boundary-sequence}
 0\longrightarrow\cO_M\longrightarrow Z
       \longrightarrow\cO_Z(Z)\longrightarrow0.
\end{equation}
The first map becomes an isomorphism after $R\zeta_*$:
by Lemma~\ref{pw:direct-images} both direct images are $\cO_N$,
and the map is nonzero on the open set where $\zeta$ is an
isomorphism.  Since $H^0(N,\cO_N)=\mathbb C$, it is an isomorphism
everywhere.  Consequently $R\zeta_*\cO_Z(Z)=0$.
For $F\in\mathcal B_{g-2}(1-g)$, the counit therefore identifies
\begin{equation}\label{pw:positive-boundary-mutation}
 \mathbf L_{L\zeta^*\mathcal B_{g-2}(1-g)}(ZL\zeta^*F)
       \simeq\cO_Z(Z)\otimes L\zeta^*F\in\mathcal A.
\end{equation}
This mutates the whole family \eqref{AWRHQERHERJQEJ} to the left
of \eqref{sEGsrhSzxfafdbRhjdt}.  The already processed higher IIa
blocks are also in $\mathcal A$ and can now be commuted to that
side: one orthogonality follows from adjunction and the other from
the current semiorthogonal order.  Hence all nonretained blocks of
I and IIa are to the left of \eqref{sEGsrhSzxfafdbRhjdt}, in
$\mathcal A$.

\smallskip
\noindent\textit{Mega-block IIb.}
Its kernels have the same formula as \eqref{sBGASFNADNM}, with
$0\leq\lambda\leq g-3$.  Retain the  family
\begin{equation}\label{kjsrHwrhRJegfw}
 \left\langle\theta^{2-g+k}\cE^{\boxtimes(g-3-2k)}
 \right\rangle_{k=\lfloor(g-3)/2\rfloor,\ldots,0}
       =L\zeta^*\mathcal B_{g-3}(2-g).
\end{equation}
For $\lambda\leq g-4$, apply \eqref{pw:left-mutation} with
\[
 \begin{gathered}
 n=\lambda-2m,\quad s=\lambda+3-g,\quad c=m-\lambda-1,\\
 n-s=g-3-2m,\qquad c+s=m+2-g.
 \end{gathered}
\]
The mutation category is the initial segment $k\geq m$ of
\eqref{kjsrHwrhRJegfw}.  
Corollary~\ref{pw:available-prefix} (a)
places the other IIb blocks on its left in $\mathcal A$.
They then commute past \eqref{sEGsrhSzxfafdbRhjdt}, by adjunction
and the current semiorthogonal order.  The two retained families
from I and II are unchanged and consecutive.

\smallskip
\noindent\textit{Mega-block IIIb.}
This part has kernels
\[
 Z^{\lambda+2-g}\theta^{m-\lambda}
       \cE^{\boxtimes(\lambda-2m)},\qquad
 0\leq\lambda\leq g-3,\quad0\leq m\leq\lfloor\lambda/2\rfloor.
\]
Retain the family with $\lambda=g-3$:
\begin{equation}\label{srgrhasrharhqer}
 \left\langle Z^{-1}\theta^{3-g+k}\cE^{\boxtimes(g-3-2k)}
 \right\rangle_{k=\lfloor(g-3)/2\rfloor,\ldots,0}
       =Z^{-1}L\zeta^*\mathcal B_{g-3}(3-g).
\end{equation}
For $\lambda\leq g-4$, multiply both the block and this family by
$Z$ and apply \eqref{pw:left-mutation}, with
\[
 \begin{gathered}
 n=\lambda-2m,\quad s=\lambda+3-g,\quad c=m-\lambda,\\
 n-s=g-3-2m,\qquad c+s=m+3-g.
 \end{gathered}
\]
The required initial segment in \eqref{srgrhasrharhqer} is $k\geq m$.
After tensoring back by $Z^{-1}$, the processed blocks lie to its
left in $Z^{-1}\mathcal A=\mathcal A'$.

The family \eqref{srgrhasrharhqer} is $Z^{-1}$ times the initial
segment $k\geq1$ of \eqref{dthehwe5jwjw6}.  Denote the corresponding
subcategory on $N$ by $\mathcal B_{g-3}(3-g)$ as above.  The second
boundary mutation is  induced by
\begin{equation}\label{pw:negative-boundary-sequence}
 0\longrightarrow Z^{-1}\longrightarrow\cO_M
       \longrightarrow\cO_Z\longrightarrow0.
\end{equation}
For $F\in\mathcal B_{g-3}(3-g)$, the unit of the left adjunction
gives
\begin{equation}\label{pw:negative-boundary-mutation}
 \mathbf R_{L\zeta^*\mathcal B_{g-3}(3-g)}(Z^{-1}L\zeta^*F)
       \simeq\cO_Z\otimes L\zeta^*F[-1]\in\mathcal A'.
\end{equation}
Here $R\zeta_*(Z\cO_Z)=0$ by
\eqref{pw:positive-boundary-sequence}.  Apply this right mutation
through the actual initial segment $k\geq1$, and commute the result
past the remaining $k=0$ block as in
Corollary~\ref{pw:available-prefix} (b).  After all these boundary blocks
have been processed, the lower IIIb blocks already in $\mathcal A'$
commute past the retained family \eqref{dthehwe5jwjw6}.  All
nonretained blocks from IIIb and IV lie to its right in $\mathcal A'$.

\smallskip
\noindent\textit{Mega-block IIIa.}
Its kernels are
\[
 \begin{gathered}
 Z^{\lambda+2-g}\theta^{m-\lambda}
       \cE^{\boxtimes(\lambda-2m)},\\
 g-2\leq\lambda\leq2(g-2),\quad
 0\leq m\leq\lfloor\lambda/2\rfloor,\quad\lambda-m\leq g-2.
 \end{gathered}
\]
Retain the  family
\begin{equation}\label{dthasgarharhwjw6}
 \left\langle\theta^{2-g+k}\cE^{\boxtimes(g-2-2k)}
 \right\rangle_{k=\lfloor(g-2)/2\rfloor,\ldots,0}
       =L\zeta^*\mathcal B_{g-2}(2-g).
\end{equation}
For $\lambda\geq g-1$, put $k_0=g-2-\lambda+m$ and use
\eqref{pw:right-mutation}, since
\[
 \begin{gathered}
 n=\lambda-2m,\quad s=\lambda+2-g,\quad c=m-\lambda,\\
 n+s=g-2-2k_0,\qquad c=2-g+k_0.
 \end{gathered}
\]
The initial segment $k\geq k_0$ of
\eqref{dthasgarharhwjw6} is available, by the same bounds as for IIa.
The resulting blocks lie in $\mathcal A'$ after the retained family.
They commute also past \eqref{dthehwe5jwjw6}, by adjunction and the
current semiorthogonal order.  

The decomposition now has the four retained families
\eqref{sEGsrhSzxfafdbRhjdt}, \eqref{kjsrHwrhRJegfw},
\eqref{dthasgarharhwjw6}, and \eqref{dthehwe5jwjw6}, in this order,
consecutive and unchanged.  All other blocks from I and II are on
their left in $\mathcal A$, and all other blocks from III and IV
are on their right in $\mathcal A'$.
Move the terminal group of nonretained blocks to the front by the Serre functor
$(-)\otimes\omega_M[\dim M]$.  It then lies in $\mathcal A$.
Thus we have obtained, by adjacent mutations and the specified
orthogonal commutations, a full semiorthogonal decomposition
\[
 \begin{aligned}
 D^b(M)=\bigl\langle &\mathcal K,
 L\zeta^*\mathcal B_{g-2}(1-g),
 L\zeta^*\mathcal B_{g-3}(2-g),\\
 &L\zeta^*\mathcal B_{g-2}(2-g),
 L\zeta^*\mathcal B_{g-1}(2-g)\bigr\rangle,
 \qquad \mathcal K\subset\mathcal A.
 \end{aligned}
\]
Here the second retained family is omitted when $g=2$.
For $F\in D^b(N)$, apply $R\zeta_*$ to the filtration of
$L\zeta^*F$ associated with this decomposition.  The $\mathcal K$
term vanishes, and $R\zeta_*L\zeta^*F=F$.
Hence the retained families generate $D^b(N)$.
Their full faithfulness and semiorthogonality were already available
from Theorem~\ref{LoomTheorem} and full faithfulness of $L\zeta^*$.
This proves Theorem~\ref{MainTheorem}.
\end{proof}

We finish this section by introducing a different semi-orthogonal decomposition of $D^b(N)$,
which has blocks given by tensor bundles $\bar\cE^{\boxtimes p}$.

\begin{figure}[htbp]
\includegraphics[width=0.8\textwidth]{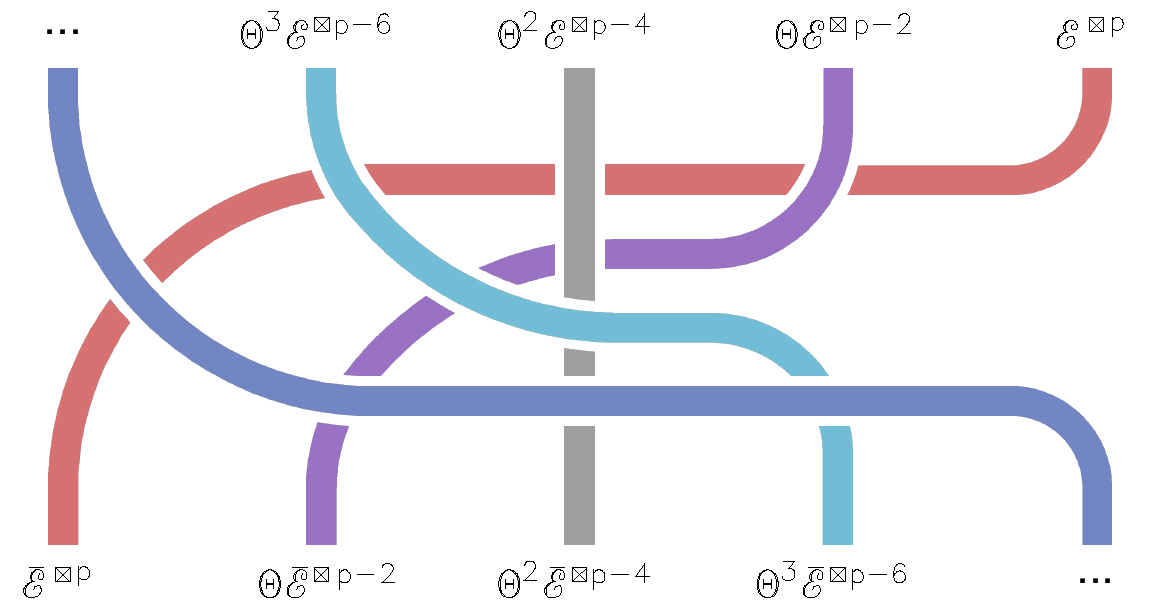}
\caption{The  mutation from ordinary to sign blocks.}\label{asfbadfbadha}
\end{figure}

\begin{theorem}\label{ComparingSODs}
Each of the four mega-blocks of Theorem~\ref{MainTheorem} can  alternatively
(and independently of the other mega-blocks) be
decomposed 
as follows:
$$
\Bigl\langle
\bigl\langle 
{\theta}^{1-g+k}\otimes\bar\cE^{\boxtimes g-2-2k}
\bigr\rangle_{0\le k\le \lfloor{g-2\over2}\rfloor},\quad
\bigl\langle 
{\theta}^{2-g+k}\otimes\bar\cE^{\boxtimes g-3-2k}
\bigr\rangle_{0\le k\le  \lfloor{g-3\over2}\rfloor},\quad{}$$
$${}\quad\bigl\langle 
{\theta}^{2-g+k}\otimes\bar\cE^{\boxtimes g-2-2k}
\bigr\rangle_{0\le k\le \lfloor{g-2\over2}\rfloor},
\quad \bigl\langle 
{\theta}^{2-g+k}\otimes\bar\cE^{\boxtimes g-1-2k}
\bigr\rangle_{0\le k\le  \lfloor{g-1\over2}\rfloor}
\Bigr\rangle.
$$ 
Within the mega-blocks with $\bar\cE$, the 
blocks are arranged in increasing order of~$k$ (instead of decreasing order of $k$ 
for the mega-blocks with $\cE$ as in Theorem~\ref{MainTheorem}). Up to a twist by a tensor power of $\theta$,
the two decompositions of each mega-block in Theorems~\ref{MainTheorem} and \ref{ComparingSODs} are related by the mutation in Figure~\ref{asfbadfbadha}.
Here the first block in the top row (and the last block in the bottom row) is given by 
$\theta^{l}$ if $p=2l$ and by $\theta^l\otimes\cE=\theta^l\otimes\bar\cE$ if $p=2l+1$.
\end{theorem}

\begin{proof}
Derived duality reverses a semiorthogonal sequence.  By
Lemma~\ref{amazinggrace} and relative duality for the source of a
Fourier--Mukai functor, it sends
$\langle\cE^{\boxtimes d}\rangle$ to
$\theta^{-d}\langle\bar\cE^{\boxtimes d}\rangle$, and sends the
latter sign block back to $\theta^{-d}\langle\cE^{\boxtimes d}\rangle$.
The source line bundles and cohomological shifts do not change these
essential images. Each mega-block in the statement is a common twist
$\theta^c\mathcal B_h(0)$, so it suffices to work with
$\mathcal B_h(0)$. For $0\le h\le g-1$ we claim
\begin{equation}\label{pw:ordinary-sign-family}
 \mathcal B_h(0)=
 \left\langle\theta^j\bar\cE^{\boxtimes(h-2j)}
 \right\rangle_{j=0,\ldots,\lfloor h/2\rfloor}.
\end{equation}
Proposition~\ref{pw:tensor-products}, with $a=0$, gives
$
 \langle\bar\cE^{\boxtimes(h-2j)}\rangle
       \subset\mathcal B_{h-2j}(0)
$.
After tensoring by $\theta^j$, \eqref{pw:prefix} places the right
side in $\mathcal B_h(0)$.  Thus the right side of
\eqref{pw:ordinary-sign-family} is contained in the left side.
Apply derived duality to this inclusion and then tensor by
$\theta^h$.  The duality formulas above interchange the two families,
so the resulting inclusion is the reverse one.  This proves equality.
The barred order in \eqref{pw:ordinary-sign-family} is
semiorthogonal, since it is the dual of the ordinary order followed
by the common twist $\theta^h$.

We prove the asserted mutation by induction on $h$.
For $h=0,1$, ordinary and sign tensor bundles agree and there is
nothing to mutate.  For $h\geq2$, the ordinary decomposition is
$
 \mathcal B_h(0)=
 \langle\theta\mathcal B_{h-2}(0),\cE^{\boxtimes h}\rangle
$.
Apply the inductively known mutation inside
$\theta\mathcal B_{h-2}(0)$, without changing that subcategory,
and left-mutate the last block through it, giving
$
 \left\langle
 \mathbf L_{\theta\mathcal B_{h-2}(0)}
       (\langle\cE^{\boxtimes h}\rangle),
 \theta\mathcal B_{h-2}(0)\right\rangle$.
On the other hand, \eqref{pw:ordinary-sign-family} is
$
 \mathcal B_h(0)=
       \langle\bar\cE^{\boxtimes h},\theta\mathcal B_{h-2}(0)\rangle$.
The first components in these decompositions are both
$\mathcal B_h(0)\cap(\theta\mathcal B_{h-2}(0))^\perp$.
They are therefore equal.  This is the remaining mutation in
Figure~\ref{asfbadfbadha}.  Twisting by $\theta^c$ applies the result
independently to each of the four mega-blocks.
\end{proof}

\end{document}